\documentclass[nosumlimits,twoside]{amsart}
\usepackage[cp850]{inputenc}
\usepackage[bookmarksnumbered,plainpages,hyperindex]{hyperref}
\input{diagrams.tex}
\diagramstyle[Postscript=dvips]
\newarrow{Eqto}{=}{=}{=}{=}{=}
\newarrow{Line}{-}{-}{-}{-}{-}
\newarrow{Dashto}{}{dash}{}{dash}>
\newarrow{Dotsto}{.}{.}{.}{.}>
\newarrow{eqto}{=}{=}{=}{=}{=}
\newtheorem{theorem}{\sc Theorem}[section]
\newtheorem{proposition}[theorem]{\sc Proposition}
\newtheorem{notation}[theorem]{\sc Notation}
\newtheorem{lemma}[theorem]{\sc Lemma}
\newtheorem{corollary}[theorem]{\sc Corollary}
\theoremstyle{definition}
\newtheorem{definition}[theorem]{\sc Definition}

\newtheorem{example}[theorem]{\sc Example}

\theoremstyle{remark}
\newtheorem{remark}[theorem]{\sc Remark}

\newtheorem{claim}[theorem]{}

\def\ker{\mathrm{Ker}}
\def\K{\mathrm{Ker}}
\def\C{\mathrm{Coker}}

\def\M{\mathcal{M}}

\def\H{\mathbb{H}}

\def\ot{\otimes}
\def\cot{\square}
\def\N{\mathbb{N}}
\def\Mt{{^C\mathcal{M}^C}}

\def\dprod{\prod}
\def\dbigoplus{\bigoplus}
\begin{document}
\title{Cotensor Coalgebras in Monoidal Categories}
\author{A. Ardizzoni, C. Menini}
\author{D. \c{S}tefan}
\thanks{This paper was written while A. Ardizzoni and C.
Menini were members of G.N.S.A.G.A. with partial financial support
from M.I.U.R..}

\begin{abstract}
We introduce the concept of cotensor coalgebra for a given
bicomodule over a coalgebra in an abelian monoidal category. Under
some further conditions we show that such a cotensor coalgebra
exists and satisfies a meaningful universal property. We prove
that this coalgebra is formally smooth whenever the comodule is
relative injective and the coalgebra itself is formally smooth.
\end{abstract}

\maketitle

\pagestyle{headings}

\section*{Introduction}

\markboth{\sc{A. ARDIZZONI, C. MENINI AND D. \c{S}TEFAN}}
{\sc{Cotensor Coalgebras in Monoidal Categories}}

Let $C$ be a coalgebra over a field $k$ and let $M$ be a
$C$-bicomodule. The cotensor coalgebra $ T^c_C(M)$ was introduced
by Nichols in \cite{Ni} as a main tool to construct some new Hopf
algebras that he called "bialgebras of type one". These bialgebras
can be reconstructed, via a bosonization procedure, from the so
called Nichols algebras, which are essentially the $H$-coinvariant
elements of the bialgebras of type one, in the case when $C=H$ is
a Hopf algebra and $M$ is a Hopf bimodule. Nichols algebras, also
named quantum symmetric algebras in \cite{Ro}, have been deeply
investigated and appear as a main step in the classification of
finite dimensional Hopf algebras problem (see, e.g., \cite{AG} and
\cite{AS}). In fact, in the case that $C=H$ is a Hopf algebra and
$M$ is a Hopf $H$-bimodule, the cotensor coalgebra $ T^c_C(M)$ is
a bialgebra that is called "quantum shuffle Hopf algebra" by Rosso
in \cite{Ro} where some fundamental properties of this bialgebra
and of its coinvariant Hopf algebra are investigated. The
coalgebra of paths of a quiver $Q$ is an instance of a cotensor
coalgebra. Namely let $Q_0$ be the set of vertices and let $Q_1$
be the set of arrows of $Q$. Then $M=KQ_1$ is a $C$-bicomodule
where $C=KQ_0$ is equipped with its natural coalgebra structure.
The cotensor coalgebra $ T^c_C(M)$ is the path coalgebra of the
quiver $Q$. In \cite{CR}, Cibils and Rosso provide the
classification of path coalgebras which admit a graded Hopf
algebra structure, allowing the quiver to be infinite. On the
other hand, in \cite{JLMS}, hereditary coalgebras with coseparable
coradical are characterized by means of a suitable cotensor
coalgebra. Moreover it is proved that if $C$ is a formally smooth
coalgebra (see Definition \ref{def: formally smooth}) and $M$ is
$\mathcal{I}$-injective (see \ref{cl:E-Proj}) then $ T^c_C(M)$ is
formally smooth.\newline In this paper we introduce the notion of
cotensor coalgebra in an abelian monoidal category. We would like
to outline that this fact is not immediate. In fact the notion of
coradical plays a fundamental role in the usual definition for
coalgebras over a field (see \cite{Ni}) while we have no coradical
substitution here. Also, having developed in \cite{AMS} the notion
of formally smooth coalgebras for abelian monoidal categories, we
wanted to obtain the second quoted result of \cite{JLMS} in this
more general setting, namely we prove Theorem \ref{teo: smooth}.
In a forthcoming paper we will deal with the other quoted result,
establishing a suitable criterion.\medskip\newline The paper is
organized as follows.\newline In Section 1 we give the definition
of abelian monoidal category (see Definition
\ref{abelianmonoidal}) and prove some preliminary results on
direct limits of direct systems in monoidal categories.\newline
Section 2 is devoted to the construction of the cotensor
coalgebra. Its coalgebra structure is obtained (see Theorem
\ref{teo: T as limit}) as a direct limit of a direct system of
certain coalgebras defined by means of Hochschild $2$-cocycles. In
Theorem \ref{coro: univ property of cotensor coalgebra} we we show
that such a cotensor coalgebra exists and satisfies a meaningful
universal property.\newline Section 3 deals with some technical
results that we will use in Section 4 where we treat formal
smoothness of the cotensor coalgebra (see Theorem \ref{teo
Smooth}).\newline Finally, in Section 5, we relate our results
with the classical one of vector spaces and prove a criterion
(Theorem \ref{teo: monoidal functor}) which establishes a strong
universal property of the cotensor coalgebra for the main abelian
monoidal categories related to a Hopf algebra.\medskip\newline
\textbf{Notation.} \ In a category $\mathcal{M}$ the set of
morphisms from $X$ to $Y$ will be denoted by $\mathcal{M}(X,Y).$
If $X$ is an object in $\mathcal{M}$ then the functor
$\mathcal{M}(X,-)$ from $\mathcal{M}$ to $\mathfrak{Sets}$
associates to any morphism $u:U\rightarrow V$ in $\mathcal{M}$ the
function that will be denoted by $\mathcal{M}(X,u).$

\section{Direct Limits in Monoidal Categories}

\begin{claim}
\label{MonCat}A \emph{monoidal category} means a category
$\mathcal{M}$ that is endowed with a functor $\otimes
:\mathcal{M}\times \mathcal{M}\rightarrow \mathcal{M}$, an object
$\mathbf{1}\in \mathcal{M}$ and functorial isomorphisms:
$a_{X,Y,Z}:(X\otimes Y)\otimes Z\rightarrow X\otimes (Y\otimes
Z),$ $l_{X}:\mathbf{1}\otimes X\rightarrow X$ and $r_{X}:X\otimes \mathbf{1}%
\rightarrow X.$ The functorial morphism $a$ is called the \emph{%
associativity constraint }and\emph{\ }satisfies the \emph{Pentagon Axiom, }%
that is the following diagram
\begin{equation*}
\begin{diagram}[grid=pentagon,size=2em]
&&&&&&(U\otimes V)\otimes (W\otimes X)&\\
&&&&&\relax\ruTo(6,2)^{a_{U\ot V,W,X}}&&&\relax\rdTo(6,2)^{a_{U,V,W\otimes X}}&\\
((U\otimes V)\otimes W)\otimes X&&&&&&&&&&&&U\otimes (V\otimes (W\otimes X)&\\
&\relax\rdTo_{a_{U,V,W}\otimes X}&&&&&&&&&&\relax\ruTo_{U\otimes a_{V,W,X}}&\\
&&(U\otimes (V\otimes W))\otimes X&&&&\relax\rTo^{a_{U,V\otimes
W,X}}&&&&U\otimes ((V\otimes W)\otimes X)&
\end{diagram}
\end{equation*}
is commutative, for every $U,\,V,$ $W,$ $X$ in $\mathcal{M}.$ The morphisms $%
l$ and $r$ are called the \emph{unit constraints} and they are
assumed to satisfy the \emph{Triangle Axiom, }i.e. the following
diagram
\begin{equation*}
\begin{diagram}
(V\otimes \mathbf{1})\otimes W && \rTo_{a_{V,\mathbf{1},W}} && V\otimes(\mathbf{1}\otimes W)\\
&\rdTo<{r_V\otimes W}&&\ldTo>{V\otimes l_W} \\ && V\otimes W&
\end{diagram}
\end{equation*}
is commutative. The object $\mathbf{1}$ is called the \emph{unit} of $%
\mathcal{M}$.For details on monoidal categories we refer to \cite[Chapter XI]%
{Ka} and \cite{Maj2}. A monoidal category is called \emph{strict}
if the associativity constraint and unit constraints are the
corresponding identity morphisms.
\end{claim}

\begin{claim}
\label{cl:CohThm}As it is noticed in \cite[p. 420]{Maj2}, the
Pentagon Axiom solves the consistency problem that appears because
there are two ways to go from $((U\otimes V)\otimes W)\otimes X$
to $U\otimes (V\otimes (W\otimes X)).$ The coherence theorem, due
to S. Mac Lane, solves the similar problem for the tensor product
of an arbitrary number of objects in $\mathcal{M}.$ Accordingly
with this theorem, we can always omit all brackets and simply
write $X_{1}\otimes \cdots \otimes X_{n}$ for any object obtained from $%
X_{1},\ldots ,X_{n}$ by using $\otimes $ and brackets. Also as a
consequence of the coherence theorem, the morphisms $a,$ $l,$ $r$
take care of themselves, so they can be omitted in any computation
involving morphisms in $\mathcal{M.}$\medskip\ The notions of
algebra, module over an algebra, coalgebra and comodule over a
coalgebra can be introduced in the general setting of monoidal
categories. For more details, see \cite{AMS}.
\end{claim}

\begin{proposition}\label{pro: bicomodules}Let $\mathcal{M}$ be a monoidal category. Assume that $\M$ is also an abelian category.
Let $C$ an $D$ be coalgebras in $\M$. Let $X,Y\in
{^{C}\mathcal{M}^D}$ and let $f:X\rightarrow Y$ be a morphism of
$(C,D)$-bicomodules. Let
$\mathbb{H}:{^{C}\mathcal{M}^D}\rightarrow \mathcal{M}$ be the
forgetful functor. Then

1) $\C (\mathbb{H}(f))$ carries a natural $(C,D)$-bicomodule
structure (such that the definition map is a morphism of
bicomodules) that makes it the cokernel of $f$ in
$^{C}\mathcal{M}^D$.

2) if the functor $C\ot (-)\ot D:\M\to {\mathcal{M}}$ is left
exact (e.g. the tensor functors are left exact), then $\K(\H(f))$
carries a natural $(C,D)$-bicomodule structure (such that the
definition map is a morphism of bicomodules) that makes it the
kernel of $f$ in ${^{C}\mathcal{M}^D}$.
\end{proposition}

\begin{proof} Since the dual category of an abelian category is
an abelian category, the conclusion follows by applying
\cite[Proposition 3.3]{Ar}.
\end{proof}

The previous proposition justifies the following definition.

\begin{definition}\label{abelianmonoidal}A monoidal category
$(\M,\ot,\mathbf{1})$ will be called an \emph{\textbf{abelian
monoidal category} } if:
\begin{enumerate}
    \item $\M$ is an abelian category
    \item both the functors $X\ot (-):\M\to\M$ and $(-)\ot
    X:\M\to\M$ are additive and left exact, for every object $X\in
    \M$.
\end{enumerate}
\end{definition}

\begin{remark}
  The above definition of abelian monoidal category is dual to \cite[Definition
  1.8]{AMS}, where algebras were investigated. Therefore, we
  should have used the terminology "dual abelian monoidal
  category", which we had not for matter of style.
\end{remark}

\begin{claim}
Let $E$ be a coalgebra in an abelian monoidal category
$\mathcal{M}$. Let us recall,
(see \cite[page 60]{Mo}), the definition of wedge product of two subobjects $%
X,Y $ of $E$ in $\M:$%
\begin{equation*}
X\wedge_E Y:=Ker[ (p _{X}\otimes p _{Y}) \circ \triangle _{E}] ,
\end{equation*}
where $p _{X}:E\rightarrow E/X$ and $p _{Y}:E\rightarrow E/Y$ are
the canonical quotient maps.
\end{claim}

\begin{proposition}\label{pro: ker is a coalgebra}
Let $(\M ,\ot ,\mathbf{1})$ be an abelian monoidal category. Let
$(C,\Delta ,\varepsilon )$ be a coalgebra in $\M$ and let $L$ be a
$C$-bicomodule. Let $f:C\rightarrow L$ be a morphism in
$\Mt$, where $C$ is regarded as a bicomodule via $\Delta $. Then $(D,\delta):=%
\K(f)$ carries a natural coalgebra structure such that $\delta$ is
a morphism of coalgebras.
\end{proposition}

\begin{proof}By Proposition \ref{pro: bicomodules}, $D$ is a
$C$-bicomodule and $\delta$ is a morphism of bicomodules. Denote
by ${\rho_D^l}$ and ${\rho_D^r}$ the left and the right
$C$-comodule structure of $D$ respectively.\newline By left
exactness of the tensor functors, we have that $(D\otimes
D,D\otimes \delta)=\K(D\otimes f)$. Consider the following
diagram:
\begin{equation*}
\begin{diagram}[h=2em,w=3em] 0& \rTo& D\ot D&\rTo^{D\ot \delta}& D\ot C
&\rTo^{D\ot f}&D\ot L\\ & && \luDotsto<{\Delta_D}&
\uTo>{\rho^r_D}\\ &&&& D
\\ \end{diagram}
\end{equation*}
We have:
\begin{equation*}
(\delta\otimes L)(D\otimes f){\rho _{D}^{r}}=(C\otimes
f)(\delta\otimes C){\rho _{D}^{r}=}(C\otimes f)\Delta
\delta={{\rho _{L}^{l}}}f\delta=0.
\end{equation*}
By left exactness of the tensor functors, $\delta\otimes L$ is a
monomorphism, so that we get $(D\otimes f){\rho _{D}^{r}}=0$ and
hence, by the universal property of the kernel, there exists a
unique morphism $\Delta _{D}:D\rightarrow
D\otimes D$ in $\M$ such that $(D\otimes \delta)\Delta _{D}=\rho _{D}^{r}$.%
\newline
Let us prove that $\Delta _{D}$ is coassociative. Since
\begin{equation}
(\delta\otimes \delta)\Delta _{D}=(\delta\otimes C)(D\otimes
\delta)\Delta _{D}=(\delta\otimes C)\rho _{D}^{r}=\Delta \delta,
\label{delta is a coalg morph}
\end{equation}
we get:
\begin{eqnarray*}
(\delta\otimes \delta\otimes \delta)(\Delta _{D}\otimes D)\Delta
_{D} &=&(\Delta \otimes
C)(\delta\otimes \delta)\Delta _{D} \\
&=&(\Delta \otimes C)\Delta \delta \\
&=&(C\otimes \Delta )\Delta \delta \\
&=&(C\otimes \Delta )(\delta\otimes \delta)\Delta
_{D}=(\delta\otimes \delta\otimes \delta)(D\otimes \Delta
_{D})\Delta _{D}.
\end{eqnarray*}
By left exactness of the tensor functors, $\delta\otimes
\delta\otimes \delta=(C\otimes C\otimes \delta)(\delta\otimes
C\otimes D)(D\otimes \delta\otimes D)$ is a monomorphism so
that we obtain $(\Delta _{D}\otimes D)\Delta _{D}=(D\otimes \Delta _{D})\Delta _{D}.$%
\newline
Set $\varepsilon _{D}:=\varepsilon _{C}\delta:D\rightarrow
\mathbf{1.}$ Then we have:
\begin{equation*}
(\delta\otimes \mathbf{1})(D\otimes \varepsilon _{D})\Delta
_{D}=(C\otimes \varepsilon _{C})(\delta\otimes \delta)\Delta
_{D}=(C\otimes \varepsilon _{C})\Delta
\delta=r_{C}^{-1}\delta=(\delta\otimes \mathbf{1})r_{D}^{-1}.
\end{equation*}
Since $(\delta\otimes \mathbf{1})$ is a monomorphism, we have
$(D\otimes \varepsilon _{D})\Delta _{D}=r_{D}^{-1}.$ Analogously
one gets $(\varepsilon _{D}\otimes D)\Delta _{D}=l_{D}^{-1}.$ Thus
$(D,\Delta _{D},\varepsilon _{D})$ is a
coalgebra and by relation (\ref{delta is a coalg morph}) and definition of $%
\varepsilon _{D}$, $\delta$ is a homomorphism of coalgebras.
\end{proof}

\begin{proposition}\label{pro: colimit=coalg}
Let $\M$ be a monoidal category with direct limits. Let
$((X_{i})_{i\in \mathbb{N}},(\xi _{i}^j)_{i,j\in \mathbb{N}})$ be
a direct system in $\M$, where, for $i\leq j$, $\xi
_{i}^j:X_{i}\rightarrow X_{j}$. Assume that $X_{i}$ is a coalgebra
and that $\xi _{i}^j$
is a homomorphism of coalgebras for any $i,j\in \mathbb{N}.$ Then $%
\underrightarrow{\lim }X_{i}$ carries a natural coalgebra
structure that makes it the direct limit of $((X_{i})_{i\in
\mathbb{N}},(\xi _{i}^j)_{i,j\in \mathbb{N}})$ as a direct system
of coalgebras.
\end{proposition}

\begin{proof}Let $(X_{i},\Delta _{X_{i}},\varepsilon _{X_{i}})$ be a coalgebra in $(\M%
,\otimes ,\mathbf{1})$ for any $i\in \mathbb{N.}$ Set $X:=\underrightarrow{%
\lim }X_{i}.$ Let $(\xi _{i}:X_{i}\rightarrow X)_{i\in
\mathbb{N}}$ be the structural morphism of the direct limit, so
that $\xi _{j}\xi _{i}^j=\xi _{i}$ for any $i\leq j$. We put
\begin{equation*}
\Delta _{i}=(\xi _{i}\otimes \xi _{i})\Delta
_{X_{i}}:X_{i}\rightarrow X\otimes X,\text{ for any }i\in
\mathbb{N.}
\end{equation*}
Since $\xi _{i}^j$ is a homomorphism of coalgebras, we have that:
\begin{equation*}
\Delta _{j}\xi^j _{i}=(\xi _{j}\otimes \xi _{j})\Delta _{X_{j}}\xi
_{i}^j=(\xi _{j}\otimes \xi _{j})(\xi _{i}^j\otimes \xi
_{i}^j)\Delta _{X_{i}}=(\xi _{i}\otimes \xi _{i})\Delta
_{X_{i}}=\Delta _{i},
\end{equation*}
so that there exists a unique morphism $\Delta :X\rightarrow
X\otimes X$ such that
\begin{equation}
\Delta \xi _{i}=\Delta _{i}=(\xi _{i}\otimes \xi _{i})\Delta
_{X_{i}}\text{ for any }i\in \mathbb{N.}  \label{xi coalg1}
\end{equation}
Let us prove that $\Delta $ is coassociative:
\begin{eqnarray*}
(X\otimes \Delta )\Delta \xi _{i} &=&(X\otimes \Delta )(\xi
_{i}\otimes \xi
_{i})\Delta _{X_{i}} \\
&=&(\xi _{i}\otimes \xi _{i}\otimes \xi _{i})(X_{i}\otimes \Delta
_{X_{i}})\Delta _{X_{i}} \\
&=&(\xi _{i}\otimes \xi _{i}\otimes \xi _{i})(\Delta
_{X_{i}}\otimes
X_{i})\Delta _{X_{i}} \\
&=&(\Delta \otimes X)(\xi _{i}\otimes \xi _{i})\Delta
_{X_{i}}=(\Delta \otimes X)\Delta \xi _{i}.
\end{eqnarray*}
Since the relation above holds true for any $i\in \mathbb{N,}$ by
the universal property of the direct limit we deduce that
$(X\otimes \Delta )\Delta =(\Delta \otimes X)\Delta .$ Now, as
$\xi _{i}^j$ is a homomorphism of coalgebras, $\varepsilon
_{X_{j}}\xi _{i}^j=\varepsilon _{X_{i}}.$ Hence, there exists a
unique morphism $\varepsilon :X\rightarrow \mathbf{1}$ such that
\begin{equation}
\varepsilon \xi _{i}=\varepsilon _{X_{i}}\text{ for any }i\in
\mathbb{N.} \label{xi coalg2}
\end{equation}
Then we have:
\begin{equation*}
(X\otimes \varepsilon )\Delta \xi _{i}=(X\otimes \varepsilon )(\xi
_{i}\otimes \xi _{i})\Delta _{X_{i}}=(\xi _{i}\otimes \mathbf{1}%
)(X_{i}\otimes \varepsilon _{X_{i}})\Delta _{X_{i}}=(\xi _{i}\otimes \mathbf{%
1})r_{X_{i}}^{-1}=r_{X}^{-1}\xi _{i}.
\end{equation*}
Since the relation above holds true for any $i\in \mathbb{N,}$ by
the universal property of direct limits we deduce that $(X\otimes
\varepsilon )\Delta =r_{X}^{-1}.$ Analogously one gets
$(\varepsilon \otimes X)\Delta =l_{X}^{-1} $. Thus $(X,\Delta
,\varepsilon )$ is a coalgebra in $\M.$ Note that relations
(\ref{xi coalg1}) and (\ref{xi coalg2}) mean that $\xi
_{i}:X_{i}\rightarrow X$ is a homomorphism of coalgebras. \newline
Let now $(C,\Delta _{C},\varepsilon _{C})$ be a coalgebra in $\M$ and let $%
(f_{i}:X_{i}\rightarrow C)_{i\in \mathbb{N}}$ be a compatible
family of morphisms of coalgebras in $\M$. Since
$(f_{i}:X_{i}\rightarrow C)_{i\in \mathbb{N}}$ is a compatible
family of morphisms in $\M,$ there exists a unique morphism
$f:X\rightarrow C$ such that $f\xi _{i}=f_{i}$ for any $i\in
\mathbb{N}$. We prove that $f$ is a homomorphism of coalgebras. We
have:
\begin{equation*}
(f\otimes f)\Delta \xi _{i}=(f\otimes f)(\xi _{i}\otimes \xi
_{i})\Delta _{X_{i}}=(f_{i}\otimes f_{i})\Delta _{X_{i}}=\Delta
_{C}f_{i}=\Delta _{C}f\xi _{i};\text{\quad }\varepsilon _{C}f\xi
_{i}=\varepsilon _{C}f_{i}=\varepsilon _{X_{i}}=\varepsilon \xi
_{i.}
\end{equation*}
Since the relations above hold true for any $i\in \mathbb{N,}$ by
the universal property of the direct limit we deduce that
\begin{equation*}
(f\otimes f)\Delta =\Delta _{C}f;\text{\quad }\varepsilon
_{C}f=\varepsilon .
\end{equation*}
\end{proof}

\begin{claim}
\label{def of delta_n}Let $X$ be an object in an abelian monoidal category $%
\left( \mathcal{M},\otimes ,\mathbf{1}\right) $. Set%
\begin{equation*}
X^{\otimes 0}=\mathbf{1},\qquad X^{\otimes 1}=X\qquad
\text{and}\qquad X^{\otimes n}=X^{\otimes n-1}\otimes X,\text{ for
every }n>1
\end{equation*}%
and for every morphism $f:X\rightarrow Y$ in $\mathcal{M}$, set%
\begin{equation*}
f^{\otimes 0}=\mathrm{Id}_{\mathbf{1}},\qquad f^{\otimes 1}=f\qquad \text{and%
}\qquad f^{\otimes n}=f^{\otimes n-1}\otimes f,\text{ for every
}n>1.
\end{equation*}%
Let $\left( C,\Delta _{C},\varepsilon _{C}\right) $ be a coalgebra in $%
\mathcal{M}$ and for every $n\in
\mathbb{N}
,$ define the $n^{\text{th}}$ iterated comultiplication of $C,$
$\Delta _{C}^{n}:C\rightarrow C^{\otimes {n+1}}$, by
\begin{equation*}
\Delta _{C}^{0}=\text{Id}_{C},\qquad \Delta _{C}^{1}=\Delta _{C}\qquad \text{%
and}\qquad \Delta _{C}^{n}=\left( \Delta _{C}^{\otimes n-1}\otimes
C\right) \Delta _{C},\text{ for every }n>1.
\end{equation*}%
Let $\delta :D\rightarrow C$ be a monomorphism which is a
homomorphism of
coalgebras in $\mathcal{M}$. Denote by $(L,p)$ the cokernel of $\delta $ in $%
\mathcal{M}$. Regard $D$ as a $C$-bicomodule via $\delta $ and
observe that, by Proposition \ref{pro: bicomodules}, $L$ is a
$C$-bicomodule and $p$ is a morphism of bicomodules. Let
\begin{equation*}
(D^{\wedge _{C}^{n}},\delta _{n}):=\ker (p^{\otimes {n}}\Delta
_{C}^{n-1})
\end{equation*}%
for any $n\in \mathbb{N}\setminus \{0\}.$ Note that $(D^{\wedge
_{C}^{1}},\delta _{1})=(D,\delta )$ and $(D^{\wedge
_{C}^{2}},\delta _{2})=D\wedge _{C}D.$ \newline In order to
simplify the notations we set $(D^{\wedge _{C}^{0}},\delta
_{0})=(0,0).$\newline Now, if we assume that $\mathcal{M}$ has
left exact tensor functors, by Proposition \ref{pro: ker is a
coalgebra}, since $p^{\otimes {n}}\Delta _{C}^{n-1}$ is a morphism
of $C$-bicomodules (as a composition of morphisms
of $C$-bicomodules), we get that $D^{\wedge _{C}^{n}}$ is a coalgebra and $%
\delta _{n}:D^{\wedge _{C}^{n}}\rightarrow C$ is a coalgebra
homomorphism for any $n>0$ and hence for any $n\in \mathbb{N}$.
\end{claim}

\begin{proposition}\label{pro: limit of delta}Let $\delta:D\to C$ be a monomorphism which is a homomorphism of coalgebras in an abelian monoidal category
$\M$. Then for any $i\leq j$ in $\mathbb{N}$ there is a (unique)
morphism $\xi_{i}^j:D^{\wedge_C ^i}\to D^{\wedge_C ^j}$ such that
\begin{equation}\label{compatibility of delta i}
\delta_j\xi_{i}^j=\delta_i.
\end{equation} Moreover $\xi_{i}^j$ is a coalgebra
homomorphism and $((D^{\wedge_C ^i})_{i\in \mathbb{N}},(\xi
_{i}^j)_{i,j\in \mathbb{N}})$ is a direct system in $\M$ whose
direct limit, if it exists, carries a natural coalgebra structure
that makes it the direct limit of $((D^{\wedge_C ^i})_{i\in
\mathbb{N}},(\xi _{i}^j)_{i,j\in \mathbb{N}})$ as a direct system
of coalgebras.
\end{proposition}

\begin{proof}Set $D^i:=D^{\wedge_C ^i}$ for any $i\in \mathbb{N}.$
Consider the following diagram:
\begin{equation*}
\begin{diagram}[h=2em,w=3em] 0& \rTo& D^{i+1}&\rTo^{\delta_{i+1}}&  C
&\rTo^{p^{\ot {i+1}}\Delta_C^{i}}& \textit{L}^{\ot {i+1}}\\ &&&
\luDotsto<{\xi_i^{i+1}}& \uTo>{\delta_i}\\ &&&& D^i
\end{diagram}
\end{equation*}
Let $i>0$. Since $\delta _{i}=\ker (p^{\otimes i}\Delta
_{C}^{i-1})$ is a coalgebra homomorphism, we have:
\begin{eqnarray*}
p^{\otimes i+1}\Delta _{C}^{i}\delta _{i} &=&p^{\otimes
i+1}(C\otimes \Delta
_{C}^{i-1})\Delta _{C}\delta _{i} \\
&=&p^{\otimes i+1}(C\otimes \Delta _{C}^{i-1})(\delta _{i}\otimes
\delta _{i})\Delta _{D^{i}} =(p\delta _{i}\otimes p^{\otimes
i}\Delta _{C}^{i-1}\delta _{i})\Delta _{D^{i}}=0.
\end{eqnarray*}
Then, for any $i\geq 1,$ by the universal property of the kernel,
there exists a unique morphism $\xi _{i}^{i+1}:D^{i}\rightarrow
D^{i+1}$ such that $\delta _{i+1}\xi _{i}^{i+1}=\delta _{i}.$ Set
$\xi _{0}^{1}=0$ and for any $j>i,$ define:
\begin{equation*}
\xi _{i}^{j}=\xi _{j-1}^{j}\xi _{j-2}^{j-1}\cdots \xi
_{i+1}^{i+2}\xi _{i}^{i+1}:D^{i}\rightarrow D^{j}.
\end{equation*}
In such a way we obviously obtain a direct system in $\M.$ Let us
prove that $\xi _{i}^{j}$ is a homomorphism of coalgebras for any
$j>i.$ It is clearly sufficient to verify this for
$j=i+1.$\newline As $\delta _{i+1}$ and $\delta _{i}$ are
coalgebra homomorphisms, we have
\begin{equation*}
(\delta _{i+1}\otimes \delta _{i+1})\Delta _{D^{i+1}}\xi
_{i}^{i+1}=\Delta _{D}\delta _{i+1}\xi _{i}^{i+1}=\Delta
_{D}\delta _{i}=(\delta _{i}\otimes \delta _{i})\Delta
_{D^{i}}=(\delta _{i+1}\otimes \delta _{i+1})(\xi
_{i}^{i+1}\otimes \xi _{i}^{i+1})\Delta _{D^{i}}.
\end{equation*}
Since the tensor functors are left exact, $\delta _{i+1}\otimes
\delta _{i+1} $ is a monomorphism so that we get $\Delta
_{D^{i+1}}\xi _{i}^{i+1}=(\xi _{i}^{i+1}\otimes \xi
_{i}^{i+1})\Delta _{D^{i}}$. Moreover we have
\begin{equation*}
\varepsilon _{D^{i+1}}\xi _{i}^{i+1}=\varepsilon _{D}\delta
_{i+1}\xi _{i}^{i+1}=\varepsilon _{D}\delta _{i}=\varepsilon
_{D^{i}}.
\end{equation*}\\
The last assertion follows by Proposition \ref{pro:
colimit=coalg}.
\end{proof}

\begin{notation}\label{notation tilde}Let $\delta:D\to C$ be a homomorphism of coalgebras in a cocomplete abelian monoidal category
$\M$. By Proposition \ref{pro: limit of delta} $((D^{\wedge_C
^i})_{i\in \mathbb{N}},(\xi _{i}^j)_{i,j\in \mathbb{N}})$ is a
direct system in $\M$ whose direct limit carries a natural
coalgebra structure that makes it the direct limit of
$((D^{\wedge_C ^i})_{i\in \mathbb{N}},(\xi _{i}^j)_{i,j\in
\mathbb{N}})$ as a direct system of coalgebras.\\From now on we
set: $(\widetilde{D}_C, (\xi_i)_{i\in \mathbb{N}})=
\underrightarrow{\lim }(D^{\wedge_C ^i})_{i\in \mathbb{N}}$, where
$\xi_i:D^{\wedge_C ^i}\to \widetilde{D}_C$ denotes the structural
morphism of the direct limit. We simply write $\widetilde{D}$ if
there is no danger of confusion. We note that, since
$\widetilde{D}$ is a direct limit of coalgebras, the canonical
(coalgebra) homomorphisms $(\delta _{i}:D^{\wedge_C ^i}\rightarrow
C)_{i\in \mathbb{N}}$, which are compatible by (\ref{compatibility
of delta i}), factorize to a unique coalgebra homomorphism
$\widetilde{\delta} :\widetilde{D}\rightarrow C$ such that
$\widetilde{\delta} \xi_i=\delta_i$ for any $i\in \mathbb{N}$.
\end{notation}

\section{Cotensor coalgebra}

\begin{claim}
Let $(C,\Delta ,\varepsilon )$ be an coalgebra in $(\mathcal{M},\otimes ,%
\mathbf{1})$ and let $\left( L,\rho _{L}^{l},\rho _{L}^{r}\right) $ be a $C$%
-bicomodule. Recall that a morphism $\zeta :L\rightarrow C\otimes
C$ is called a Hochschild $2$-cocyle whenever
\begin{equation*}
b^{2}(\zeta )=(\zeta \otimes C)\circ \rho _{L}^{r}-(C\otimes
\Delta )\circ \zeta +(\Delta \otimes C)\circ \zeta -(C\otimes
\zeta )\circ \rho _{L}^{l}
\end{equation*}%
is zero. See \ref{standard complex}.
\end{claim}

\begin{definition}
Let $(C,\Delta _{C},\varepsilon _{C})$ be a coalgebra in $(\mathcal{L}%
,\otimes ,\mathbf{1})$ and let $(L,\rho _{l},\rho _{r})$ be a
bicomodule over $C$. A \emph{Hochschild extension }of $C$ with
cokernel $L$, is an exact sequence in $\mathcal{M}$:
\begin{equation}
0\rTo C\rTo^{\sigma }E{\rTo^{p}}L\rTo0  \tag{$E$}
\end{equation}%
that satisfy the following conditions:

a) $(E,\Delta_{E},\varepsilon_{E})$ is a coalgebra in
$\mathcal{M}$ ;

b) $\sigma $ is a homomorphism of coalgebras that has a retraction $\pi $ in $%
\mathcal{M}$;

c) $C\wedge_E C=E$, that is $(p\otimes p)\Delta =0$;

d) the morphisms $\rho_l$ and $\rho_r$ fulfill the following
relations
\begin{equation*}
\rho _{l} p = (\pi \otimes p) \Delta_{E} \qquad \rho _{r} p=
(p\otimes \pi) \Delta_{E}.
\end{equation*}
\end{definition}

The following result will lead to the definition of a coalgebra
structure for the cotensor coalgebra.

\begin{lemma}
\label{X lem dual 1.5.7}Let $(C,\Delta ,\varepsilon )$ be a
coalgebra in an abelian monoidal category $\mathcal{M}$ and let
$(L,\rho _{L}^{r},\rho _{L}^{l})$ a $C$-bicomodule. Suppose that
$\zeta :L\rightarrow C\otimes C $ is a morphism in $\mathcal{M}$.
\newline Define $\Delta _{\zeta }:C\oplus L\rightarrow \left(
C\oplus L\right) \otimes \left( C\oplus L\right) $ and
$\varepsilon _{\zeta }:C\oplus L\rightarrow \mathbf{1}$ by:
\begin{eqnarray}
\Delta _{\zeta }:= &&\left( i_{C}\otimes i_{C}\right) \Delta
p_{C}+\left[ \left( i_{L}\otimes i_{C}\right) \rho _{L}^{r}+\left(
i_{C}\otimes i_{L}\right) \rho _{L}^{l}-\left( i_{C}\otimes
i_{C}\right) \zeta \right]
p_{L},  \label{delta} \\
\varepsilon _{\zeta }:= &&\varepsilon p_{C}+l_{\mathbf{1}}\left(
\varepsilon \otimes \varepsilon \right) \zeta p_{L},
\label{epsilon}
\end{eqnarray}%
where $i_{C},\ i_{L}$ are the canonical injections in $C\oplus L$ and $%
p_{C},\ p_{L}$ are the canonical projections. Then $\Delta _{\zeta
}$ is a
coassociative comultiplication if and only if $\zeta $ is a Hochschild $2$%
-cocycle. Moreover, in this case, $\left( C\oplus L,\Delta _{\zeta
},\varepsilon _{\zeta }\right) $ is a coalgebra and
$i_{C}:C\rightarrow \left( C\oplus L,\Delta _{\zeta },\varepsilon
_{\zeta }\right) $ is a Hochschild extension of $C$ with cokernel
$(L,p_{L}).$ This extension will be denoted by $E_{\zeta }$.
\end{lemma}

\begin{proof}
The dual of an abelian monoidal category is an abelian monoidal
category in the sense of \cite[Definition 1.8]{AMS}.  The
conclusion follows by applying \cite[Lemma 2.5]{AMS}.
\end{proof}

Let $\mathcal{M}$ be an abelian monoidal category and let $C$ be a
coalgebra in $(\M,\ot,\mathbf{1})$. Given a right $C$-bicomodule
$(V,\rho^r_V)$ and a left $C$-comodule $(W,\rho^l_W)$, their
cotensor product over $C$ in $\mathcal{M}$ is defined to be the
equalizer $(V\cot_{C}W,{\chi})$ of the couple of morphism
$(\rho^r_V\ot W,V\ot \rho^l_W)$:
\begin{diagram}
0&\rTo&V\cot_{C}W &\rTo^{\chi}&V\otimes W&\pile{ \rTo^{\rho^r_V\ot
W}\\\rTo_{V\ot \rho^l_W}} &V\otimes C\otimes W.
\end{diagram}
Since the tensor functors are left exact, in view of Proposition
\ref{pro: bicomodules}, then $V\Box _{C}W$ is also a
$C$-bicomodule, namely it is $C$-sub-bicomodule of $V\ot W$,
whenever $V$ and $W$ are $C$-bicomodules. Furthermore, in this
case, the category $({^C\M^C},\cot_C,C)$ is still an abelian
monoidal category; the associative and unit constraints are
induced by the ones in $\M$ (the proof is dual to \cite[Theorem
1.11]{AMS}). Therefore, also using $\cot_C$, one can forget about
brackets. Moreover the functors $M\cot_C(-):{^C\M}\to \M$ and
$(-)\cot_C M:{\M^C}\to\M$ are left exact for any $M\in\M.$

We will write $\cot$ instead of $\cot_C$, whenever there is no
danger of misunderstanding.

\begin{claim}
\label{claim 4.2}Let $(C,\Delta ,\varepsilon )$ be a coalgebra in
a cocomplete abelian monoidal category $\mathcal{M}$ and let $(M,\rho _{M}^{r},\rho _{M}^{l})$ be a $C$--$\,$%
{}{}bicomodule. Set
\begin{equation*}
M^{\square 0}=C,M^{\square 1}=M\text{\quad and\quad }M^{\square
n}=M^{\square n-1}\square M\text{ for any }n>1
\end{equation*}%
and define $(C^{n}(M))_{n\in \mathbb{N}}$ by
\begin{equation*}
C^{0}(M)=0,C^{1}(M)=C\text{\quad and\quad
}C^{n}(M)=C^{n-1}(M)\oplus M^{\square n-1}\text{ for any }n>1.
\end{equation*}%
Let $\sigma _{i}^{i+1}:C^{i}(M)\rightarrow C^{i+1}(M)$ be the
canonical inclusion and for any $j>i,$ define:
\begin{equation*}
\sigma _{i}^{j}=\sigma _{j-1}^{j}\sigma _{j-2}^{j-1}\cdots \sigma
_{i+1}^{i+2}\sigma _{i}^{i+1}:C^{i}(M)\rightarrow C^{j}(M).
\end{equation*}%
Then $(({C^{i}}(M))_{i\in \mathbb{N}},(\sigma _{i}^{j})_{i,j\in
\mathbb{N}})$ is a direct system in $\mathcal{M}$. We set
\begin{equation*}
T_{C}^{c}(M)=\bigoplus_{n\in \mathbb{N}}M^{\square n}=C\oplus
M\oplus M^{\square 2}\oplus M^{\square 3}\oplus \cdots
\end{equation*}%
and we denote by $\sigma _{i}:C^{i}(M)\rightarrow T_{C}^{c}(M)$
the canonical inclusion.\medskip\newline Throughout let
\begin{eqnarray*}
\pi _{n}^{m} &:&C^{n}(M)\rightarrow C^{m}(M)\text{ }(m\leq n),\text{\qquad }%
\pi _{n}:T_{C}^{c}(M)\rightarrow C^{n}(M), \\
p_{n}^{m} &:&C^{n}(M)\rightarrow M^{\square m}\text{ }(m<n),\text{\qquad }%
p_{n}:T_{C}^{c}(M)\rightarrow M^{\square n},
\end{eqnarray*}
be the canonical projections and let
\begin{eqnarray*}
\sigma _{m}^{n} &:&C^{m}(M)\rightarrow C^{n}(M)\text{ }(m\leq n),\text{%
\qquad }\sigma _{n}:C^{n}(M)\rightarrow T_{C}^{c}(M), \\
i_{m}^{n} &:&M^{\square m}\rightarrow C^{n}(M)\text{ }(m<n),\text{\qquad }%
i_{m}:M^{\square m}\rightarrow T_{C}^{c}(M),
\end{eqnarray*}%
be the canonical injection for any $m,n\in \mathbb{N}$.\newline
For technical reasons we set $\pi _{n}^{m}=0,$ $\sigma _{m}^{n}=0$ for any $n<m$ and $%
p_{n}^{m}=0$, $i_{m}^{n}=0$ for any $n\leq m.$ Then, we have the
following relations:%
\begin{equation*}
p_{n}\sigma _{k}=p_{k}^{n},\qquad p_{n}i_{k}=\delta _{n,k}\mathrm{Id}%
_{M^{\square k}},\qquad \pi _{n}i_{k}=i_{k}^{n}.
\end{equation*}%
Moreover, we have:
\begin{equation*}
\begin{tabular}{lll}
$\pi _{n}^{m}\sigma _{k}^{n}=\sigma _{k}^{m},\text{ if }k\leq m\leq n,$ & $%
\qquad $and$\qquad $ & $\pi _{n}^{m}\sigma _{k}^{n}=\pi _{k}^{m},\text{ if }%
m\leq k\leq n,$ \\
$p_{n}^{m}\pi _{k}^{n}=p_{k}^{m},\text{ if }m<n\leq k,$ & $\qquad $and$%
\qquad $ & $\sigma _{n}^{m}i_{k}^{n}=i_{k}^{m},\text{ if }k<n\leq m,$ \\
$p_{n}^{m}\sigma _{k}^{n}=p_{k}^{m},\text{ if }m<k\leq n,$ & $\qquad $and$%
\qquad $ & $\pi _{n}^{m}i_{k}^{n}=i_{k}^{m},\text{ if }k<m\leq n,$ \\
$p_{n}^{m}\pi _{n}=p_{m},\text{ if }m<n,$ & $\qquad $and$\qquad $
& $\sigma
_{n}i_{m}^{n}=i_{m},\text{ if }m<n\text{,}$ \\
$\pi _{n}\sigma _{k}=\sigma _{k}^{n},\text{ if }k\leq n,$ & $\qquad $and$%
\qquad $ & $\pi _{n}\sigma _{k}=\pi _{k}^{n},\text{ if }n\leq k,$ \\
$p_{n}^{m}i_{m}^{n}=\mathrm{Id}_{M^{\square m}},\text{ if }m<n.$ &
&
\end{tabular}%
\end{equation*}
In the other cases, these compositions are zero.
\end{claim}

\begin{proposition}
\label{pro: V=lim V^n} Let $\mathcal{M}\ $ be a cocomplete abelian category. Let $(V_{i})_{i\in \mathbb{N}}$ be a family of objects in $\mathcal{%
M}$ and let $\left( V,v_{i}\right) =\oplus _{i\in
\mathbb{N}
}V_{i}$ be the direct sum of the family $(V_{i})_{i\in
\mathbb{N}}$. Then
\begin{equation*}
(V,\nabla _{i=0}^{n}(v_{i}))=\underrightarrow{\lim }(\oplus
_{i=0}^{n}V_{i}),
\end{equation*}%
where $\nabla _{i=0}^{n}(v_{i}):\oplus _{i=0}^{n}V_{i}\to V$
denotes the codiagonal morphism associated to the family $\left(
v_{i}\right) _{i=0}^{n}.$
\end{proposition}

\begin{proof}
Set $V^{n}:=\oplus _{i=0}^{n}V_{i},$ for any $n\in \mathbb{N}$, and let $%
w_{m}^{n}:V^{m}\rightarrow V^{n}$ be the canonical inclusion for
$m\leq n.$
Let $(f_{n}:V^{n}\rightarrow X)_{n}$ be a compatible family of morphisms in $%
\mathcal{M}$, i.e. $f_{n}w_{m}^{n}=f_{m}$ for any $m\leq n.$ Let $%
v_{m}^{n}:V_{m}\rightarrow V^{n}$ be the canonical inclusion for every $%
m\leq n$ and let $v_{m}^{n}=0$ otherwise. Note that the morphism
$\nabla _{i=0}^{n}(v_{i}):V^{n}\rightarrow V$ is uniquely defined
by the following relation:
\begin{equation*}
\left[ \nabla _{i=0}^{n}(v_{i})\right] v_{m}^{n}=v_{m},\text{ for every }%
m\leq n.\text{ }
\end{equation*}%
Observe that, for every $m\leq n\leq t,$ we have%
\begin{equation*}
f_{t}v_{m}^{t}=f_{t}w_{n}^{t}v_{m}^{n}=f_{n}v_{m}^{n}
\end{equation*}%
so that, by the universal property of the direct sum, there exists
a unique morphism $f:V\rightarrow X$ such that
\begin{equation}
fv_{m}=f_{n}v_{m}^{n},  \label{unive f}
\end{equation}
for any $m\in \mathbb{N}$, where $n\in \mathbb{N}$ and $m\leq n.$
Thus
\begin{equation*}
f_{n}v_{m}^{n}=fv_{m}=f\left[ \nabla _{i=0}^{n}(v_{i})\right]
v_{m}^{n}\text{ for every }m\leq n.
\end{equation*}%
By the universal property of $V^{n}:=\oplus _{i=0}^{n}V_{i},$
$f_{n}$ is the unique morphism that composed with $v_{m}^{n}$
gives $f_{n}v_{m}^{n}$ for any $m\leq n.$ We get that
\begin{equation*}
f_{n}=f\left[ \nabla _{i=0}^{n}(v_{i})\right] ,\text{ for every
}n\in \mathbb{N}.
\end{equation*}%
In order to conclude that $V=\underrightarrow{\lim }V^{i},$ it
remains to prove that $f:V\rightarrow X$ is the unique morphism
with this property. Let $g:V\rightarrow X$ be a morphism such that
$f_{n}=g\left[ \nabla _{i=0}^{n}(v_{i})\right] $ for every $n\in
\mathbb{N}$. Then
\begin{equation*}
fv_{m}=f_{n}v_{m}^{n}=g\left[ \nabla _{i=0}^{n}(v_{i})\right]
v_{m}^{n}=gv_{m}\text{ for every }m,n\in \mathbb{N},m\leq n.
\end{equation*}%
By uniqueness of $f$ with respect to (\ref{unive f}), we get
$g=f.$
\end{proof}

\begin{corollary}\label{pro: T=lim c^n}Let $(C,\Delta ,\varepsilon )$ be a
coalgebra in a cocomplete abelian monoidal category $\M$ and let $(M,\rho _{M}^{r},\rho _{M}^{l})$ be a $%
C$-bicomodule.
Then $$ (T^c_C(M),(\sigma_n)_{n\in \N})=\underrightarrow{%
\lim }C^{n}(M).$$
\end{corollary}

\begin{proof}
  Just observe that $C^n(M)=\oplus
_{m=0}^{n-1}M^{\cot m}$ and $\sigma_n=\nabla _{m=0}^{n-1}(i_{m}).$
\end{proof}

\begin{claim}
 Note that $M^{\cot n}$ is a $C$-bicomodule via $\rho^l_n:=\rho^l_M\cot M^{\cot n-1}$ and
$\rho^r_n:=M^{\cot n-1}\cot \rho^r_M$.

Our next aim is to define, for any $n\in \N\setminus\{0\}$, a
Hochschild 2-cocycle
$$\zeta^n :M^{\cot n}\rightarrow C^{n}(M)\otimes C^{n}(M).$$
Then we will apply Lemma \ref{X lem dual 1.5.7} to obtain that,
for any $n>0$, $C^{n+1}(M)=C^{n}(M)\oplus M^{\cot n}$ can be
endowed with a coalgebra structure $ \left( C^{n+1}(M),\Delta
_{\zeta ^{n}},\varepsilon _{\zeta^{n} }\right) $ in $\mathcal{M}$
such that the canonical inclusion
$\sigma_{n}^{n+1}:C^{n}(M)\rightarrow C^{n+1}(M)$ is a Hochschild
extension of $C^{n}(M)$ with cokernel $M^{\cot n}$. Then, by
Proposition \ref{pro: colimit=coalg}, $ T^c_C(M)$ will carry a
natural coalgebra structure that makes it the direct limit of
$((C^{i}(M))_{i\in \mathbb{N}},(\sigma _{i}^j)_{i,j\in
\mathbb{N}})$ as a direct system of
coalgebras.\medskip\newline
Let $\chi _{M}:M\cot M\rightarrow M\otimes M$ be the canonical
inclusion and define
\begin{equation}\zeta ^{1}=0\text{\quad and\quad }\zeta
^{n}=-\sum_{t=1}^{n-1}(i_{t}^{n}\otimes i_{n-t}^{n})(M^{\cot
t-1}\cot \chi_M\cot M^{\cot n-1-t}), \forall n>1.   \label{def of
zeta}
\end{equation} where we identify $C\cot X$ and $X\cot C$ with $X,$ for
any $C$-bicomodule $X$.
\end{claim}

\begin{proposition}
\label{Pro: cotenso coalgebra}Let $(C,\Delta ,\varepsilon )$ be a
coalgebra in a cocomplete abelian monoidal category $\mathcal{M}$
and let $(M,\rho _{M}^{r},\rho _{M}^{l})$ be a $C$-bicomodule. Let
\begin{gather*}
\Delta (1)=\Delta ,\qquad \text{and}\qquad \varepsilon (1)=\varepsilon , \\
\overline{\rho }_{1}^{l}=\rho _{M}^{l},\qquad \text{and}\qquad \overline{%
\rho }_{1}^{r}=\rho _{M}^{r},
\end{gather*}%
and for every $n\geq 2$ set%
\begin{equation*}
\Delta (n)=\Delta _{\zeta ^{n-1},}\qquad \text{and}\qquad
\varepsilon (n)=\varepsilon _{\zeta ^{n-1}},
\end{equation*}%
as defined in (\ref{delta}), (\ref{epsilon}), and let%
\begin{equation*}
\overline{\rho }_{n}^{l}=(\sigma _{1}^{n}\otimes M^{\square
n})(\rho
_{M}^{l}\square M^{\square n-1}),\qquad \text{and}\qquad \overline{\rho }%
_{n}^{r}=(M^{\square n}\otimes \sigma _{1}^{n})(M^{\square
n-1}\square \rho _{M}^{r}).
\end{equation*}%
Then, for any $n\geq 1,$ we have

\begin{enumerate}
\item[a)] $(C^{n}(M),\Delta (n),\varepsilon (n))\ $is a coalgebra.

\item[b)] $(M^{\square n},\overline{\rho }_{n}^{l},\overline{\rho
}_{n}^{r})$ is a $C^{n}(M)$-bicomodule such that the morphism
$\zeta ^{n}:M^{\square n}\rightarrow C^{n}(M)\otimes C^{n}(M)$,
given by (\ref{def of zeta}), defines a Hochschild 2-cocycle.

\item[c)] $\varepsilon (n)=\varepsilon _{C}\pi _{n}^{1}.$

\item[d)] For every $1\leq t\leq n-1$ we have%
\begin{equation}
\Delta (n)i_{t}^{n}=\left( i_{t}^{n}\otimes \sigma _{t}^{n}\right) \overline{%
\rho }_{t}^{r}+\left( \sigma _{t}^{n}\otimes i_{t}^{n}\right)
\overline{\rho }_{t}^{l}-\left( \sigma _{t}^{n}\otimes \sigma
_{t}^{n}\right) \zeta ^{t}. \label{formula delta}
\end{equation}

\item[e)] $\Delta (n)$ fulfils the following relations:
\begin{eqnarray*}
\Delta (n)i_{0}^{n} &=&\left( i_{0}^{n}\otimes i_{0}^{n}\right) \Delta  \\
\Delta (n)i_{1}^{n} &=&\left( i_{1}^{n}\otimes i_{0}^{n}\right)
\rho
_{M}^{r}+\left( i_{0}^{n}\otimes i_{1}^{n}\right) \rho _{M}^{l} \\
\text{and for }2 &\leq &t\leq n-1 \\
\Delta (n)i_{t}^{n} &=&\left( i_{t}^{n}\otimes i_{0}^{n}\right)
(M^{\square t-1}\square \rho _{M}^{r})+\left( i_{0}^{n}\otimes
i_{t}^{n}\right) (\rho
_{M}^{l}\square M^{\square t-1})+ \\
&&+\sum_{r=1}^{t-1}(i_{r}^{n}\otimes i_{t-r}^{n})(M^{\square
r-1}\square \chi _{M}\square M^{\square t-1-r}).
\end{eqnarray*}
\end{enumerate}
\end{proposition}

\begin{proof}
Set $C^{n}=C^{n}(M)$ for any $n\geq 1$. Recall that $M^{\square n}$ is a $C$%
-bicomodule via
\begin{equation*}
\rho _{n}^{l}:=\rho _{M}^{l}\square M^{\square n-1},\qquad
\text{and}\qquad \rho _{n}^{r}:=M^{\square n-1}\square \rho
_{M}^{r}.
\end{equation*}%
Let us prove all the statements of the theorem by induction on $n\mathbb{%
\geq }1.$

If $n=1,$ then $C^{1}=C$ is a coalgebra and $M^{\square 1}=M$ is a $C^{1}$%
-bicomodule by hypothesis. Obviously $\zeta ^{1}=0$ fulfills (\ref{2 cocycle}%
). We have $\varepsilon (1)=\varepsilon =\varepsilon _{C}\pi
_{1}^{1}$ and, since $i_{0}^{1}=\mathrm{Id}_{C}$ and
$i_{1}^{1}=0,$ we get
\begin{equation*}
\Delta (1)i_{0}^{1}=\Delta =\left( i_{0}^{1}\otimes
i_{0}^{1}\right) \Delta ,\qquad \Delta (1)i_{1}^{1}=0=\left(
i_{1}^{1}\otimes i_{0}^{1}\right) \rho _{M}^{r}+\left(
i_{0}^{1}\otimes i_{1}^{1}\right) \rho _{M}^{l}.
\end{equation*}

Let $n\geq 2$. Assume all the assertions hold true for any $1\leq
t<n.$
\newline
Thus $(C^{n-1},\Delta (n-1),\varepsilon (n-1))$ is a coalgebra in $\mathcal{M%
}$, $(M^{\square n-1},\overline{\rho }_{n-1}^{l},\overline{\rho
}_{n-1}^{r})$
is a $C^{n-1}$-bicomodule and $\zeta ^{n-1}$ is a $2$-cocycle. By Lemma \ref%
{X lem dual 1.5.7} applied to $"C"=C^{n-1}$ and $"M"=M^{\square n-1}$, then $%
(C^{n},\Delta (n),\varepsilon (n))$ is a coalgebra$.$ Moreover
$\sigma _{t-1}^{t}:C^{t-1}\rightarrow C^{t}$ is a Hochschild
extension of $C^{t-1}$ with cokernel $M^{\square t-1}$ for any
$1\leq t<n$. \newline
Since $(M^{\square n},\rho _{n}^{l},\rho _{n}^{r})$ is a $C$-bicomodule and $%
\sigma _{1}^{n}:C\rightarrow C^{n}$ is a coalgebra homomorphism
(as a
composition of coalgebra homomorphism), then $(M^{\square n},\overline{\rho }%
_{n}^{l},\overline{\rho }_{n}^{r})$ is a $C^{n}$-bicomodule, where
\begin{equation*}
\overline{\rho }_{n}^{l}=(\sigma _{1}^{n}\otimes M^{\square
n})(\rho
_{M}^{l}\square M^{\square n-1}),\qquad \text{and}\qquad \overline{\rho }%
_{n}^{r}=(M^{\square n}\otimes \sigma _{1}^{n})(M^{\square
n-1}\square \rho _{M}^{r}).
\end{equation*}%
Recall that, by definition, we have:
\begin{multline*}
\Delta (n)=\Delta _{\zeta ^{n-1}} \\
=\left( \sigma _{n-1}^{n}\otimes \sigma _{n-1}^{n}\right) \Delta
(n-1)\pi _{n}^{n-1}+\left[
\begin{array}{c}
\left( i_{n-1}^{n}\otimes \sigma _{n-1}^{n}\right) \overline{\rho }%
_{n-1}^{r}+\left( \sigma _{n-1}^{n}\otimes i_{n-1}^{n}\right)
\overline{\rho
}_{n-1}^{l}+ \\
-\left( \sigma _{n-1}^{n}\otimes \sigma _{n-1}^{n}\right) \zeta ^{n-1}%
\end{array}%
\right] p_{n}^{n-1},
\end{multline*}%
For any $0\leq t\leq n-1,$ we have that%
\begin{eqnarray*}
&&\left[ \left( i_{n-1}^{n}\otimes \sigma _{n-1}^{n}\right) \overline{\rho }%
_{n-1}^{r}+\left( \sigma _{n-1}^{n}\otimes i_{n-1}^{n}\right)
\overline{\rho }_{n-1}^{l}-\left( \sigma _{n-1}^{n}\otimes \sigma
_{n-1}^{n}\right) \zeta
^{n-1}\right] p_{n}^{n-1}i_{t}^{n} \\
&=&\left[ \left( i_{n-1}^{n}\otimes \sigma _{n-1}^{n}\right) \overline{\rho }%
_{n-1}^{r}+\left( \sigma _{n-1}^{n}\otimes i_{n-1}^{n}\right)
\overline{\rho }_{n-1}^{l}-\left( \sigma _{n-1}^{n}\otimes \sigma
_{n-1}^{n}\right) \zeta ^{n-1}\right] \delta _{t,n-1}
\end{eqnarray*}%
so that we obtain
\begin{equation}
\Delta (n)i_{t}^{n}=\left( \sigma _{n-1}^{n}\otimes \sigma
_{n-1}^{n}\right) \Delta (n-1)\pi _{n}^{n-1}i_{t}^{n}+\left[
\begin{array}{c}
\left( i_{n-1}^{n}\otimes \sigma _{n-1}^{n}\right) \overline{\rho }%
_{n-1}^{r}+\left( \sigma _{n-1}^{n}\otimes i_{n-1}^{n}\right)
\overline{\rho
}_{n-1}^{l} \\
-\left( \sigma _{n-1}^{n}\otimes \sigma _{n-1}^{n}\right) \zeta ^{n-1}%
\end{array}%
\right] \delta _{t,n-1}.  \label{Chissa}
\end{equation}

For $t=0<n-1,$ we have:
\begin{equation*}
\Delta (n)i_{0}^{n}=\left( \sigma _{n-1}^{n}\otimes \sigma
_{n-1}^{n}\right) \Delta (n-1)i_{0}^{n-1},
\end{equation*}%
so that, inductively we get:
\begin{equation*}
\Delta (n)i_{0}^{n}=\left( \sigma _{n-1}^{n}\otimes \sigma
_{n-1}^{n}\right) \left( i_{0}^{n-1}\otimes i_{0}^{n-1}\right)
\Delta _{C}=\left( i_{0}^{n}\otimes i_{0}^{n}\right) \Delta _{C}.
\end{equation*}%
Let us prove that for every $t\in
\mathbb{N}
,$ such that $0<t\leq n-1$ we have (\ref{formula delta}).\newline
Now, we apply (\ref{Chissa}). If $t=n-1,$ we get%
\begin{equation*}
\Delta (n)i_{t}^{n}=\left( i_{n-1}^{n}\otimes \sigma
_{n-1}^{n}\right) \overline{\rho }_{n-1}^{r}+\left( \sigma
_{n-1}^{n}\otimes i_{n-1}^{n}\right) \overline{\rho
}_{n-1}^{l}-\left( \sigma _{n-1}^{n}\otimes \sigma
_{n-1}^{n}\right) \zeta ^{n-1}.
\end{equation*}%
If $0<t<n-1$, we get%
\begin{eqnarray*}
\Delta (n)i_{t}^{n} &=&\left( \sigma _{n-1}^{n}\otimes \sigma
_{n-1}^{n}\right) \Delta (n-1)i_{t}^{n-1} \\
&=&\left( \sigma _{n-1}^{n}\otimes \sigma _{n-1}^{n}\right) \left[
\left( i_{t}^{n-1}\otimes \sigma _{t}^{n-1}\right) \overline{\rho
}_{t}^{r}+\left( \sigma _{t}^{n-1}\otimes i_{t}^{n-1}\right)
\overline{\rho }_{t}^{l}-\left(
\sigma _{t}^{n-1}\otimes \sigma _{t}^{n-1}\right) \zeta ^{t}\right]  \\
&=&\left( i_{t}^{n}\otimes \sigma _{t}^{n}\right) \overline{\rho }%
_{t}^{r}+\left( \sigma _{t}^{n}\otimes i_{t}^{n}\right) \overline{\rho }%
_{t}^{l}-\left( \sigma _{t}^{n}\otimes \sigma _{t}^{n}\right)
\zeta ^{t}.
\end{eqnarray*}%
Thus we have obtained (\ref{formula delta}). Note that, for $t=1,$
by
definition of $\overline{\rho }_{1}^{l},\overline{\rho }_{1}^{r}$ and since $%
\zeta ^{1}=0,$ one gets
\begin{equation*}
\Delta (n)i_{1}^{n}=\left( i_{1}^{n}\otimes \sigma _{1}^{n}\right) \overline{%
\rho }_{1}^{r}+\left( \sigma _{1}^{n}\otimes i_{1}^{n}\right)
\overline{\rho }_{1}^{l}=\left( i_{1}^{n}\otimes i_{0}^{n}\right)
\rho _{M}^{r}+\left( i_{0}^{n}\otimes i_{1}^{n}\right) \rho
_{M}^{l}.
\end{equation*}%
\newline
For $2\leq t\leq n-1$, by definition of $\overline{\rho }_{n}^{l},\overline{%
\rho }_{n}^{r}$ and of $\zeta ^{t},$ from (\ref{formula delta}),
we get:
\begin{eqnarray*}
\Delta (n)i_{t}^{n} &=&\left( i_{t}^{n}\otimes \sigma
_{t}^{n}\right) (M^{\square t}\otimes \sigma _{1}^{t})(M^{\square
t-1}\square \rho _{M}^{r})+\left( \sigma _{t}^{n}\otimes
i_{t}^{n}\right) (\sigma
_{1}^{t}\otimes M^{\square t})(\rho _{M}^{l}\square M^{\square t-1})+ \\
&&+\sum_{r=1}^{t-1}\left( \sigma _{t}^{n}\otimes \sigma
_{t}^{n}\right) (i_{r}^{t}\otimes i_{t-r}^{t})(M^{\square
r-1}\square \chi _{M}\square M^{\square t-1-r}),
\end{eqnarray*}%
so that we obtain:
\begin{eqnarray*}
\Delta (n)i_{t}^{n} &=&\left( i_{t}^{n}\otimes i_{0}^{n}\right)
(M^{\square t-1}\square \rho _{M}^{r})+\left( i_{0}^{n}\otimes
i_{t}^{n}\right) (\rho
_{M}^{l}\square M^{\square t-1})+ \\
&&+\sum_{r=1}^{t-1}(i_{r}^{n}\otimes i_{t-r}^{n})(M^{\square
r-1}\square \chi _{M}\square M^{\square t-1-r}).
\end{eqnarray*}%
Moreover, by definition, we have:
\begin{equation}
\varepsilon (n)=\varepsilon _{\zeta ^{n-1}}=\varepsilon (n-1)\pi
_{n}^{n-1}+l_{\mathbf{1}}[\varepsilon (n-1)\otimes \varepsilon
(n-1)]\zeta ^{n-1}p_{n}^{n-1}.  \label{formula epsilon(n)}
\end{equation}%
Since $n\geq 2$ and $\varepsilon (n-1)=\varepsilon _{C}\pi _{n-1}^{1},$ by (%
\ref{formula epsilon(n)}), we have
\begin{equation*}
\varepsilon (n)=\varepsilon _{C}\pi _{n-1}^{1}\pi _{n}^{n-1}+l_{\mathbf{1}%
}(\varepsilon _{C}\pi _{n-1}^{1}\otimes \varepsilon _{C}\pi
_{n-1}^{1})\zeta ^{n-1}p_{n}^{n-1}=\varepsilon _{C}\pi
_{n}^{1}+l_{\mathbf{1}}(\varepsilon _{C}\otimes \varepsilon
_{C})(\pi _{n-1}^{1}\otimes \pi _{n-1}^{1})\zeta
^{n-1}p_{n}^{n-1}.
\end{equation*}

By definition of $\zeta ^{n-1},$ if $n=2$ we have $\zeta ^{1}=0$ so that $%
(\pi _{n-1}^{1}\otimes \pi _{n-1}^{1})\zeta ^{n-1}=0.$ If $n\geq
3,$ we have
\begin{equation*}
(\pi _{n-1}^{1}\otimes \pi _{n-1}^{1})\zeta
^{n-1}=\sum_{t=1}^{n-2}(\pi _{n-1}^{1}\otimes \pi
_{n-1}^{1})(i_{t}^{n-1}\otimes i_{n-1-t}^{n-1})(M^{\square
t-1}\square \chi _{M}\square M^{\square n-2-t})=0.
\end{equation*}%
We conclude that $(\pi _{n}^{1}\otimes \pi _{n}^{1})\zeta ^{n-1}=0$ for any $%
n\geq 2$ and hence $\varepsilon (n)=\varepsilon _{C}\pi
_{n}^{1}.$\newline Recall that $\zeta ^{n}$ is a $2$-cocycle means
that it verifies the following relation:
\begin{equation}
0=b^{2}(\zeta ^{n})=(\zeta ^{n}\otimes C^{n})\overline{\rho }%
_{n}^{r}-[C^{n}\otimes \Delta (n)]\zeta ^{n}+[\Delta (n)\otimes
C^{n}]\zeta ^{n}-(C^{n}\otimes \zeta ^{n})\overline{\rho
}_{n}^{l}.  \label{2 cocycle}
\end{equation}

Now, since $\sigma _{1}^{n}=i_{0}^{n}$ and $n\geq 2,$ we have
\begin{eqnarray*}
-(\zeta ^{n}\otimes C^{n})\overline{\rho }_{n}^{r}
&=&\sum_{t=1}^{n-1}[(i_{t}^{n}\otimes i_{n-t}^{n})(M^{\square
t-1}\square \chi _{M}\square M^{\square n-1-t})\otimes
C^{n}](M^{\square n}\otimes
\sigma _{1}^{n})(M^{\square n-1}\square \rho _{M}^{r}) \\
&=&\sum_{t=1}^{n-1}(i_{t}^{n}\otimes i_{n-t}^{n}\otimes \sigma
_{1}^{n})(M^{\square t-1}\square \chi _{M}\square M^{\square
n-1-t}\otimes
C)(M^{\square n-1}\square \rho _{M}^{r}) \\
&=&\sum_{t=1}^{n-1}\left( i_{t}^{n}\otimes i_{n-t}^{n}\otimes
i_{0}^{n}\right) (M^{\square t}\otimes M^{\square n-1-t}\square
\rho _{M}^{r})(M^{\square t-1}\square \chi _{M}\square M^{\square
n-1-t})
\end{eqnarray*}
and {\small
\begin{align*}
& -[C^{n}\otimes \Delta (n)]\zeta ^{n} \\
& =[C^{n}\otimes \Delta (n)]\sum_{t=1}^{n-1}(i_{t}^{n}\otimes
i_{n-t}^{n})(M^{\square t-1}\square \chi _{M}\square M^{\square n-1-t}) \\
& =\sum_{t=1}^{n-1}(i_{t}^{n}\otimes \Delta
(n)i_{n-t}^{n})(M^{\square
t-1}\square \chi _{M}\square M^{\square n-1-t}) \\
& =\sum_{t=1}^{n-1}\{i_{t}^{n}\otimes \left[
\begin{array}{c}
\left( i_{n-t}^{n}\otimes i_{0}^{n}\right) (M^{\square
n-t-1}\square \rho _{M}^{r})+\left( i_{0}^{n}\otimes
i_{n-t}^{n}\right) (\rho _{M}^{l}\square
M^{\square n-t-1})+ \\
+\sum_{r=1}^{n-t-1}(i_{r}^{n}\otimes i_{n-t-r}^{n})(M^{\square
r-1}\square
\chi _{M}\square M^{\square n-t-1-r})%
\end{array}%
\right] \}(M^{\square t-1}\square \chi _{M}\square M^{\square n-1-t}) \\
& =\sum_{t=1}^{n-1}\left( i_{t}^{n}\otimes i_{n-t}^{n}\otimes
i_{0}^{n}\right) (M^{\square t}\otimes M^{\square n-t-1}\square
\rho
_{M}^{r})(M^{\square t-1}\square \chi _{M}\square M^{\square n-1-t})+ \\
& +\sum_{t=1}^{n-1}\left( i_{t}^{n}\otimes i_{0}^{n}\otimes
i_{n-t}^{n}\right) (M^{\square t}\otimes \rho _{M}^{l}\square
M^{\square
n-t-1})(M^{\square t-1}\square \chi _{M}\square M^{\square n-1-t})+ \\
& +\sum_{t=1}^{n-1}\sum_{r=1}^{n-t-1}(i_{t}^{n}\otimes
i_{r}^{n}\otimes i_{n-t-r}^{n})(M^{\square t}\otimes M^{\square
r-1}\square \chi _{M}\square M^{\square n-t-1-r})(M^{\square
t-1}\square \chi _{M}\square M^{\square
n-1-t}) \\
& =-(\zeta ^{n}\otimes C^{n})\overline{\rho
}_{n}^{r}+\sum_{t=1}^{n-1}\left( i_{t}^{n}\otimes i_{0}^{n}\otimes
i_{n-t}^{n}\right) (M^{\square t}\otimes \rho _{M}^{l}\square
M^{\square n-t-1})(M^{\square t-1}\square \chi
_{M}\square M^{\square n-1-t})+ \\
& +\sum_{t=1}^{n-1}\sum_{r=1}^{n-t-1}(i_{t}^{n}\otimes
i_{r}^{n}\otimes i_{n-t-r}^{n})(M^{\square t}\otimes M^{\square
r-1}\square \chi _{M}\square M^{\square n-t-1-r})(M^{\square
t-1}\square \chi _{M}\square M^{\square n-1-t})
\end{align*}%
} Analogously one has
\begin{equation*}
-(C^{n}\otimes \zeta ^{n})\overline{\rho }_{n}^{l}=%
\sum_{t=1}^{n-1}(i_{0}^{n}\otimes i_{t}^{n}\otimes
i_{n-t}^{n})(\rho _{M}^{l}\square M^{\square t-1}M^{\square
n-t})(M^{\square t-1}\square \chi _{M}\square M^{\square n-1-t})
\end{equation*}%
and {\small
\begin{eqnarray*}
&&\lbrack \Delta (n)\otimes C^{n}]\zeta ^{n} \\
&=&-(C^{n}\otimes \zeta ^{n})\overline{\rho
}_{n}^{l}+\sum_{t=1}^{n-1}\left( i_{t}^{n}\otimes i_{0}^{n}\otimes
i_{n-t}^{n}\right) (M^{\square t-1}\square \rho _{M}^{r}\otimes
M^{\square n-t})(M^{\square t-1}\square \chi
_{M}\square M^{\square n-1-t})+ \\
&&+\sum_{t=1}^{n-1}\sum_{r=1}^{t-1}(i_{r}^{n}\otimes
i_{t-r}^{n}\otimes i_{n-t}^{n})(M^{\square r-1}\square \chi
_{M}\square M^{\square t-1-r}\otimes M^{\square n-t})(M^{\square
t-1}\square \chi _{M}\square M^{\square n-1-t}).
\end{eqnarray*}%
} Now, for any $1\leq t\leq n-1$, by definition of $\chi _{M}$ we
have
\begin{eqnarray*}
&&\left( i_{t}^{n}\otimes i_{0}^{n}\otimes i_{n-t}^{n}\right)
(M^{\square t-1}\square \rho _{M}^{r}\otimes M^{\square
n-t})(M^{\square t-1}\square
\chi _{M}\square M^{\square n-1-t}) \\
&=&\left( i_{t}^{n}\otimes i_{0}^{n}\otimes i_{n-t}^{n}\right)
(\left[ M^{\square t-1}\square (\rho _{M}^{r}\otimes M)\chi
_{M}\square M^{\square
n-1-t}\right] \\
&=&\left( i_{t}^{n}\otimes i_{0}^{n}\otimes i_{n-t}^{n}\right)
(\left[ M^{\square t-1}\square (M\otimes \rho _{M}^{l})\chi
_{M}\square M^{\square
n-1-t}\right] \\
&=&\left( i_{t}^{n}\otimes i_{0}^{n}\otimes i_{n-t}^{n}\right)
(M^{\square t}\otimes \rho _{M}^{l}\square M^{\square
n-1-t})(M^{\square t-1}\square \chi _{M}\square M^{\square n-1-t})
\end{eqnarray*}%
an also
\begin{eqnarray*}
&&\sum_{t=1}^{n-1}\sum_{r=1}^{t-1}(i_{r}^{n}\otimes
i_{t-r}^{n}\otimes i_{n-t}^{n})(M^{\square r-1}\square \chi
_{M}\square M^{\square t-1-r}\otimes M^{\square n-t})(M^{\square
t-1}\square \chi _{M}\square
M^{\square n-1-t}) \\
&=&\sum_{t=1}^{n-1}\sum_{r=1}^{t-1}(i_{r}^{n}\otimes
i_{t-r}^{n}\otimes i_{n-t}^{n})(M^{\square r}\otimes M^{\square
t-1-r}\square \chi _{M}\square M^{\square n-1-t})(M^{\square
r-1}\square \chi _{M}\square M^{\square n-1-r})
\\
&=&\sum_{t=1}^{n-1}\sum_{\substack{ r+j=t  \\
r,j>0}}(i_{r}^{n}\otimes i_{j}^{n}\otimes i_{n-t}^{n})(M^{\square
r}\otimes M^{\square j-1}\square \chi _{M}\square M^{\square
n-1-t})(M^{\square r-1}\square \chi _{M}\square
M^{\square n-1-r}) \\
&=&\sum_{\substack{ r+j+k=n  \\ k,r,j>0}}(i_{r}^{n}\otimes
i_{j}^{n}\otimes i_{k}^{n})(M^{\square r}\otimes M^{\square
j-1}\square \chi _{M}\square M^{\square k-1})(M^{\square
r-1}\square \chi _{M}\square M^{\square n-1-r})
\\
&=&\sum_{t=1}^{n-1}\sum_{\substack{ r+k=t  \\
r,k>0}}(i_{t}^{n}\otimes i_{r}^{n}\otimes i_{k}^{n})(M^{\square
t}\otimes M^{\square r-1}\square \chi _{M}\square M^{\square
k-1})(M^{\square t-1}\square \chi _{M}\square
M^{\square n-1-t}) \\
&=&\sum_{t=1}^{n-1}\sum_{r=1}^{n-t-1}(i_{t}^{n}\otimes
i_{r}^{n}\otimes i_{n-t-r}^{n})(M^{\square t}\otimes M^{\square
r-1}\square \chi _{M}\square M^{\square n-t-1-r})(M^{\square
t-1}\square \chi _{M}\square M^{\square n-1-t})
\end{eqnarray*}%
Then $\zeta ^{n}$ satisfies (\ref{2 cocycle}).
\end{proof}

\begin{theorem}\label{teo: T as limit}Let $(C,\Delta ,\varepsilon )$ be a
coalgebra in a cocomplete abelian monoidal category $\M$ with left
exact tensor functors and let $(M,\rho _{M}^{r},\rho _{M}^{l})$ be
a $C$-bicomodule.
  $ (T^c_C(M),(\sigma_i)_{i\in \N})$ carries a
natural coalgebra structure that makes it the direct limit of
$((C^{i}(M))_{i\in \mathbb{N}},(\sigma _{i}^j)_{i,j\in
\mathbb{N}})$ as a direct system of coalgebras.
\end{theorem}

\begin{proof}By Proposition \ref{Pro: cotenso coalgebra}, for any $n\in %
\N\setminus \{0\}$, the canonical inclusion
$\sigma_{n-1}^n:C^{n-1}(M)\rightarrow C^{n}(M)$ is a Hochschild
extension of $C^{n-1}(M)$ with cokernel $M^{\cot n-1}$. In
particular $\sigma_{n-1}^n$ is a coalgebra homomorphism for any
$n\in \N.$ Then, by \ref{claim 4.2}, $\sigma_{m}^n$ is a coalgebra
homomorphism for any $m,n\in \N.$ Now, in view of Corollary
\ref{pro: T=lim c^n} and Proposition \ref{pro: colimit=coalg},
$T^c_C(M)$ carries a natural coalgebra structure that makes it the
direct limit of $((C^{i}(M))_{i\in \mathbb{N}},(\sigma
_{i}^j)_{i,j\in \mathbb{N}})$ as a direct system of coalgebras.
\end{proof}

\begin{lemma}
\label{lem: p_n}Let $(C,\Delta ,\varepsilon )$ be a coalgebra in a
cocomplete abelian monoidal category $\M$ and let $(M,\rho
_{M}^{r},\rho _{M}^{l})$ be a $C$-bicomodule. Let
$T:=T_{C}^{c}(M).$ Then, for any $m,n\geq 1$ the following
relations hold true:
\begin{eqnarray}
(p_{m}\otimes p_{n})\Delta _{T} &=&(M^{\square m-1}\square \chi
_{M}\square M^{\square n-1})p_{m+n}\text{, for any }m,n\geq 1;
\label{formula lem: p_n}
\\
(p_{m}\otimes p_{0})\Delta _{T} &=&(M^{\square m-1}\square \rho
_{M}^{r})p_{m}\text{, for any }m\geq 1;  \label{formula lem: p_1} \\
(p_{0}\otimes p_{n})\Delta _{T} &=&(\rho _{M}^{l}\square
M^{\square
n-1})p_{n}\text{, for any }n\geq 1;  \label{formula lem: p_1bis} \\
(p_{0}\otimes p_{0})\Delta _{T} &=&\Delta _{C}p_{0}.
\label{formula lem: p_0}
\end{eqnarray}
\end{lemma}

\begin{proof}
By construction, the comultiplication $\Delta _{T}$ of $T$ is
uniquely defined by the following relation
\begin{equation*}
\Delta _{T}\sigma _{i}=(\sigma _{i}\otimes \sigma _{i})\Delta
(i),\text{ for any }i\in \mathbb{N},
\end{equation*}
so that, for any $m,n,k\geq 0$ we have
\begin{eqnarray*}
(p_{m}\otimes p_{n})\Delta _{T}i_{k} &=&(p_{m}\otimes p_{n})\Delta
_{T}\sigma _{k+1}i_{k}^{k+1} \\
&=&(p_{m}\otimes p_{n})(\sigma _{k+1}\otimes \sigma _{k+1})\Delta
(k+1)i_{k}^{k+1}=(p_{k+1}^{m}\otimes p_{k+1}^{n})\Delta
(k+1)i_{k}^{k+1}.
\end{eqnarray*}
Recall that:
\begin{eqnarray*}
\Delta (k+1)i_{0}^{k+1} &=&\left( i_{0}^{k+1}\otimes
i_{0}^{k+1}\right)
\Delta _{C} \\
\Delta (k+1)i_{1}^{k+1} &=&\left( i_{1}^{k+1}\otimes
i_{0}^{k+1}\right) \rho
_{M}^{r}+\left( i_{0}^{k+1}\otimes i_{1}^{k+1}\right) \rho _{M}^{l} \\
\text{ and for }2 &\leq &t \\
\Delta (k+1)i_{t}^{k+1} &=&\left( i_{t}^{k+1}\otimes
i_{0}^{k+1}\right) (M^{\square t-1}\square \rho _{M}^{r})+\left(
i_{0}^{k+1}\otimes
i_{t}^{k+1}\right) (\rho _{M}^{l}\square M^{\square t-1})+ \\
&&+\sum_{r=1}^{t-1}(i_{r}^{k+1}\otimes i_{t-r}^{k+1})(M^{\square
r-1}\square \chi _{M}\square M^{\square t-1-r}).
\end{eqnarray*}
Then, for $k=0,1$ we have respecively
\begin{eqnarray*}
(p_{m}\otimes p_{n})\Delta _{T}i_{0} &=&(p_{1}^{m}\otimes
p_{1}^{n})\Delta (1)i_{0}^{1}=\left\{
\begin{tabular}{l}
$0,$ for $m\geq 1$ or $n\geq 1;$ \\
$\Delta (1)i_{0}^{1},$ for $m=n=0;$%
\end{tabular}
\right.  \\
(p_{m}\otimes p_{n})\Delta _{T}i_{1} &=&(p_{2}^{m}\otimes
p_{2}^{n})\Delta
(2)i_{1}^{2} \\
&=&(p_{2}^{m}\otimes p_{2}^{n})[\left( i_{1}^{2}\otimes
i_{0}^{2}\right) \rho _{M}^{r}+\left( i_{0}^{2}\otimes
i_{1}^{2}\right) \rho _{M}^{l}]=\left\{
\begin{tabular}{l}
$0,$ for $m\geq 1$ and $n\geq 1;$ \\
$\delta _{m,1}\rho _{M}^{r},$ for $n=0;$ \\
$\delta _{n,1}\rho _{M}^{l},$ for $m=0;$
\end{tabular}
\right.
\end{eqnarray*}
while, for $k\geq 2,$we have:
\begin{eqnarray*}
&&(p_{m}\otimes p_{n})\Delta _{T}i_{k} \\
&=&(p_{k+1}^{m}\otimes p_{k+1}^{n})\Delta (k+1)i_{k}^{k+1} \\
&=&(p_{k+1}^{m}\otimes p_{k+1}^{n})\left[
\begin{array}{c}
\left( i_{k}^{k+1}\otimes i_{0}^{k+1}\right) (M^{\square
k-1}\square \rho _{M}^{r})+\left( i_{0}^{k+1}\otimes
i_{k}^{k+1}\right) (\rho _{M}^{l}\square
M^{\square k-1}) \\
+\sum_{r=1}^{k-1}(i_{r}^{k+1}\otimes i_{k-r}^{k+1})(M^{\square
r-1}\square
\chi _{M}\square M^{\square k-1-r})%
\end{array}%
\right]  \\
&=&\left\{
\begin{tabular}{l}
$\delta _{n,k-m}(M^{\square m-1}\square \chi _{M}\square
M^{\square n-1}),$
for $m\geq 1$ and $n\geq 1;$ \\
$\delta _{0,k-m}(M^{\square k-1}\square \rho _{M}^{r}),$ for $n=0;$ \\
$\delta _{n,k}(\rho _{M}^{l}\square M^{\square k-1}),$ for $m=0;$%
\end{tabular}%
\ \right.
\end{eqnarray*}
so that, for any $k\geq 0,$ we get:
\begin{eqnarray*}
(p_{m}\otimes p_{n})\Delta _{T}i_{k} &=&\delta _{n+m,k}(M^{\square
m-1}\square \chi _{M}\square M^{\square n-1})=(M^{\square
m-1}\square \chi
_{M}\square M^{\square n-1})p_{m+n}i_{k}\text{, for any }m,n\geq 1; \\
(p_{m}\otimes p_{0})\Delta _{T}i_{k} &=&\delta _{m,k}(M^{\square
k-1}\square \rho _{M}^{r})=(M^{\square m-1}\square \rho
_{M}^{r})p_{m}i_{k}\text{, for
any }m\geq 1; \\
(p_{0}\otimes p_{n})\Delta _{T}i_{k} &=&\delta _{n,k}(\rho
_{M}^{l}\square M^{\square k-1})=(\rho _{M}^{l}\square M^{\square
n-1})p_{n}i_{k}\text{, for
any }n\geq 1; \\
(p_{0}\otimes p_{0})\Delta _{T}i_{k} &=&\delta _{0,k}\Delta
(1)=\Delta _{C}p_{0}i_{k}.
\end{eqnarray*}
Therefore, we conclude.
\end{proof}

\begin{proposition}\label{pro: ciocco}
Let $(C,\Delta ,\varepsilon )$ be a coalgebra in a cocomplete
abelian monoidal
category $\M$  and let $%
(M,\rho _{M}^{r},\rho _{M}^{l})$ be a $C$-bicomodule. Let $%
T:=T_{C}^{c}(M).$ Let $E$ be a coalgebra and let $\alpha
:E\rightarrow T$ and $\beta :E\rightarrow T$ be coalgebra
homomorphisms. If $p_{1}\alpha =p_{1}\beta ,$ then $p_{n}\alpha
=p_{n}\beta $ for any $n\geq 1.$
\end{proposition}

\begin{proof}
Let us prove it by induction on $n\geq 1,$ the case $n=1$ being
true by assumption. Thus let $n\geq 2$ be such that $p_{n}\alpha
=p_{n}\beta .$ Then
\begin{align*}
 (\chi _{M}\square M^{\square n-1})p_{n+1}\alpha &\overset{(\ref{formula lem: p_n})}{
=}(p_{1}\otimes p_{n})\Delta _{T}\alpha  \\
& =(p_{1}\alpha \otimes p_{n}\alpha )\Delta _{T} \\
& =(p_{1}\beta \otimes p_{n}\beta )\Delta _{T} \\
& =(p_{1}\otimes p_{n})\Delta _{T}\beta \overset{(\ref{formula
lem: p_n})}{=}(\chi _{M}\square M^{\square n-1})p_{n+1}\beta .
\end{align*}
Since $\chi _{M}$ is a monomorphism, then, by left exactness of
the tensor
functors, $\chi _{M}\square M^{\square n-1}$ is a monomorphism too, so that $%
p_{n+1}\alpha =p_{n+1}\beta .$
\end{proof}

Our next aim is to prove the universal property of $T^c_C(M)$.

\begin{lemma}
\label{lem w_E=w_D tilde}Let $(C,\Delta ,\varepsilon )$ be a
coalgebra in a cocomplete abelian monoidal category $\M$ and let
$(M,\rho _{M}^{r},\rho _{M}^{l})$ be a $C$-bicomodule. Let $\delta
:D\rightarrow E$ be a monomorphism which is a morphism of
coalgebras such that the canonical morphism $\widetilde{\delta }:\widetilde{D%
}\rightarrow E$ of Notation \ref{notation tilde} is a
monomorphism. Then we have
\begin{equation*}
D^{\wedge _{\widetilde{D}}^{n}}=(D^{\wedge _{E}^{n}},\xi _{n}).
\end{equation*}
\end{lemma}

\begin{proof}
Denote by $(L,p)$ the cokernel of $\delta :D\rightarrow E$ in $\M$ and by $(%
\widetilde{L},\widetilde{p})$ the cokernel of $\xi
_{1}:D\rightarrow \widetilde{D}$ in $\M.$ Since
$p\widetilde{\delta }\xi _{1}=p\delta =0,$ by
the universal property of the cokernel, there exists a unique morphism $%
\lambda :\widetilde{L}\rightarrow L$ such that
\begin{equation*}
\lambda \widetilde{p}=p\widetilde{\delta }.
\end{equation*}
Since $\widetilde{\delta }$ is a monomorphism, so is $\lambda .$
Therefore, we have that
\begin{equation*}
D^{\wedge _{\widetilde{D}}^{n}}:=\ker (\widetilde{p}^{\otimes n}\Delta _{%
\widetilde{D}}^{n-1})=\ker (\lambda ^{\otimes
n}\widetilde{p}^{\otimes
n}\Delta _{\widetilde{D}}^{n-1})=\ker (p^{\otimes n}\widetilde{\delta }%
^{\otimes n}\Delta _{\widetilde{D}}^{n-1})=\ker (p^{\otimes
n}\Delta _{E}^{n-1}\widetilde{\delta }).
\end{equation*}
We conclude by observing that, since, by definition, $(D^{\wedge
_{E}^{n}},\delta _{n}=\widetilde{\delta }\xi _{n}):=\ker
(p^{\otimes n}\Delta _{E}^{n-1})$ where $\widetilde{\delta }$ is a
monomorphism, the following relation holds true
\begin{equation*}
(D^{\wedge _{E}^{n}},\xi _{n})=\ker (p^{\otimes n}\Delta _{E}^{n-1}%
\widetilde{\delta }).
\end{equation*}
\end{proof}

\begin{theorem}
\label{teo: pre univ property of cotensor coalgebra} Let
$(C,\Delta
,\varepsilon )$ be a coalgebra in a cocomplete abelian monoidal category $%
\mathcal{M}$ and let $(M,\rho _{M}^{r},\rho _{M}^{l})$ be a
$C$-bicomodule. Let $\delta :D\rightarrow E$ be a monomorphism
which is a homomorphism of
coalgebras such that the canonical morphism $\widetilde{\delta }:\widetilde{D%
}\rightarrow E$ of Notation \ref{notation tilde} is a monomorphism. Let $%
f_{C}:\widetilde{D}\rightarrow C$ be a coalgebra homomorphism and let $f_{M}:%
\widetilde{D}\rightarrow M$ be a morphism of $C$-bicomodules such that $%
f_{M}\xi _{1}=0$, where $\widetilde{D}$ is a $C$-bicomodule via
$f_{C}.$ Then there is a unique morphism
$f:\widetilde{D}\rightarrow T_{C}^{c}(M)$ such that
\begin{equation*}
f\xi _{n}=\sigma _{n}f_{n},\text{ for any }n\in \mathbb{N},
\end{equation*}%
where
\begin{equation}
f_{n}=\sum_{t=0}^{n}i_{t}^{n}f_{M}^{\square t}\overline{\Delta }_{\widetilde{%
D}}^{t-1}\xi _{n}  \label{formula: f_n}
\end{equation}%
and $\overline{\Delta
}_{\widetilde{D}}^{n}:\widetilde{D}\rightarrow
\widetilde{D}^{\square n+1}$ is the $n^{\text{th}}$-iteration of $\overline{%
\Delta }_{\widetilde{D}}$ $(\overline{\Delta }_{\widetilde{D}}^{-1}=f_{C},%
\overline{\Delta }_{\widetilde{D}}^{0}=\text{Id}_{\widetilde{D}},\overline{%
\Delta }_{\widetilde{D}}^{1}=\overline{\Delta }_{\widetilde{D}}:\widetilde{D}%
\rightarrow \widetilde{D}\square \widetilde{D}).$ \newline
Moreover:

1) $f$ is a coalgebra homomorphism;

2) $p_{0}f=f_{C}$ and $p_{1}f=f_{M}$, where
$p_{n}:T_{C}^{c}(M)\rightarrow M^{\square n}$ denotes the
canonical projection. \newline Furthermore, any coalgebra
homomorphism $f:\widetilde{D}\rightarrow T_{C}^{c}(M)$ that
fulfils 2) satisfies the following relation:
\begin{equation}
p_{k}f=f_{M}^{\square k}\overline{\Delta
}_{\widetilde{D}}^{k-1}\text{ for any }k\in \mathbb{N}.
\label{relation f}
\end{equation}
\end{theorem}

\begin{proof}
Set $T:=T_{C}^{c}(M)$. Following (\ref{def of delta_n}), denote by
$(L,p)$ the cokernel of $\delta :D\rightarrow E$ in $\mathcal{M}$
and let
\begin{equation*}
(D^{\wedge _{E}^{n}},\delta _{n}):=\mathrm{Ker\,}(p^{\otimes
n}\Delta _{E}^{n-1})
\end{equation*}%
for any $n\in \mathbb{N}\setminus \{0\}$, where $\Delta
^{n}:E\rightarrow
E^{\otimes {n+1}}$ is the $n^{\text{th}}$ iterated comultiplication of $E$ $%
(\Delta ^{0}=\mathrm{Id}_{E},\Delta ^{1}:=\Delta _{E}).$ Since
$f_{C}\xi
_{1}:D\rightarrow C$ is a coalgebra homomorphism, then $D$ becomes a $C$%
-bicomodule and $\xi _{1}$ a morphism of $C$-bicomodules. Set
$C^{n}=C^{n}(M) $ and $D^{n}=D^{\wedge _{E}^{n}}$ for any $n\in
\mathbb{N}\setminus \{0\}.$ Define $f_{n}:D^{n}\rightarrow C^{n}$,
for every $n\in \N$, as in (\ref{formula: f_n}). \newline
Note that, for every $n\geq 1,$ $i_{n}^{n}=0$ so that we have%
\begin{equation*}
f_{n}=\sum_{t=0}^{n-1}i_{t}^{n}f_{M}^{\square t}\overline{\Delta }_{%
\widetilde{D}}^{t-1}\xi _{n}.
\end{equation*}%
\newline
Let us prove that $f_{n+1}\xi _{n}^{n+1}=\sigma _{n}^{n+1}f_{n}$ for any $%
n\in \mathbb{N}.$ We have that
\begin{eqnarray*}
f_{n+1}\xi _{n}^{n+1} &=&\sum_{t=0}^{n}i_{t}^{n+1}f_{M}^{\square t}\overline{%
\Delta }_{\widetilde{D}}^{t-1}\xi _{n+1}\xi _{n}^{n+1} \\
&=&\sum_{t=0}^{n}i_{t}^{n+1}f_{M}^{\square t}\overline{\Delta }_{\widetilde{D%
}}^{t-1}\xi _{n} \\
&=&\sum_{t=0}^{n-1}\sigma _{n}^{n+1}i_{t}^{n}f_{M}^{\square t}\overline{%
\Delta }_{\widetilde{D}}^{t-1}\xi _{n}+i_{n}^{n+1}f_{M}^{\square n}\overline{%
\Delta }_{\widetilde{D}}^{n-1}\xi _{n} \\
&=&\sigma _{n}^{n+1}f_{n}+i_{n}^{n+1}f_{M}^{\square n}\overline{\Delta }_{%
\widetilde{D}}^{n-1}\xi _{n}.
\end{eqnarray*}%
Let $(\widetilde{L},\widetilde{p})$ be the cokernel of $\xi
_{1}:D\rightarrow
\widetilde{D}$ in $\mathcal{M}.$ Let $\chi _{\widetilde{L}}:\widetilde{%
L}\square \widetilde{L}\rightarrow \widetilde{L}\otimes
\widetilde{L}$ be the canonical injection. Define $\chi _{\widetilde{L}}^{n}:%
\widetilde{L}^{\square n}\rightarrow \widetilde{L}^{\otimes n}$, for every $n\in \N$, by setting $\chi _{%
\widetilde{L}}^{0}=\mathrm{Id}_{C},$ $\chi _{\widetilde{L}}^{1}=\mathrm{Id}_{%
\widetilde{L}},$ $\chi _{\widetilde{L}}^{2}=\chi _{\widetilde{L}}$
 and $\chi _{\widetilde{L}}^{n}=(\widetilde{L}%
^{\otimes n-2}\otimes \chi _{\widetilde{L}})(\chi _{\widetilde{L}%
}^{n-1}\square \widetilde{L})$ for any $n>2.$ Since the tensor
functors are left exact, $\chi _{\widetilde{L}}^{n}$ is a
monomorphism. By Lemma \ref{lem w_E=w_D tilde}, we have
\begin{equation*}
(D^{\wedge _{E}^{n}},\xi _{n})=D^{\wedge _{\widetilde{D}}^{n}}=\mathrm{Ker\,}%
(\widetilde{p}^{\otimes n}\Delta _{\widetilde{D}}^{n-1})
\end{equation*}%
Thus, for any $n\geq 1,$ we have
\begin{equation*}
\mathrm{Ker\,}(\widetilde{p}^{\otimes n}\Delta _{\widetilde{D}}^{n-1})=%
\mathrm{Ker\,}(\chi _{\widetilde{L}}^{n}\circ \widetilde{p}^{\square n}%
\overline{\Delta }_{\widetilde{D}}^{n-1})=\mathrm{Ker\,}(\widetilde{p}%
^{\square n}\overline{\Delta }_{\widetilde{D}}^{n-1}).
\end{equation*}%
so that we get
\begin{equation}
\widetilde{p}^{\square n}\overline{\Delta
}_{\widetilde{D}}^{n-1}\xi _{n}=0. \label{annullamento delta_n}
\end{equation}%
Now, since $f_{M}\xi _{1}=0,$ there is a unique morphism of $C$%
-bicomodules $\lambda :\widetilde{L}\rightarrow M$ such that
$\lambda \widetilde{p}=f_{M}.$ Thus:
\begin{equation}
f_{M}^{\square n}\overline{\Delta }_{\widetilde{D}}^{n-1}\xi
_{k}=\lambda
^{\square n}(\widetilde{p}^{\square n}\overline{\Delta }_{\widetilde{D}%
}^{n-1}\xi _{n})\xi _{k}^{n}\overset{(\ref{annullamento
delta_n})}{=}0,\text{ for any }k\leq n.  \label{formula boh}
\end{equation}%
We conclude that $f_{n+1}\xi _{n}^{n+1}=\sigma _{n}^{n+1}f_{n}$ for any $%
n\in \mathbb{N}$ so that $(\sigma _{n}f_{n}:D^{n}\rightarrow
T)_{n}$ is a compatible family of morphisms in $\mathcal{M}$. Thus
there is a unique morphism $f:\widetilde{D}\rightarrow T$ such
that
\begin{equation*}
f\xi _{k}=\sigma _{k}f_{k},\text{ for any }k\in \mathbb{N}.
\end{equation*}%
We have that
\begin{equation}
p_{k}^{n}\sum_{t=0}^{k}i_{t}^{k}f_{M}^{\square t}\overline{\Delta }_{%
\widetilde{D}}^{t-1}\xi _{k}=f_{M}^{\square n}\overline{\Delta }_{\widetilde{%
D}}^{n-1}\xi _{k}  \label{formula mah}
\end{equation}%
for any $0\leq n<k$. Note that, for $k\leq n,$ $p_{k}^{n}=0$ and
the right side of (\ref{formula mah}) is zero by (\ref{formula
boh}). Thus, the relation above holds true for any $k,n\in
\mathbb{N}.$ Then we get:
\begin{equation*}
p_{n}f\xi _{k}=p_{n}\sigma
_{k}f_{k}=p_{k}^{n}f_{k}=p_{k}^{n}\sum_{t=0}^{k}i_{t}^{k}f_{M}^{\square t}%
\overline{\Delta }_{\widetilde{D}}^{t-1}\xi _{k}=f_{M}^{\square n}\overline{%
\Delta }_{\widetilde{D}}^{n-1}\xi _{k},\text{ for any }k,n\in
\mathbb{N}.
\end{equation*}%
We conclude that $p_{n}f=f_{M}^{\square n}\overline{\Delta
}_{E}^{n-1}$, for
any $n\in \mathbb{N}.$ In particular, for $n=0,1$, we get $p_{0}f=\overline{%
\Delta }_{\widetilde{D}}^{-1}=f_{C}$ and $p_{1}f=f_{M}\overline{\Delta }_{%
\widetilde{D}}^{0}=f_{M}.$\newline We have now to prove that $f$
is a coalgebra homomorphism. Let us check that
$f_{n}:D^{n}\rightarrow C^{n}$ is a coalgebra homomorphism for
every $n\in \mathbb{N}.$ \newline For $n=0,$ $f_{n}=0$ and there
is nothing to prove. \newline Assume $n\geq 1.$ By Proposition
\ref{Pro: cotenso coalgebra}, we get:
\begin{eqnarray*}
&&\Delta (n)f_{n} \\
&=&\sum_{t=0}^{n}\Delta (n)i_{t}^{n}f_{M}^{\square t}\overline{\Delta }_{%
\widetilde{D}}^{t-1}\xi _{n} \\
&=&\Delta (n)i_{0}^{n}f_{M}^{\square 0}\overline{\Delta }_{\widetilde{D}%
}^{-1}\xi _{n}+\Delta (n)i_{1}^{n}f_{M}^{\square 1}\overline{\Delta }_{%
\widetilde{D}}^{0}\xi _{n}+\sum_{t=2}^{n}\Delta (n)i_{t}^{n}f_{M}^{\square t}%
\overline{\Delta }_{\widetilde{D}}^{t-1}\xi _{n} \\
&=&\Delta (n)i_{0}^{n}f_{C}\xi _{n}+\Delta (n)i_{1}^{n}f_{M}\xi
_{n}+\sum_{t=2}^{n}\Delta (n)i_{t}^{n}f_{M}^{\square t}\overline{\Delta }_{%
\widetilde{D}}^{t-1}\xi _{n} \\
&=&(i_{0}^{n}\otimes i_{0}^{n})\Delta _{C}f_{C}\xi _{n}+\left[
(i_{1}^{n}\otimes i_{0}^{n})\rho _{M}^{r}+(i_{0}^{n}\otimes
i_{1}^{n})\rho
_{M}^{l}\right] f_{M}\xi _{n}+ \\
&&+\sum_{t=2}^{n}\left[
\begin{array}{c}
\left( i_{t}^{n}\otimes i_{0}^{n}\right) (M^{\square t-1}\square
\rho _{M}^{r})+\left( i_{0}^{n}\otimes i_{t}^{n}\right) (\rho
_{M}^{l}\square
M^{\square t-1})+ \\
\sum_{r=1}^{t-1}(i_{r}^{n}\otimes i_{t-r}^{n})(M^{\square
r-1}\square \chi
_{M}\square M^{\square t-1-r})%
\end{array}%
\right] f_{M}^{\square t}\overline{\Delta }_{\widetilde{D}}^{t-1}\xi _{n} \\
&=&(i_{0}^{n}f_{C}\otimes i_{0}^{n}f_{C})\Delta _{\widetilde{D}}\xi _{n}+%
\left[ (i_{1}^{n}\otimes i_{0}^{n})\rho _{M}^{r}+(i_{0}^{n}\otimes
i_{1}^{n})\rho _{M}^{l}\right] f_{M}\xi _{n}+ \\
&&+\sum_{t=2}^{n}\left[ \left( i_{t}^{n}\otimes i_{0}^{n}\right)
(M^{\square t-1}\square \rho _{M}^{r})+\left( i_{0}^{n}\otimes
i_{t}^{n}\right) (\rho
_{M}^{l}\square M^{\square t-1})\right] f_{M}^{\square t}\overline{\Delta }_{%
\widetilde{D}}^{t-1}\xi _{n}+ \\
&&+\sum_{t=2}^{n}\sum_{r=1}^{t-1}(i_{r}^{n}\otimes
i_{t-r}^{n})(M^{\square
r-1}\square \chi _{M}\square M^{\square t-1-r})f_{M}^{\square t}\overline{%
\Delta }_{\widetilde{D}}^{t-1}\xi _{n} \\
&=&(i_{0}^{n}f_{C}\otimes i_{0}^{n}f_{C})\Delta _{\widetilde{D}}\xi _{n}+ \\
&&+\sum_{t=1}^{n}\left[ \left( i_{t}^{n}\otimes i_{0}^{n}\right)
(M^{\square t-1}\square \rho _{M}^{r})+\left( i_{0}^{n}\otimes
i_{t}^{n}\right) (\rho
_{M}^{l}\square M^{\square t-1})\right] f_{M}^{\square t}\overline{\Delta }_{%
\widetilde{D}}^{t-1}\xi _{n}+ \\
&&+\sum_{t=2}^{n}\sum_{r=1}^{t-1}(i_{r}^{n}\otimes
i_{t-r}^{n})(M^{\square
r-1}\square \chi _{M}\square M^{\square t-1-r})f_{M}^{\square t}\overline{%
\Delta }_{\widetilde{D}}^{t-1}\xi _{n}.
\end{eqnarray*}
Now:
\begin{eqnarray*}
\sum_{t=1}^{n}\left( i_{t}^{n}\otimes i_{0}^{n}\right) (M^{\square
t-1}\square \rho _{M}^{r})f_{M}^{\square t}\overline{\Delta }_{\widetilde{D}%
}^{t-1}\xi _{n} &=&\sum_{t=1}^{n}\left( i_{t}^{n}\otimes
i_{0}^{n}\right) (f_{M}^{\square t-1}\square \rho
_{M}^{r}f_{M})\overline{\Delta }_{\widetilde{D}}^{t-1}\xi
_{n} \\
&=&\sum_{t=1}^{n}\left( i_{t}^{n}\otimes i_{0}^{n}\right)
[f_{M}^{\square
t-1}\square (f_{M}\otimes C)\rho _{\widetilde{D}}^{r}]\overline{\Delta }_{%
\widetilde{D}}^{t-1}\xi _{n} \\
&=&\sum_{t=1}^{n}\left( i_{t}^{n}\otimes i_{0}^{n}\right)
(f_{M}^{\square
t-1}\square f_{M}\otimes C)(\widetilde{D}^{\square t-1}\square \rho _{%
\widetilde{D}}^{r})\overline{\Delta }_{\widetilde{D}}^{t-1}\xi _{n} \\
&=&\sum_{t=1}^{n}\left( i_{t}^{n}\otimes i_{0}^{n}\right)
(f_{M}^{\square t}\otimes C)[\widetilde{D}^{\square t-1}\square
(\widetilde{D}\otimes f_{C})\Delta
_{\widetilde{D}}]\overline{\Delta }_{\widetilde{D}}^{t-1}\xi
_{n} \\
&=&\sum_{t=1}^{n}\left( i_{t}^{n}\otimes i_{0}^{n}\right)
(f_{M}^{\square
t}\otimes f_{C})[\widetilde{D}^{\square t-1}\square \Delta _{\widetilde{D}}]%
\overline{\Delta }_{\widetilde{D}}^{t-1}\xi _{n} \\
&=&\sum_{t=1}^{n}\left( i_{t}^{n}\otimes i_{0}^{n}\right)
(f_{M}^{\square
t}\otimes f_{C})[\overline{\Delta }_{\widetilde{D}}^{t-1}\otimes \widetilde{D%
}]\Delta _{\widetilde{D}}\xi _{n} \\
&=&\sum_{t=1}^{n}\left( i_{t}^{n}f_{M}^{\square t}\overline{\Delta }_{%
\widetilde{D}}^{t-1}\otimes i_{0}^{n}f_{C}\right) \Delta
_{\widetilde{D}}\xi _{n}
\end{eqnarray*}%
Analogously one gets
\begin{equation*}
\sum_{t=1}^{n}\left( i_{0}^{n}\otimes i_{t}^{n}\right) (\rho
_{M}^{l}\square M^{\square t-1})f_{M}^{\square t}\overline{\Delta
}_{\widetilde{D}}^{t-1}\xi
_{n}=\sum_{t=1}^{n}\left( i_{0}^{n}f_{C}\otimes i_{t}^{n}f_{M}^{\square t}%
\overline{\Delta }_{\widetilde{D}}^{t-1}\right) \Delta
_{\widetilde{D}}\xi _{n}
\end{equation*}%
Moreover we have:
\begin{eqnarray*}
&&\sum_{t=2}^{n}\sum_{r=1}^{t-1}(i_{r}^{n}\otimes
i_{t-r}^{n})(M^{\square
r-1}\square \chi _{M}\square M^{\square t-1-r})f_{M}^{\square t}\overline{%
\Delta }_{\widetilde{D}}^{t-1}\xi _{n} \\
&=&\sum_{t=2}^{n}\sum_{r=1}^{t-1}(i_{r}^{n}\otimes
i_{t-r}^{n})(f_{M}^{\square r}\otimes f_{M}^{\square t-r})(\widetilde{D}%
^{\square r-1}\square \chi _{\widetilde{D}}\square
\widetilde{D}^{\square
t-1-r})\overline{\Delta }_{\widetilde{D}}^{t-1}\xi _{n} \\
&=&\sum_{t=2}^{n}\sum_{r=1}^{t-1}(i_{r}^{n}f_{M}^{\square
r}\otimes
i_{t-r}^{n}f_{M}^{\square t-r})(\overline{\Delta }_{\widetilde{D}%
}^{r-1}\otimes \overline{\Delta }_{\widetilde{D}}^{t-1-r})\Delta _{%
\widetilde{D}}\xi _{n} \\
&=&\sum_{t=2}^{n}\sum_{r=1}^{t-1}(i_{r}^{n}f_{M}^{\square
r}\overline{\Delta
}_{\widetilde{D}}^{r-1}\otimes i_{t-r}^{n}f_{M}^{\square t-r}\overline{%
\Delta }_{\widetilde{D}}^{t-1-r})\Delta _{\widetilde{D}}\xi _{n}
\end{eqnarray*}%
So we get
\begin{eqnarray*}
\Delta (n)f_{n} &=&(i_{0}^{n}f_{C}\otimes i_{0}^{n}f_{C})\Delta _{\widetilde{%
D}}\xi _{n}+ \\
&&+\sum_{t=1}^{n}\left( i_{t}^{n}f_{M}^{\square t}\overline{\Delta }_{%
\widetilde{D}}^{t-1}\otimes i_{0}^{n}f_{C}\right) \Delta
_{\widetilde{D}}\xi
_{n}+\sum_{t=1}^{n}\left( i_{0}^{n}f_{C}\otimes i_{t}^{n}f_{M}^{\square t}%
\overline{\Delta }_{\widetilde{D}}^{t-1}\right) \Delta
_{\widetilde{D}}\xi
_{n}+ \\
&&+\sum_{t=2}^{n}\sum_{r=1}^{t-1}(i_{r}^{n}f_{M}^{\square
r}\overline{\Delta
}_{\widetilde{D}}^{r-1}\otimes i_{t-r}^{n}f_{M}^{\square t-r}\overline{%
\Delta }_{\widetilde{D}}^{t-1-r})\Delta _{\widetilde{D}}\xi _{n}
\end{eqnarray*}%
On the other hand we have
\begin{eqnarray*}
&&(f_{n}\otimes f_{n})\Delta _{D^{n}} \\
&=&(\sum_{t=0}^{n}i_{t}^{n}f_{M}^{\square t}\overline{\Delta }_{\widetilde{D}%
}^{t-1}\xi _{n}\otimes \sum_{k=0}^{n}i_{k}^{n}f_{M}^{\square k}\overline{%
\Delta }_{\widetilde{D}}^{k-1}\xi _{n})\Delta _{D^{n}} \\
&=&(\sum_{t=0}^{n}i_{t}^{n}f_{M}^{\square t}\overline{\Delta }_{\widetilde{D}%
}^{t-1}\otimes \sum_{k=0}^{n}i_{k}^{n}f_{M}^{\square k}\overline{\Delta }_{%
\widetilde{D}}^{k-1})\Delta _{\widetilde{D}}\xi _{n} \\
&=&(i_{0}^{n}f_{M}^{\square 0}\overline{\Delta
}_{\widetilde{D}}^{-1}\otimes
i_{0}^{n}f_{M}^{\square 0}\overline{\Delta }_{\widetilde{D}}^{-1})\Delta _{%
\widetilde{D}}\xi _{n}+ \\
&&+\sum_{t=1}^{n}(i_{t}^{n}f_{M}^{\square t}\overline{\Delta }_{\widetilde{D}%
}^{t-1}\otimes i_{0}^{n}f_{M}^{\square 0}\overline{\Delta }_{\widetilde{D}%
}^{-1})\Delta _{\widetilde{D}}\xi
_{n}+\sum_{k=1}^{n}(i_{0}^{n}f_{M}^{\square 0}\overline{\Delta }_{\widetilde{%
D}}^{0-1}\otimes i_{k}^{n}f_{M}^{\square k}\overline{\Delta }_{\widetilde{D}%
}^{k-1})\Delta _{\widetilde{D}}\xi _{n} \\
&&+(\sum_{t=1}^{n}i_{t}^{n}f_{M}^{\square t}\overline{\Delta }_{\widetilde{D}%
}^{t-1}\otimes \sum_{k=1}^{n}i_{k}^{n}f_{M}^{\square k}\overline{\Delta }_{%
\widetilde{D}}^{k-1})\Delta _{\widetilde{D}}\xi _{n} \\
&=&(i_{0}^{n}f_{C}\otimes i_{0}^{n}f_{C})\Delta
_{\widetilde{D}}\xi
_{n}+\sum_{t=1}^{n}(i_{t}^{n}f_{M}^{\square t}\overline{\Delta }_{\widetilde{%
D}}^{t-1}\otimes i_{0}^{n}f_{C})\Delta _{\widetilde{D}}\xi
_{n}+\sum_{k=1}^{n}(i_{0}^{n}f_{C}\otimes i_{k}^{n}f_{M}^{\square k}%
\overline{\Delta }_{\widetilde{D}}^{k-1})\Delta _{\widetilde{D}}\xi _{n}+ \\
&&+\sum_{t=1}^{n}\sum_{k=1}^{n}(i_{t}^{n}f_{M}^{\square t}\overline{\Delta }%
_{\widetilde{D}}^{t-1}\otimes i_{k}^{n}f_{M}^{\square k}\overline{\Delta }_{%
\widetilde{D}}^{k-1})\Delta _{\widetilde{D}}\xi _{n} \\
&=&\Delta (n)f_{n}-\sum_{t=2}^{n}\sum_{r=1}^{t-1}(i_{r}^{n}f_{M}^{\square r}%
\overline{\Delta }_{\widetilde{D}}^{r-1}\otimes
i_{t-r}^{n}f_{M}^{\square t-r}\overline{\Delta
}_{\widetilde{D}}^{t-1-r})\Delta _{\widetilde{D}}\xi
_{n}+ \\
&&+\sum_{t=1}^{n}\sum_{k=1}^{n}(i_{t}^{n}f_{M}^{\square t}\overline{\Delta }%
_{\widetilde{D}}^{t-1}\otimes i_{k}^{n}f_{M}^{\square k}\overline{\Delta }_{%
\widetilde{D}}^{k-1})\Delta _{\widetilde{D}}\xi _{n}.
\end{eqnarray*}
But
\begin{align*}
& (i_{t}^{n}f_{M}^{\square t}\overline{\Delta
}_{\widetilde{D}}^{t-1}\otimes
i_{k}^{n}f_{M}^{\square k}\overline{\Delta }_{\widetilde{D}}^{k-1})\Delta _{%
\widetilde{D}}\xi _{n} \\
& =(i_{t}^{n}f_{M}^{\square t}\otimes i_{k}^{n}f_{M}^{\square k})(\overline{%
\Delta }_{\widetilde{D}}^{t-1}\otimes \overline{\Delta }_{\widetilde{D}%
}^{k-1})\Delta _{\widetilde{D}}\xi _{n} \\
& =(i_{t}^{n}f_{M}^{\square t}\otimes i_{k}^{n}f_{M}^{\square k})(\widetilde{%
D}^{\square t-1}\square \chi _{\widetilde{D}}\square
\widetilde{D}^{\square
k-1})\overline{\Delta }_{\widetilde{D}}^{t+k-1}\xi _{n} \\
& =(i_{t}^{n}\otimes i_{k}^{n})(f_{M}^{\square t}\otimes f_{M}^{\square k})(%
\widetilde{D}^{\square t-1}\square \chi _{\widetilde{D}}\square \widetilde{D}%
^{\square k-1})\overline{\Delta }_{\widetilde{D}}^{t+k-1}\xi _{n} \\
& =(i_{t}^{n}\otimes i_{k}^{n})(M^{\square t-1}\square \chi
_{M}\square
M^{\square k-1})f_{M}^{\square t+k}\overline{\Delta }_{\widetilde{D}%
}^{t+k-1}\xi _{n}.
\end{align*}%
By $(\ref{formula boh}),$ the last term is zero whenever $t+k>n,$
so that:
\begin{eqnarray*}
\sum_{t=1}^{n}\sum_{k=1}^{n}(i_{t}^{n}f_{M}^{\square t}\overline{\Delta }_{%
\widetilde{D}}^{t-1}\otimes i_{k}^{n}f_{M}^{\square k}\overline{\Delta }_{%
\widetilde{D}}^{k-1})\Delta _{\widetilde{D}}\xi _{n} &=&\sum_{\substack{ %
1\leq t,k\leq n \\ t+k\leq n}}(i_{t}^{n}f_{M}^{\square t}\overline{\Delta }_{%
\widetilde{D}}^{t-1}\otimes i_{k}^{n}f_{M}^{\square k}\overline{\Delta }_{%
\widetilde{D}}^{k-1})\Delta _{\widetilde{D}}\xi _{n} \\
&=&\sum_{t=2}^{n}\sum_{r=1}^{t-1}(i_{r}^{n}f_{M}^{\square
r}\overline{\Delta
}_{\widetilde{D}}^{r-1}\otimes i_{t-r}^{n}f_{M}^{\square t-r}\overline{%
\Delta }_{\widetilde{D}}^{t-1-r})\Delta _{\widetilde{D}}\xi _{n},
\end{eqnarray*}%
and hence
\begin{equation*}
(f_{n}\otimes f_{n})\Delta _{D^{n}}=\Delta (n)f_{n}.
\end{equation*}%
Furthermore $\varepsilon (0)f_{0}=0=\varepsilon _{D^{0}},$ while, for every $%
n\geq 1$, we have
\begin{eqnarray*}
\varepsilon (n)f_{n} &=&\sum_{t=0}^{n}\varepsilon
(n)i_{t}^{n}f_{M}^{\square
t}\overline{\Delta }_{\widetilde{D}}^{t-1}\xi _{n} \\
&=&\sum_{t=0}^{n}\varepsilon _{C}\pi _{n}^{1}i_{t}^{n}f_{M}^{\square t}%
\overline{\Delta }_{\widetilde{D}}^{t-1}\xi _{n}=\varepsilon
_{C}\pi
_{n}^{1}i_{0}^{n}f_{M}^{\square 0}\overline{\Delta }_{\widetilde{D}%
}^{0-1}\xi _{n}=\varepsilon _{C}i_{0}^{1}f_{C}\xi _{n}=\varepsilon
_{C}f_{C}\xi _{n}=\varepsilon _{D^{n}}.
\end{eqnarray*}%
We conclude that $f_{n}$ is a coalgebra homomorphism. Now, by construction, $%
f$ is the unique morphism such that $f\xi _{k}=\sigma _{k}f_{k},$ for any $%
k\in \mathbb{N}.$ By Proposition \ref{pro: limit of delta},
$((D^{\wedge _{C}^{i}})_{i\in \mathbb{N}},(\xi _{i}^{j})_{i,j\in
\mathbb{N}})$ is a direct system in $\mathcal{M}$ whose direct
limit $\widetilde{D}$ carries a natural coalgebra structure that
makes it the direct limit of $((D^{\wedge _{C}^{i}})_{i\in
\mathbb{N}},(\xi _{i}^{j})_{i,j\in \mathbb{N}})$ as a direct
system of coalgebras. Since $\sigma _{k}$ is a coalgebra
homomorphism so is $\sigma _{k}f_{k}$ and hence $f$ is a coalgebra
homomorphism.\newline Assume now that $g:E\rightarrow T$ is
another coalgebra homomorphism such that $p_{0}g=f_{C}$ and
$p_{1}g=f_{M}.$ Then, by Proposition \ref{pro: ciocco}, we have
$p_{n}g=p_{n}f$ for any $n\in \mathbb{N}.$
\end{proof}

\begin{lemma}
\label{lem: complete}Let $(X_{i})_{i\in \N}$ be a family of
objects in a cocomplete and complete abelian category $\M$
satisfying AB5. Let $Y$ be an object in $\M$ and $f:Y\rightarrow
\oplus _{i\in \N}X_{i}$ be a morphism such that
\begin{equation*}
p_{k}f=0\text{ for any }k\in \N,
\end{equation*}
where $p_{k}:\oplus X_{i}\rightarrow X_{k}$ denotes the canonical
projection. Then $f=0.$
\end{lemma}

\begin{proof}
By \cite[Corollary 8.10, page 61]{Po}, $\M$ is a $C_2$-category so
that the conclusion follows by \cite[Proposition, page 54]{Po}.
\end{proof}

\begin{theorem}
\label{coro: univ property of cotensor coalgebra} Let $(C,\Delta
,\varepsilon )$ be a coalgebra in a cocomplete and complete
abelian monoidal category $\M$ satisfying AB5. Let $(M,\rho
_{M}^{r},\rho _{M}^{l})$ be a $C$-bicomodule. Let $\delta
:D\rightarrow E$ be a monomorphism which is a homomorphism of coalgebras. Let $%
f_{C}:\widetilde{D}\rightarrow C$ be a coalgebra homomorphism and let $f_{M}:%
\widetilde{D}\rightarrow M$ be a morphism of $C$-bicomodules such
that $f_{M}\xi _{1}=0$, where $\widetilde{D}$ is a bicomodule via
$f_{C}.$ Then there is a unique coalgebra homomorphism
$f:\widetilde{D}\rightarrow
T_{C}^{c}(M)$ such that $p_{0}f=f_{C}$ and $p_{1}f=f_{M}$, where $%
p_{n}:T_{C}^{c}(M)\rightarrow M^{\square n}$ denotes the canonical
projection.
\begin{equation*}
\begin{diagram}[h=2em,w=3em]
T_{C}^{c}(M)&\rTo^{p_1}&M\\\dTo^{p_0}&\luDotsto^{f}&\uTo>{f_M}&\luTo^0\\C&%
\lTo_{f_C}&\widetilde{D}&\lTo _{\xi_1}&D \end{diagram}
\end{equation*}
\end{theorem}

\begin{proof}
Since $\M$ satisfies AB5, the morphism $\widetilde{\delta
}:\widetilde{D}\rightarrow E$ of Notation \ref{notation tilde} is
a monomorphism, so
that, by applying Theorem \ref{teo: pre univ property of cotensor coalgebra}%
, there is a coalgebra homomorphism $f:\widetilde{D}\rightarrow T_{C}^{c}(M)$ such that $%
p_{0}f=f_{C}$ and $p_{1}f=f_{M}$. Furthermore any such a coalgebra
homomorphism fulfils the following relation
\begin{equation*}
p_{k}f=f_{M}^{\square k}\overline{\Delta
}_{\widetilde{D}}^{k-1}\text{ for any }k\in \N,
\end{equation*}
where $\overline{\Delta
}_{\widetilde{D}}^{n}:\widetilde{D}\rightarrow
\widetilde{D}^{\square n+1}$ is the $n^{\text{th}}$-iteration of $\overline{%
\Delta }_{\widetilde{D}}$ as defined in Theorem \ref{teo: pre univ
property of cotensor coalgebra}.\newline Now, let
$f,g:\widetilde{D}\rightarrow T$ be morphisms such that
\begin{equation*}
p_{k}f=p_{k}g\text{ for any }k\in \N.
\end{equation*}
By Lemma \ref{lem: complete} $f=g.$
\end{proof}

\begin{lemma}
\label{lem: wedge of kernels}Let $E$ be a coalgebra in an abelian
monoidal category $\mathcal{M}$. Let $f:E\rightarrow L $ and
$g:E\rightarrow M$ be morphism in $\mathcal{M}$. Then
\begin{equation*}
\mathrm{Ker\,} (f)\wedge _{E}\mathrm{Ker\,} (g)=Ker[(f\otimes
g)\circ \triangle _{E}].
\end{equation*}
\end{lemma}

\begin{proof}
Let $(X,i_{X})=\mathrm{Ker\,} (f)$ and let
$(Y,i_{Y})=\mathrm{Ker\,} (g).$ Let $p_{X}:E\rightarrow E/X$ and
$p_{Y}:E\rightarrow E/Y$ be the canonical quotient maps. Since
$fi_{X}=0,$ by the universal property of the cokernel,
there exists a unique morphism $\gamma _{X}:E/X\rightarrow L$ such that $%
\gamma _{X}p_{X}=f.$ Moreover, we have $(E/X,p_{X})=\mathrm{coker}(i_{X})=%
\mathrm{coker}(\mathrm{Ker\,} (f))=\mathrm{coim}(f).$ As
$\mathcal{M}$ is
abelian, it is straightforaword to prove that $(E/X,\gamma _{X})=\mathrm{Im}%
(f).$ In particular $\gamma _{X}$ is a monomorphism. Analogously
one gets a monomorphism $\gamma _{Y}:E/Y\rightarrow M$ such that
$\gamma _{Y}p_{Y}=g.$ Since $\mathcal{M}$ has left exact tensor
functors, then $\gamma _{X}\otimes \gamma _{Y}$ is a monomorphism,
so that, by definition, we get:
\begin{equation*}
X\wedge _{E}Y:=Ker[(p _{X}\otimes p _{Y})\circ \triangle
_{E}]=Ker[(\gamma _{X}\otimes \gamma _{Y})(p _{X}\otimes p
_{Y})\circ \triangle _{E}]=Ker[(f\otimes g)\circ \triangle _{E}].
\end{equation*}
\end{proof}

\begin{proposition}
\label{pro: D^2}Let $\delta :D\rightarrow C$ be a monomorphism
which is a coalgebra homomorphism in an abelian monoidal category
$\mathcal{M}$. Then we have
\begin{equation}
D^{\wedge _{C}^{m}}\wedge _{C}D^{\wedge _{C}^{n}}=D^{\wedge
_{C}^{m+n}}. \label{formula 2 pro: D^2}
\end{equation}
\end{proposition}

\begin{proof}
Set $(L,p)=\mathrm{coker}(\sigma ).$ Let
\begin{equation*}
(D^{\wedge _{C}^{n}},\delta _{n}):=\mathrm{Ker\,}(p^{\otimes
n}\Delta _{C}^{n-1})\text{.}
\end{equation*}%
By Lemma \ref{lem: wedge of kernels}, we have:
\begin{equation*}
D^{\wedge _{C}^{m+n}}=\mathrm{Ker\,}[p^{\otimes m+n}\Delta _{C}^{m+n-1}]=%
\mathrm{Ker\,}[(p^{\otimes m}\Delta _{C}^{m-1}\otimes p^{\otimes
n}\Delta _{C}^{n-1})\Delta _{C}]=D^{\wedge _{C}^{m}}\wedge
_{C}D^{\wedge _{C}^{n}}.
\end{equation*}
\end{proof}

\begin{theorem}\label{pro: wedge limit}
Let $(C,\Delta ,\varepsilon )$ be a coalgebra in a cocomplete and
complete abelian monoidal category $\mathcal{M}$ satisfying $AB5.$
Let $(M,\rho
_{M}^{r},\rho _{M}^{l})$ be a $C$-bicomodule. Let $T:=T_{C}^{c}(M)$. Then%
\begin{equation*}
(C^{n}\left( M\right),\sigma_n) =C^{\wedge _{T}^{n}},
\end{equation*}%
for every $n\in
\mathbb{N}
$.
\end{theorem}

\begin{proof}
Let $n\in
\mathbb{N}
$. Let
\begin{equation*}
\lambda _{a}^{n}:M^{\square a}\rightarrow \oplus _{b\geq
n}M^{\square b}
\end{equation*}%
be defined by%
\begin{equation*}
\lambda _{a}^{n}=\left\{
\begin{tabular}{ll}
the canonical injection & if $a\geq n,$ \\
$0$ & otherwise.%
\end{tabular}%
\right.
\end{equation*}%
Define
\begin{equation*}
\nu _{n}:\oplus _{a\geq n}M^{\square a}\rightarrow T
\end{equation*}%
as the codiagonal map of the family $\left( i_{a}\right) _{a\geq
n}$ so that
we have%
\begin{equation*}
\nu _{n}\circ \lambda _{a}^{n}=\left\{
\begin{tabular}{ll}
$i_{a}$ & for every $a\geq n$ \\
$0$ & otherwise.%
\end{tabular}%
\right.
\end{equation*}%
Define
\begin{equation*}
\tau _{n}:T\rightarrow \oplus _{a\geq n}M^{\square a}
\end{equation*}%
as the codiagonal map of the family $\left( \lambda
_{a}^{n}\right) _{a\in
\mathbb{N}
},$ that is%
\begin{equation*}
\tau _{n}i_{a}=\lambda _{a}^{n}\text{, for every }a\in
\mathbb{N}
.
\end{equation*}%
Thus, we have%
\begin{equation*}
\tau _{n}\nu _{n}\lambda _{a}^{n}=\tau _{n}i_{a}=\lambda _{a}^{n}\text{ for }%
a\geq n,
\end{equation*}%
so that%
\begin{equation*}
\tau _{n}\nu _{n}=\mathrm{Id}_{\oplus _{a\geq n}M^{\square
a}}\text{ for every }n\in
\mathbb{N}
.
\end{equation*}%
Let $C^{n}=C^{n}\left( M\right) $ for every $n\in
\mathbb{N}
.$ Let us prove the following sequence%
\begin{equation*}
0\rightarrow C^{n}\overset{\sigma _{n}}{\longrightarrow }T\overset{\tau _{n}}%
{\longrightarrow }\oplus _{a\geq n}M^{\square a}\rightarrow 0
\end{equation*}%
is exact. We check that $\left( \oplus _{a\geq n}M^{\square
a},\tau _{n}\right) =\mathrm{Coker}\left( \sigma _{n}\right)
$.\newline Since $\tau _{n}\nu _{n}=\mathrm{Id},$ it is clear that
$\tau _{n}$ is an epimorphism and that $\nu _{n}$ is a
monomorphism.\newline
From%
\begin{equation*}
\tau _{n}\sigma _{n}i_{a}^{n}=\tau _{n}i_{a}=\lambda
_{a}^{n}=0,\text{ for every }0\leq a\leq n-1,
\end{equation*}%
we deduce that%
\begin{equation*}
\tau _{n}\sigma _{n}=0\text{ for every }n\in
\mathbb{N}
.
\end{equation*}%
Let $f:T\rightarrow X$ be a morphism such that $f\sigma _{n}=0$ for every $%
n\in
\mathbb{N}
.$ Thus, for every $0\leq a\leq n-1,$ we have%
\begin{equation*}
fi_{a}=f\sigma _{n}i_{a}^{n}=0
\end{equation*}%
Set%
\begin{equation*}
\overline{f}=f\nu _{n}
\end{equation*}%
and let us prove that $f=\overline{f}\tau _{n}.$ From%
\begin{equation*}
\overline{f}\tau _{n}i_{a}=f\nu _{n}\lambda _{a}^{n}=\left\{
\begin{tabular}{ll}
$fi_{a}$ & for every $a\geq n$ \\
$0$ & otherwise.%
\end{tabular}%
\right.
\end{equation*}%
we deduce that $\overline{f}\tau _{n}i_{a}=fi_{a},$ for every
$a\in
\mathbb{N}
,$ and hence $\overline{f}\tau _{n}=f.$

Let us prove that $C^{n}=C^{\wedge _{T}^{n}},$ for every $n\in
\mathbb{N}
.$ \newline The case $n=0$ is trivial. Let us prove the equality
above for every $n\geq 1 $ by induction on $n.$\newline If $n=1,$
by definition, we have $C^{1}=C=C^{\wedge _{T}^{1}}.$\newline Let
$n\geq 2$ and assume that $C^{\wedge _{T}^{n-1}}=C^{n-1}.$ By
Proposition \ref{pro: D^2} and Lemma \ref{lem: wedge of
kernels}, we have%
\begin{equation*}
C^{\wedge _{T}^{n}}=C^{\wedge _{T}^{n-1}}\wedge _{T}C^{\wedge
_{T}^{1}}=C^{n-1}\wedge _{T}C^{1}=\mathrm{Ker\,}\left( \tau
_{n-1}\right) \wedge _{T}\mathrm{Ker\,}\left( \tau _{1}\right)
=\mathrm{Ker\,}\left[ \left( \tau _{n-1}\otimes \tau _{1}\right)
\Delta _{T}\right]
\end{equation*}
so that%
\begin{equation*}
0\rightarrow C^{\wedge _{T}^{n}}\hookrightarrow T\overset{\left(
\tau _{n-1}\otimes \tau _{1}\right) \Delta _{T}}{\longrightarrow
}\left( \oplus _{a\geq n-1}M^{\square a}\right) \otimes \left(
\oplus _{a\geq 1}M^{\square a}\right)
\end{equation*}%
is an exact sequence. In order to conclude, it is enough to check that the following sequence%
\begin{equation*}
0\rightarrow C^{n}\overset{\sigma _{n}}{\longrightarrow
}T\overset{\left( \tau _{n-1}\otimes \tau _{1}\right) \Delta
_{T}}{\longrightarrow }\left( \oplus _{a\geq n-1}M^{\square
a}\right) \otimes \left( \oplus _{a\geq 1}M^{\square a}\right)
\end{equation*}%
is exact. For every $0\leq a\leq n-1,$ we have%
\begin{equation*}
\left( \tau _{n-1}\otimes \tau _{1}\right) \Delta _{T}\sigma
_{n}i_{a}^{n}=\left( \tau _{n-1}\otimes \tau _{1}\right) \left(
\sigma _{n}\otimes \sigma _{n}\right) \Delta \left( n\right)
i_{a}^{n}
\end{equation*}%
By Proposition \ref{Pro: cotenso coalgebra} we can write
\begin{equation*}
\Delta \left( n\right) i_{a}^{n}=\sum_{r=0}^{a}(i_{r}^{n}\otimes
i_{a-r}^{n})f_{r,a-r}
\end{equation*}%
where $f_{ij}:M^{\square a}\rightarrow M^{\square i}\otimes
M^{\square j}$
are suitable morphisms. Thus we get%
\begin{eqnarray*}
\left( \tau _{n-1}\otimes \tau _{1}\right) \Delta _{T}\sigma
_{n}i_{a}^{n} &=&\sum_{r=0}^{a}\left( \tau _{n-1}\otimes \tau
_{1}\right) \left( \sigma
_{n}\otimes \sigma _{n}\right) (i_{r}^{n}\otimes i_{a-r}^{n})f_{r,a-r} \\
&=&\sum_{r=0}^{a}\left( \tau _{n-1}i_{r}\otimes \tau
_{1}i_{a-r}\right) f_{r,a-r}=\sum_{r=0}^{a}\left( \lambda
_{r}^{n-1}\otimes \lambda _{a-r}^{1}\right) f_{r,a-r}=0
\end{eqnarray*}%
Therefore $ \left( \tau _{n-1}\otimes \tau _{1}\right) \Delta
_{T}\sigma _{n}=0 $, for every $n\in \N.$\\ Let $g:Y\rightarrow T$
be a morphism such that
\begin{equation*}
\left( \tau _{n-1}\otimes \tau _{1}\right) \Delta _{T}g=0.
\end{equation*}%
Now, for every $c\in \N$ and for every $a\geq b,$ we have%
\begin{equation*}
p_{a}\nu _{b}\tau _{b}i_{c}=p_{a}\nu _{b}\lambda _{c}^{b}=\left\{
\begin{tabular}{ll}
$p_{a}i_{c}$ & for every $c\geq b$ \\
$0=p_{a}i_{c}$ & otherwise.%
\end{tabular}%
\right.
\end{equation*}%
so that%
\begin{equation*}
p_{a}\nu _{b}\tau _{b}=p_{a},\text{ for every }a\geq b.
\end{equation*}%
Thus, for every $a\geq n-1$ and $b\geq 1,$ by Lemma \ref{lem: p_n}, we have%
\begin{equation*}
0=\left( p_{a}\nu _{n-1}\tau _{n-1}\otimes p_{b}\nu _{1}\tau
_{1}\right) \Delta _{T}g=\left( p_{a}\otimes p_{b}\right) \Delta
_{T}g{=}\left( M^{\square a-1}\square \chi _{M}\square M^{\square
b}\right) p_{a+b}g.
\end{equation*}%
By left exactness of the tensor functors, $M^{\square a-1}\square
\chi
_{M}\square M^{\square b}$ is a monomorphism so that%
\begin{equation*}
p_{a+b}g=0.
\end{equation*}%
We conclude that%
\begin{equation*}
p_{c}g=0,\text{ for every }c\geq n.
\end{equation*}%
Set%
\begin{equation*}
\overline{g}=\pi _{n}g.
\end{equation*}%
and let us prove that $g=\sigma _{n}\overline{g}.$ By Lemma
\ref{lem:
complete} this is the case if and only if%
\begin{equation*}
p_{a}g=p_{a}\sigma _{n}\overline{g},\text{ for every }a\in
\mathbb{N}
.
\end{equation*}%
We have%
\begin{equation*}
p_{a}\sigma _{n}\overline{g}=p_{n}^{a}\pi _{n}g=\left\{
\begin{tabular}{ll}
$p_{a}g$ & for every $a<n$ \\
$0=p_{a}g$ & otherwise.%
\end{tabular}%
\right.
\end{equation*}
\end{proof}

\section{Technicalities}

The main aim of this section is to proof Theorem \ref{teo alfa
beta sigma} that will be our main tool in the proof of Theorem
\ref{teo Smooth}.

\begin{lemma}
\label{lem: schiacciamento}Let $\mathcal{M}$ be a monoidal
category with left exact direct limits. Let $((A_{i})_{i\in
\mathbb{N}},(\alpha _{i}^{j})_{i,j\in \mathbb{N}})$ and
$((B_{i})_{i\in \mathbb{N}},(\beta _{i}^{j})_{i,j\in \mathbb{N}}$
be direct systems in $\mathcal{M}$, where, for $i\leq j$, $\alpha
_{i}^{j}:A_{i}\rightarrow A_{j}$ and $\beta
_{i}^{j}:B_{i}\rightarrow B_{j}$. Let $(\gamma
_{i}:A_{i}\rightarrow
B_{i})_{i\in \mathbb{N}}$ be a direct system of monomorphisms. Let $%
(A,(\alpha _{i})_{i\in \mathbb{N}})=\underrightarrow{\lim }A_{i}$ and let $%
(\beta _{i}:B_{i}\rightarrow A)_{i\in \mathbb{N}}$ be a compatible
family of monomorphisms such that $\beta _{i}\gamma _{i}=\alpha
_{i}$ for any $i\in
\mathbb{N.}$ Then $(A,(\beta _{i})_{i\in \mathbb{N}})=\underrightarrow{\lim }%
B_{i}.$
\end{lemma}

\begin{proof}
Since direct limits are left exact in $\mathcal{M}$, the canonical morphism $%
\underrightarrow{\lim }\beta _{i}:\underrightarrow{\lim
}B_{i}\rightarrow A$ is a monomorphism. Moreover, since
$\underrightarrow{\lim }\beta _{i}\circ
\underrightarrow{\lim }\gamma _{i}=\underrightarrow{\lim }\alpha _{i}=%
\mathrm{Id}_{A},$ we have that $\underrightarrow{\lim }\beta _{i}$
is also an epimorphism and hence an isomorphism.
\end{proof}

\begin{claim}
Let $\left( \mathcal{M},\otimes ,\mathbf{1}\right) $ be a
cocomplete abelian monoidal category. Recall that a \emph{graded
}coalgebra in $\mathcal{M}$ is a coalgebra $\left( B,\Delta
,\varepsilon \right) $ endowed with a family $\left( B_{i},\beta
_{i}\right) $ of subobjects of $B$, such that
\begin{equation*}
B=\oplus _{i\in \N}B_{i}
\end{equation*}%
and there esists a family $\left( \Delta _{i}\right) _{_{i\in\N}}$
of morphisms
\begin{equation*}
\Delta _{i}:B_{i}\rightarrow \left( B\otimes B\right) _{i}=\oplus
_{a+b=i}\left( B_{a}\otimes B_{b}\right) ,
\end{equation*}%
such that
\begin{equation*}
\Delta \beta _{i}=\left[ \nabla _{a+b=i}\left( \beta _{a}\otimes
\beta _{b}\right) \right] \Delta _{i}
\end{equation*}%
and%
\begin{equation*}
\varepsilon \beta _{i}=0,\text{ for every }i\geq 1.
\end{equation*}%
Here $\nabla _{a+b=i}\left( \beta _{a}\otimes \beta _{b}\right) $
denotes the codiagonal morphism associated to the family $\left(
\beta _{a}\otimes \beta _{b}\right) _{a+b=i}.$\newline In
particular it follows that $\varepsilon $ restricts to a morphism
\begin{equation*}
\varepsilon _{0}:B_{0}\rightarrow \mathbf{1}
\end{equation*}%
such that%
\begin{equation*}
\varepsilon =\varepsilon _{0}\beta _{0}
\end{equation*}%
and $\left( B_{0},\Delta _{0},\varepsilon _{0}\right) $ is a coalgebra in $%
\mathcal{M}.$ Moreover $\beta _{0}$ is a coalgebra homomorphism.
\end{claim}

\begin{proposition}
\label{pro: lim for graded coalg}Let $\left( \mathcal{M},\otimes ,\mathbf{1}%
\right) $ be a cocomplete abelian monoidal category. Let $B=\oplus
_{i\in
\mathbb{N}
}B_{i}$ be a graded coalgebra. Denote by $\left( L,p\right) $ the
cokernel
of $\beta _{0}$ in $\mathcal{M}$. Then%
\begin{equation}
p^{\otimes n+1}\Delta _{B}^{n}\beta _{b}=0,\text{ for every }0\leq
b\leq n. \label{formula lem: c^n in wedge}
\end{equation}%
Moreover%
\begin{equation*}
B=\underrightarrow{\lim }(B_{0}^{\wedge _{B}^{i}})_{i\in
\mathbb{N}}.
\end{equation*}
\end{proposition}

\begin{proof}
Denote by $\beta _{i}:B_{i}\rightarrow B$ the canonical inclusion
and denote by $\tau _{i}:B\rightarrow B_{i}$ the canonical
projection, for every $i\in \N$. Since $\beta _{0}$ is a coalgebra
homomorphism and $\beta _{0}$ is a
monomorphism, we can consider%
\begin{equation*}
(B_{0}^{\wedge _{B}^{n}},\delta _{n}):=\mathrm{Ker\,}(p^{\otimes
n}\Delta _{B}^{n-1})\text{.}
\end{equation*}%
Denote by $\xi _{i}^{j}:B_{0}^{\wedge _{B}^{i}}\rightarrow
B_{0}^{\wedge
_{B}^{j}}$ the canonical inclusion, for every $j\geq i.$ Hence we have%
\begin{equation*}
\delta _{j}\xi _{i}^{j}=\delta _{i}.
\end{equation*}%
In order to prove (\ref{formula lem: c^n in wedge}), we proceed by
induction on $n\geq 0.$ For $n=0,$ then $b=0$ and we have
\begin{equation*}
p^{\otimes n+1}\Delta _{B}^{n}\beta _{b}=p\Delta _{B}^{0}\beta
_{0}=p\beta _{0}=0.
\end{equation*}%
Let $n\geq 1$ and assume $p^{\otimes i+1}\Delta _{B}^{i}\beta
_{j}=0,$ for
every $0\leq j\leq i\leq n-1.$ For every $0\leq c\leq n,$ we have%
\begin{eqnarray*}
p^{\otimes n+1}\Delta _{B}^{n}\beta _{c} &=&\left( p^{\otimes
n-1}\otimes p^{\otimes 2}\right) \left( \Delta _{B}^{n-2}\otimes
\Delta _{B}\right)
\Delta _{B}\beta _{c} \\
&=&\left( p^{\otimes n-1}\Delta _{B}^{n-2}\otimes p^{\otimes
2}\Delta
_{B}\right) \left[ \nabla _{a+b=c}\left( \beta _{a}\otimes \beta _{b}\right) %
\right] \Delta _{c} \\
&=&\left[ \nabla _{a+b=c}\left( p^{\otimes n-1}\Delta
_{B}^{n-2}\beta _{a}\otimes p^{\otimes 2}\Delta _{B}\beta
_{b}\right) \right] \Delta _{c}=0.
\end{eqnarray*}%
By definition of $(B_{0}^{\wedge _{B}^{n}},\delta _{n}),$ there
exists a
unique morphism%
\begin{equation*}
\gamma _{n}:\oplus _{i=0}^{n}B_{i}\rightarrow B_{0}^{\wedge
_{B}^{n+1}}
\end{equation*}%
such that%
\begin{equation*}
\delta _{n+1}\gamma _{n}=\nabla _{i=0}^{n}\beta _{i}.
\end{equation*}%
Since each $\beta _{i}$ cosplits, then $\nabla _{i=0}^{n}\beta
_{i}$ is a monomorphism. Thus also $\gamma _{n}$ is a
monomorphism. Denote by $\beta _{a}^{b}:\oplus
_{i=0}^{a}B_{i}\rightarrow \oplus _{i=0}^{b}B_{i}$ the
canonical injection when $a\leq b.$ Then we have%
\begin{equation*}
\delta _{n+2}\gamma _{n+1}\beta _{n}^{n+1}=\left( \nabla
_{i=0}^{n+1}\beta _{i}\right) \beta _{n}^{n+1}=\nabla
_{i=0}^{n}\beta _{i}=\delta _{n+1}\gamma _{n}=\delta _{n+2}\xi
_{n+1}^{n+2}\gamma _{n}
\end{equation*}%
Since $\delta _{n+2}$ is a monomorphism, we get that
\begin{equation*}
\gamma _{n+1}\beta _{n}^{n+1}=\xi _{n+1}^{n+2}\gamma _{n}
\end{equation*}%
for every $n\in
\mathbb{N}
.$ Thus $(\gamma _{n}:\oplus _{i=0}^{n}B_{i}\rightarrow
B_{0}^{\wedge
_{B}^{n+1}})_{n\in \mathbb{N}}$ defines a direct system of monomorphisms in $%
\mathcal{M}$. Now, as, by Proposition \ref{pro: V=lim V^n},
$(B,(\nabla _{i=0}^{n}\beta _{i})_{n\in
\mathbb{N}})=\underrightarrow{\lim }\left( \oplus
_{i=0}^{n}B_{i}\right) $, by Lemma \ref{lem: schiacciamento} we
have that $\left( B,(\delta _{n})_{n\in \mathbb{N}}\right)
=\underrightarrow{\lim }(B_{0}^{\wedge _{B}^{i}})_{i\in
\mathbb{N}}.$
\end{proof}

\begin{proposition}
\label{pro: ker of tensor} Let $i\in \{1,2\}.$ Let
$f_{i}:X_{i}\rightarrow Y_{i}$ be morphisms in an abelian monoidal
category $\mathcal{M}$. Let $\sigma _{i}:Y_{i}\rightarrow X_{i}$
such that $f_{i}\sigma _{i}=\mathrm{Id}_{Y_{i}}$. Then
\begin{equation*}
\mathrm{Ker\,}(f_{1}\otimes f_{2})=[\mathrm{Ker\,}(f_{1})\otimes
X_{2}]+[X_{1}\otimes \mathrm{Ker\,}(f_{2})].
\end{equation*}
\end{proposition}

\begin{proof}
Let $(K_{i},k_{i})=\mathrm{Ker\,}(f_{i})$ for $i=1,2$. Let $\nu
_{1}:K_{1}\otimes X_{2}\rightarrow (K_{1}\otimes X_{2})\oplus
(X_{1}\otimes K_{2})$ and $\nu _{2}:X_{1}\otimes K_{2}\rightarrow
(K_{1}\otimes X_{2})\oplus (X_{1}\otimes K_{2})$ be the canonical
inclusions. Then, by the universal property of coproducts, there
is a unique morphism $\tau :(K_{1}\otimes X_{2})\oplus
(X_{1}\otimes K_{2})\rightarrow X_{1}\otimes X_{2}$ such that
\begin{equation}
\tau \nu _{1}=k_{1}\otimes X_{2}\text{\qquad and\qquad }\tau \nu
_{2}=X_{1}\otimes k_{2}.  \label{fromula cioccorana}
\end{equation}%
By definition, one has $(K_{1}\otimes X_{2})+(X_{1}\otimes K_{2})=\mathrm{Im}%
(\tau )=\mathrm{Ker\,}(\pi ),$ where $(C,\pi )=\mathrm{coker}(\tau ).$%
\newline
Thus, in order to prove our statement, we will show that $(C,\pi
)=(Y_{1}\otimes Y_{2},f_{1}\otimes f_{2}).$ By (\ref{fromula
cioccorana}), we have
\begin{eqnarray*}
(f_{1}\otimes f_{2})\tau \nu _{1} &=&(f_{1}\otimes
f_{2})(k_{1}\otimes
X_{2})=0, \\
(f_{1}\otimes f_{2})\tau \nu _{2} &=&(f_{1}\otimes
f_{2})(X_{1}\otimes k_{2})=0,
\end{eqnarray*}%
so that $(f_{1}\otimes f_{2})\tau =0$. By the universal property
of cokernels, we obtain a unique morphism $\alpha :C\rightarrow
Y_{1}\otimes Y_{2}$ such that $\alpha \pi =f_{1}\otimes f_{2}.$
\newline Define $\beta :Y_{1}\otimes Y_{2}\rightarrow C$ by $\beta
:=\pi (\sigma _{1}\otimes \sigma _{2}).$ Let us prove that $\beta
$ is a two-sided inverse of $\alpha $. Clearly one has
\begin{equation*}
\alpha \beta =\alpha \pi (\sigma _{1}\otimes \sigma
_{2})=(f_{1}\otimes f_{2})(\sigma _{1}\otimes \sigma
_{2})=Id_{X_{1}\otimes X_{2}}.
\end{equation*}%
Now, since $f_{i}\sigma _{i}=\mathrm{Id}_{Y_{i}},$ there is a
unique
morphism $\rho _{i}:X_{i}\rightarrow K_{i}$ such that $\rho _{i}k_{i}=%
\mathrm{Id}_{K_{i}}$ and
\begin{equation}
k_{i}\rho _{i}+\sigma _{i}f_{i}=\mathrm{Id}_{X_{i}},\text{ for any
}i\in \{1,2\}.  \label{formula torta di riso}
\end{equation}%
Then we have:
\begin{eqnarray*}
\beta \alpha \pi &=&\beta (f_{1}\otimes f_{2}) \\
&=&\pi (\sigma _{1}f_{1}\otimes \sigma _{2}f_{2}) \\
&\overset{(\ref{formula torta di riso})}{=}&\pi \lbrack \sigma
_{1}f_{1}\otimes (\mathrm{Id}_{X_{2}}-k_{2}\rho _{2})] \\
&=&\pi (\sigma _{1}f_{1}\otimes \mathrm{Id}_{X_{2}})-\pi (\sigma
_{1}f_{1}\otimes k_{2}\rho _{2}) \\
&=&\pi \lbrack (\mathrm{Id}_{X_{1}}-k_{1}\rho _{1})\otimes \mathrm{Id}%
_{X_{2}}]-\pi (X_{1}\otimes k_{2})(\sigma _{1}f_{1}\otimes \rho _{2}) \\
&\overset{(\ref{fromula cioccorana})}{=}&\pi -\pi (k_{1}\rho
_{1}\otimes
\mathrm{Id}_{X_{2}})-\pi \tau \nu _{2}(\sigma _{1}f_{1}\otimes \rho _{2}) \\
&=&\pi -\pi (k_{1}\otimes \mathrm{Id}_{X_{2}})(\rho _{1}\otimes \mathrm{Id}%
_{X_{2}}) \\
&\overset{(\ref{fromula cioccorana})}{=}&\text{ }\pi -\pi \tau \nu
_{1}(\rho _{1}\otimes \mathrm{Id}_{X_{2}}) \\
&=&\pi
\end{eqnarray*}%
Since $\pi $ is an epimorphism we conclude that $\beta \alpha =\mathrm{Id}%
_{C}$ and hence that $\alpha $ is an isomorphism. Thus $(C,\pi
)=(Y_{1}\otimes Y_{2},f_{1}\otimes f_{2})$.
\end{proof}

\begin{proposition}
\label{pro: D^2 seconda}Let $\delta :D\rightarrow C$ be a morphism
that cosplits in $\mathcal{M}$. If $\delta $ is a coalgebra
homomorphism, then we have
\begin{equation}
D^{\wedge _{C}^{2}}=D\wedge _{C}D=\Delta _{C}^{-1}(D\otimes
C+C\otimes D). \label{formula pro: D^2}
\end{equation}
\end{proposition}

\begin{proof}
Set $(L,p)=\mathrm{coker}(\sigma ).$ Let
\begin{equation*}
(D^{\wedge _{C}^{n}},\delta _{n}):=\mathrm{Ker\,}(p^{\otimes
n}\Delta _{C}^{n-1})\text{.}
\end{equation*}%
We have
\begin{eqnarray*}
D^{\wedge _{C}^{2}} &=&\mathrm{Ker\,}[(p\otimes p)\Delta _{C}] \\
&\overset{(* )}{=}&\Delta _{C}^{-1}[\mathrm{Ker\,}(p\otimes p)] \\
&\overset{(** )}{=}&\Delta _{C}^{-1}\left\{
[\mathrm{Ker\,}(p)\otimes C]+[C\otimes \mathrm{Ker\,}(p)]\right\}
=\Delta _{C}^{-1}[(D\otimes C)+(C\otimes D)].
\end{eqnarray*}%
where in (*) we have applied (\cite[Proposition 5.1, page 90]{St})
and in (**) Proposition \ref{pro: ker of tensor}.
\end{proof}

\begin{lemma}
\label{lem: banana}Let $\sigma :A\rightarrow A^{\prime }$ be a
coalgebra
homomorphism in an abelian monoidal category $\mathcal{M}.$ Assume that $%
\sigma $ is a monomorphism. Let $i\in \{1,2\}$ and let $\alpha
_{i}:X_{i}\rightarrow A,$ $\alpha _{i}^{\prime }:X_{i}^{\prime
}\rightarrow
A,$ $\sigma _{i}:X_{i}\rightarrow X_{i}^{\prime }$ be morphism such that $%
\sigma \alpha _{i}=\alpha _{i}^{\prime }\sigma _{i}.$ Let $(L_{i},p_{i})=%
\mathrm{coker}(\alpha _{i})$ and $(L_{i}^{\prime },p_{i}^{\prime })=\mathrm{%
coker}(\alpha _{i}^{\prime })$ and let
\begin{equation*}
(X_{1}\wedge _{A}X_{2},\lambda
_{1,2})=\mathrm{Ker\,}[(p_{1}\otimes p_{2})\Delta
_{A}]\text{\qquad and\qquad }(X_{1}^{\prime }\wedge
_{A'}X_{2}^{\prime },\lambda _{1,2}^{\prime
})=\mathrm{Ker\,}[(p_{1}^{\prime }\otimes p_{2}^{\prime })\Delta
_{A^{\prime }}].
\end{equation*}%
Then there is a (unique) morphism
\begin{equation*}
\lambda =\lambda (A,X_{1},X_{2};A^{\prime },X_{1}^{\prime
},X_{2}^{\prime }):X_{1}\wedge _{A}X_{2}\rightarrow X_{1}^{\prime
}\wedge _{A^{\prime }}X_{2}^{\prime }
\end{equation*}%
such that $\lambda _{1,2}^{\prime }\lambda =\sigma \lambda
_{1,2}.$ Moreover $\lambda $ is a monomorphism.
\end{lemma}

\begin{proof}
Since
\begin{equation*}
p_{i}^{\prime }\sigma \alpha _{i}=p_{i}^{\prime }\alpha
_{i}^{\prime }\sigma _{i}=0,
\end{equation*}%
by the universal property of cokernels, there is a unique morphism
$\tau _{i}:L_{i}\rightarrow L_{i}^{\prime }$ such that $\tau
_{i}p_{i}=p_{i}^{\prime }\sigma .$ Then we have:
\begin{equation*}
(p_{1}^{\prime }\otimes p_{2}^{\prime })\Delta _{A^{\prime
}}\sigma \lambda _{1,2}=(p_{1}^{\prime }\otimes p_{2}^{\prime
})(\sigma \otimes \sigma )\Delta _{A}\lambda _{1,2}=(\tau
_{1}\otimes \tau _{2})(p_{1}\otimes p_{2})\Delta _{A}\lambda
_{1,2}=0.
\end{equation*}%
By the universal property of kernels, there is a unique morphism
$\lambda :X_{1}\wedge _{A}X_{2}\rightarrow X_{1}^{\prime }\wedge
_{A^{\prime }}X_{2}^{\prime }$ such that $\lambda _{1,2}^{\prime
}\lambda =\sigma \lambda _{1,2}.$ Clearly, as $\sigma $ and
$\lambda _{1,2}$ are monomorphisms, $\lambda $ is a monomorphism
too.
\end{proof}

\begin{lemma}
\label{lem: direct systems}Let $((X_{i})_{i\in \mathbb{N}},(\xi
_{i}^{j})_{i,j\in \mathbb{N}})$ and let $((Y_{i})_{i\in
\mathbb{N}},(\zeta
_{i}^{j})_{i,j\in \mathbb{N}})$ be direct systems in a monoidal category $%
\mathcal{M}$, where, for $i\leq j$, $\xi _{i}^{j}:X_{i}\rightarrow
X_{j}$ and $\zeta _{i}^{j}:Y_{i}\rightarrow Y_{j}.$ Let $\sigma
:A\rightarrow B$ be a coalgebra homomorphism and let $(\alpha
_{i}:X_{i}\rightarrow A)_{i\in \mathbb{N}}$ and $(\beta
_{i}:Y_{i}\rightarrow B)_{i\in \mathbb{N}}$ be compatible families
of morphisms in $\mathcal{M}.$ Let $\lambda _{i}:X_{i}\rightarrow
Y_{i}$ be a morphism such that $\beta _{i}\lambda _{i}=\sigma
\alpha _{i},$ for any $i\in \mathbb{N.}$ If $\beta _{i}$ is a
monomorphism, for any $i\in \mathbb{N}$, we have that $(\lambda
_{i}:X_{i}\rightarrow Y_{i})_{i\in \mathbb{N}}$ is a direct system
of morphisms in $\mathcal{M}.$
\end{lemma}

\begin{proof}
For any $i\leq j,$ we have that:
\begin{equation*}
\beta _{j}\lambda _{j}\xi _{i}^{j}=\sigma \alpha _{j}\xi
_{i}^{j}=\sigma \alpha _{i}=\beta _{i}\lambda _{i}=\beta _{j}\zeta
_{i}^{j}\lambda _{i}.
\end{equation*}%
Since $\beta _{j}$ is a monomorphism for any $j\in \mathbb{N}$, we
conclude
that $\lambda _{j}\xi _{i}^{j}=\zeta _{i}^{j}\lambda _{i}$ i.e. that $%
(\lambda _{i}:X_{i}\rightarrow Y_{i})_{i\in \mathbb{N}}$ is a
direct system of morphisms in $\mathcal{M}.$
\end{proof}

\begin{lemma}
\label{lem: sub direct systems}Let $\mathcal{M}$ be a cocomplete
monoidal category with left exact direct limits. Let
$((X_{i})_{i\in \mathbb{N}},(\xi _{i}^{j})_{i,j\in \mathbb{N}})$
be a direct system in $\mathcal{M}$, where, for $i\leq j$, $\xi
_{i}^{j}:X_{i}\rightarrow X_{j}$. \newline
Let $\gamma :\mathbb{N}\rightarrow \mathbb{N}$ be an injection. Then $%
((X_{\gamma (i)})_{i\in \mathbb{N}},(\xi _{\gamma (i)}^{\gamma
(j)})_{i,j\in \mathbb{N}})$ is a direct system in $\mathcal{M}$.
Let $(X,(\lambda
_{i})_{i\in \mathbb{N}})=\underrightarrow{\lim }X_{\gamma (i)}$, where $%
\lambda _{i}:$ $X_{\gamma (i)}\rightarrow X$ for any $i\in
\mathbb{N}$.
\newline
Then $(X,(\xi _{i})_{i\in \mathbb{N}})=\underrightarrow{\lim }X_{i},$ where $%
\xi _{i}:X_{i}\rightarrow X$ is defined by $\xi _{i}:=\lambda
_{j}\xi
_{i}^{\gamma (j)}:X_{i}\rightarrow X$, where $j\in \mathbb{N}$ is such that $%
\gamma (j)\geq i.$
\end{lemma}

\begin{proof}
Clearly $((X_{\gamma (i)})_{i\in \mathbb{N}},(\xi _{\gamma
(i)}^{\gamma (j)})_{i,j\in \mathbb{N}})$ is a direct system. Let
us prove the last assertion. Let $j,j^{\prime }\in \mathbb{N}$
such that $\gamma (j^{\prime })\geq \gamma (j)\geq i.$ Then
\begin{equation*}
\lambda _{j^{\prime }}\xi _{i}^{\gamma (j^{\prime })}=\lambda
_{j^{\prime }}\xi _{\gamma (j)}^{\gamma (j^{\prime })}\xi
_{i}^{\gamma (j)}=\lambda _{j}\xi _{i}^{\gamma (j)},
\end{equation*}%
so that $\xi _{i}$ is well defined. Note that
\begin{equation}
\xi _{\gamma (j)}=\lambda _{j}\xi _{\gamma (j)}^{\gamma
(j)}=\lambda _{j}. \label{formula pappardelle}
\end{equation}%
Moreover, for any $i\leq j,$ and for any $t\in \mathbb{N}$ such
that $\gamma (t)\geq j$ we have:
\begin{equation*}
\xi _{j}\xi _{i}^{j}=\lambda _{t}\xi _{j}^{\gamma (t)}\xi
_{i}^{j}=\lambda _{t}\xi _{i}^{\gamma (t)}=\xi _{i},
\end{equation*}%
so that $(\xi _{i}:X_{i}\rightarrow X)_{i\in \mathbb{N}}$ is a
direct system of morphisms. Let now $(f_{i}:X_{i}\rightarrow
Y)_{i\in \mathbb{N}}$ be a compatible family of morphisms in
$\mathcal{M}$. Then $(f_{\gamma (i)}:X_{\gamma (i)}\rightarrow
Y)_{i\in \mathbb{N}}$ is a compatible family
of morphisms in $\mathcal{M}$ so that there exists a unique morphism $%
f:X\rightarrow Y$ such that $f\lambda _{i}=f_{\gamma (i)}$ for any
$i\in \mathbb{N.}$ For any $i\in \mathbb{N}$ and for any $j\in
\mathbb{N}$ is such that $\gamma (j)\geq i,$we obtain
\begin{equation*}
f\xi _{i}=f\lambda _{j}\xi _{i}^{\gamma (j)}=f_{\gamma (i)}\xi
_{i}^{\gamma (j)}=f_{i}.
\end{equation*}%
Let $g:X\rightarrow Y$ be another morphism such that $g\xi
_{i}=f_{i}.$ Then we have
\begin{equation*}
g\lambda _{i}\overset{(\ref{formula pappardelle})}{=}g\xi _{\gamma
(i)}=f_{\gamma (i)}.
\end{equation*}%
By uniqueness of $f$ we get $g=f,$ so that $(X,(\xi _{i})_{i\in \mathbb{N}})=%
\underrightarrow{\lim }X_{i}$.
\end{proof}

\begin{lemma}
\label{lem: horrible}Let $(\mathcal{M},\otimes ,\mathbf{1})$ be a
cocomplete abelian monoidal category satisfying $AB5,$ with left
exact direct limits and left and right exact tensor functors. Let
$(\mathbb{I},\leq )$ be a directed partially ordered set. Let
$((X_{i})_{i\in \mathbb{I}},(\xi _{i}^{j})_{i,j\in \mathbb{I}})$
be a direct system in $\mathcal{M}$, where,
for $i\leq j$, $\xi _{i}^{j}:X_{i}\rightarrow X_{j}.$ Let $%
(w_{i}:X_{i}\rightarrow W)_{i\in \mathbb{I}}$ be a compatible
family of
monomorphisms in $\mathcal{M}$. Let $(X,(\xi _{i})_{i\in \mathbb{I}})=%
\underrightarrow{\lim }X_{i}.$ Let $w:X\rightarrow W$ be the
unique morphism such that $w\xi _{i}=w_{i}$ for every $i$. Let
$\xi :\oplus X_{i}\rightarrow
X$ be the unique morphism such that $\xi \varepsilon _{i}=\xi _{i}$ for any $%
i\in \mathbb{I}$ and let $\omega :\oplus X_{i}\rightarrow W$ be
the unique
morphism such that $\omega \varepsilon _{i}=w_{i}$ for any $i\in \mathbb{I}$%
, where $\varepsilon _{i}:X_{i}\rightarrow \oplus X_{i}$ is the
canonical inclusion. Then:
\begin{equation*}
w\xi =\omega .
\end{equation*}%
Moreover $w$ is a monomorphism and $\xi $ is an epimorphism.
\end{lemma}

\begin{proof}
Since $w\xi _{i}=w_{i}$, the $\xi _{i}$'s are monomorphisms.
Clearly we have
\begin{equation*}
w\xi \varepsilon _{i}=w\xi _{i}=w_{i}=\omega \varepsilon _{i}\text{ for any }%
i\in \mathbb{I}
\end{equation*}%
and hence
\begin{equation*}
w\xi =\omega
\end{equation*}%
Moreover, regarding $(w_{i}:X_{i}\rightarrow W)_{i\in \mathbb{I}}$
as a direct system of monomorphism, in view of $AB5$, we have that
$w$ is a monomorphism and $\xi $ is an epimorphism.
\end{proof}

\begin{lemma}
\label{lem: colimit and tensor}Let $(\mathcal{M},\otimes
,\mathbf{1})$ be a cocomplete abelian monoidal category satisfying
$AB5$ and with left and right exact tensor functors. Let
$(\mathbb{I},\leq )$ be a directed partially
ordered set. Let $((X_{i})_{i\in \mathbb{I}},(\xi _{i}^{j})_{i,j\in \mathbb{I%
}})$ be a direct system in $\mathcal{M}$, where, for $i\leq j$,
$\xi
_{i}^{j}:X_{i}\rightarrow X_{j}.$ If $\oplus X_{i}$ commutes with $\otimes $%
, then $\underrightarrow{\lim }X_{i}$ does.
\end{lemma}

\begin{proof}
Let $((X_{i})_{i\in \mathbb{I}},(\xi _{i}^{j})_{i,j\in
\mathbb{I}})$ be a direct system in $\mathcal{M}$, where, for
$i\leq j$, $\xi
_{i}^{j}:X_{i}\rightarrow X_{j}.$ Let $(X,(\xi _{i})_{i\in \mathbb{I}})=%
\underrightarrow{\lim }X_{i}$ and let $Y$ be an object in
$\mathcal{M}$. By \cite[Lemma 1.2, page 115]{St}, the $\xi _{i}$'s
are monomorphisms. Also, by the universal property of the
coproduct, there is a unique morphism $\xi :\oplus
X_{i}\rightarrow X$ such that
\begin{equation}
\xi \varepsilon _{i}=\xi _{i}\text{ for any }i\in \mathbb{I}
\label{form2}
\end{equation}%
where $\varepsilon _{i}:X_{i}\rightarrow \oplus X_{i}$ is the
canonical inclusion. Moreover, $\xi $ is an epimorphism. \newline
Assume that
\begin{equation}
((\oplus X_{i})\otimes Y,\varepsilon _{i}\otimes Y)=\oplus
(X_{i}\otimes Y), \label{formula coproduct}
\end{equation}%
Let $\gamma _{i}:X_{i}\otimes Y\rightarrow \underrightarrow{\lim }%
(X_{i}\otimes Y)$ be the canonical morphism. By the universal
property of
coproduct and by (\ref{formula coproduct}), there is a unique morphism $%
\gamma :(\oplus X_{i})\otimes Y\rightarrow \underrightarrow{\lim }%
(X_{i}\otimes Y)$ such that
\begin{equation}
\gamma (\varepsilon _{i}\otimes Y)=\gamma _{i}\text{ for any }i\in \mathbb{I}%
\text{.}  \label{form1}
\end{equation}%
In an analogous way, by the universal property of direct limits,
there is a unique morphism $\Lambda :\underrightarrow{\lim
}(X_{i}\otimes Y)\rightarrow X\otimes Y$ such that
\begin{equation}
\Lambda \gamma _{i}=\xi _{i}\otimes Y\text{ for any }i\in
\mathbb{I}\text{.} \label{form3}
\end{equation}%
It is easy to see that we can apply Lemma \ref{lem: horrible} to
the present situation and get:
\begin{equation*}
\Lambda \gamma =\xi \otimes Y.
\end{equation*}%
where $\Lambda $ is a monomorphism and $\gamma $ is an
epimorphism. Moreover, since the tensor functor is left exact and
$\xi $ is an epimorphism, we get that $\xi \otimes Y$ is an
epimorphism. Hence $\Lambda $ is an epimorphism too.
\end{proof}

\begin{theorem}
\label{teo alfa beta sigma} Let $(C,\Delta ,\varepsilon )$ be a
coalgebra in a cocomplete abelian monoidal category $\mathcal{M}$
satisfying $AB5,$ with left and right exact tensor functors.
Assume that denumerable coproducts commute with $\otimes $. Let
$\alpha :C\rightarrow A$ and $\sigma :A\rightarrow B$ be
monomorphisms which are coalgebra homomorphisms and let $\beta
=\sigma
\alpha .$ Assume that $\sigma $ cosplits in $\mathcal{M}$. Let $p_{\alpha }=%
\mathrm{coker}(\alpha )$ in $\mathcal{M}$, let $(C^{\wedge
_{A}^{n}},\alpha _{n}):=\mathrm{Ker\,}(p_{\alpha }^{\otimes
n}\Delta _{A}^{n-1})$ and assume that $\alpha _{n}$ cosplits in
$\mathcal{M}$, for every $n\in \N.$ If $\widetilde{C}_{A}=A$ and $B=A^{\wedge _{B}2},$ then $\widetilde{C}%
_{B}=B.$
\end{theorem}

\begin{proof}
For any morphism $\eta $ we set $(L_{\eta },p_{\eta
})=\mathrm{coker}(\eta )$ in $\mathcal{M}.$ By Proposition
\ref{pro: D^2 seconda}, we get
\begin{equation*}
B=A^{\wedge _{B}2}\overset{(\ref{formula pro: D^2})}{=}\Delta
_{B}^{-1}(A\otimes B+B\otimes A).
\end{equation*}%
Let
\begin{equation*}
(C^{\wedge _{A}^{n}},\alpha _{n}):=\mathrm{Ker\,}(p_{\alpha
}^{\otimes n}\Delta _{A}^{n-1})\text{\qquad and\qquad }(C^{\wedge
_{B}^{n}},\beta _{n}):=\mathrm{Ker\,}(p_{\beta }^{\otimes n}\Delta
_{B}^{n-1}).
\end{equation*}%
By assumption $(A,(\alpha _{n})_{n\in \mathbb{N}})=\widetilde{C}_{A}=%
\underrightarrow{\lim }C^{\wedge _{A}^{n}}.$ Then, by Lemma
\ref{lem: colimit and tensor}, we obtain:
\begin{equation*}
A\otimes B=(\underrightarrow{\lim }C^{\wedge _{A}^{n}})\otimes B=%
\underrightarrow{\lim }(C^{\wedge _{A}^{n}}\otimes B).
\end{equation*}%
We have:
\begin{eqnarray*}
B &=&\Delta _{B}^{-1}[\underrightarrow{\lim }(C^{\wedge _{A}^{m}}\otimes B)+%
\underrightarrow{\lim }(B\otimes C^{\wedge _{A}^{n}})] \\
&=&\Delta _{B}^{-1}\underrightarrow{\lim }[(C^{\wedge
_{A}^{n}}\otimes
B)+(B\otimes C^{\wedge _{A}^{n}})] \\
&{=}&\underrightarrow{\lim }\Delta _{B}^{-1}[(C^{\wedge
_{A}^{n}}\otimes
B)+(B\otimes C^{\wedge _{A}^{n}})] \\
&=&\underrightarrow{\lim }(C^{\wedge _{A}^{n}}\wedge _{B}C^{\wedge
_{A}^{n}})
\end{eqnarray*}%
where in the second equality we have used that in an
$AB5$-category direct limits of direct systems of subobjects are
just sums of their respective
families; in the third we have used a well known property of $AB5$%
-categories (see \cite[Proposition 1.1, page 114]{St}); in the
last equality we have used Proposition \ref{pro: D^2 seconda} in
the case $\delta =\sigma \alpha _{n}:C^{\wedge
_{A}^{n}}\rightarrow B.$ Note that, by applying inductively Lemma
\ref{lem: banana}, we obtain, for any $n\in \mathbb{N},$ a
monomorphism
\begin{equation*}
\lambda _{A,B}^{n}:C^{\wedge _{A}^{n}}\rightarrow C^{\wedge
_{B}^{n}}
\end{equation*}%
such that
\begin{equation}
\beta _{n}\lambda _{A,B}^{n}=\sigma \alpha _{n}.
\label{compatibility uffaa}
\end{equation}%
By Lemma \ref{lem: direct systems}, $(\lambda _{A,B}^{n}:C^{\wedge
_{A}^{n}}\rightarrow C^{\wedge _{B}^{n}})_{i\in \mathbb{N}}$ is a
direct system of monomorphisms in $\mathcal{M}.$ Let $m\leq n.$
Note that, if we denote by $\xi _{A,m}^{n}:C^{\wedge
_{A}^{m}}\rightarrow C^{\wedge _{A}^{n}}$ and by $\xi
_{B,m}^{n}:C^{\wedge _{B}^{m}}\rightarrow C^{\wedge _{B}^{n}}$ the
canonical morphisms, this means that $\lambda _{A,B}^{n}\xi
_{A,m}^{n}=\xi _{B,m}^{n}\lambda _{A,B}^{m}$.\newline
Let $(L_{n},p_{n})=\mathrm{coker}(\sigma _{n}\alpha _{n})$ and $%
(L_{n}^{\prime },p_{n}^{\prime })=\mathrm{coker}(\beta _{n})$ and
let
\begin{equation*}
(C^{\wedge _{A}^{n}}\wedge _{B}C^{\wedge _{A}^{n}},\lambda _{n})=\mathrm{%
Ker\,}[(p_{n}\otimes p_{n})\Delta _{B}]\text{\qquad and\qquad
}(C^{\wedge
_{B}^{n}}\wedge _{B}C^{\wedge _{B}^{n}},\lambda _{n}^{\prime })=\mathrm{Ker\,%
}[(p_{n}^{\prime }\otimes p_{n}^{\prime })\Delta _{B}].
\end{equation*}%
In view of (\ref{compatibility uffaa}), we can apply Lemma
\ref{lem: banana} to obtain a monomorphism
\begin{equation*}
\lambda (n)=\lambda (B,C^{\wedge _{A}^{n}},C^{\wedge
_{A}^{n}};B,C^{\wedge _{B}^{n}},C^{\wedge _{B}^{n}}):C^{\wedge
_{A}^{n}}\wedge _{B}C^{\wedge _{A}^{n}}\rightarrow C^{\wedge
_{B}^{n}}\wedge _{B}C^{\wedge _{B}^{n}}
\end{equation*}%
such that
\begin{equation*}
\lambda _{n}^{\prime }\lambda (n)=\lambda _{n}.
\end{equation*}%
By Lemma \ref{lem: direct systems}, $(\lambda (n):C^{\wedge
_{A}^{n}}\wedge _{B}C^{\wedge _{A}^{n}}\rightarrow C^{\wedge
_{B}^{n}}\wedge _{B}C^{\wedge
_{B}^{n}})_{n\in \mathbb{N}}$ is a direct system of monomorphisms in $%
\mathcal{M}.$ Since, by the foregoing, $B=\underrightarrow{\lim
}(C^{\wedge
_{A}^{n}}\wedge _{B}C^{\wedge _{A}^{n}}),$ by Lemma \ref{lem: schiacciamento}%
, applied in the case $\gamma _{i}=\lambda (i)$ for any $i\in
\mathbb{N}$, we obtain that
\begin{equation*}
(B,(\lambda _{n}^{\prime })_{n\in \mathbb{N}})=\underrightarrow{\lim }%
(C^{\wedge _{B}^{n}}\wedge _{B}C^{\wedge _{B}^{n}}).
\end{equation*}%
As
\begin{equation*}
(C^{\wedge _{B}^{n}}\wedge _{B}C^{\wedge _{B}^{n}},\lambda _{n}^{\prime })%
\overset{(\ref{formula 2 pro: D^2})}{=}(C^{\wedge _{B}^{2n}},\beta
_{2n})
\end{equation*}%
we obtain
\begin{equation*}
(B,(\beta _{2n})_{n\in \mathbb{N}})=\underrightarrow{\lim
}C^{\wedge _{B}^{2n}}.
\end{equation*}%
Now, apply Lemma \ref{lem: sub direct systems} in the case when $\gamma :%
\mathbb{N}\rightarrow \mathbb{N}$ is defined by setting $\gamma
(n)=2n$ for every $n\in \mathbb{N.}$ Then we get
\begin{equation*}
(B,(\beta _{n})_{n\in \mathbb{N}})=\underrightarrow{\lim }C^{\wedge _{B}n}=%
\widetilde{C}_{B}.
\end{equation*}
\end{proof}

\section{Formal Smoothness of the Cotensor Coalgebra}
The main aim of this section is to prove Theorem \ref{teo Smooth}
which asserts that the cotensor coalgebra $T_{C}^{c}(M)$ is
formally smooth whenever $C$ is a formally smooth coalgebra in a
cocomplete and complete abelian monoidal category $\M$ satisfying
$AB5,$ and $M$ is an $\mathcal{I}$-injective $C$-bicomodule in
$\M$.

\begin{definition}\label{standard complex}
Let $(C,\Delta ,\varepsilon )$ be a coalgebra in $(\mathcal{M},\otimes ,%
\mathbf{1})$ and let $\left( L,\rho _{L}^{l},\rho _{L}^{r}\right) $ be a $C$%
-bicomodule. Let us consider \emph{the standard complex}:
\begin{equation*}
0\longrightarrow \mathcal{M}(L,\mathbf{1})\overset{b^{0}}{\longrightarrow }%
\mathcal{M}(L,{C})\overset{b^{1}}{\longrightarrow }\mathcal{M}(L,{C\otimes C}%
)\overset{b^{2}}{\longrightarrow }\mathcal{M}(L,{C\otimes C\otimes C})%
\overset{b^{3}}{\longrightarrow }\cdots
\end{equation*}%
Let $n\in \N.$ For every $0\leq i\leq n+1$ and for every $f\in {\mathcal{M}}%
(L,C^{\otimes n}),$ we define,\\ for $n=0:$%
\begin{equation*}
b_{0}^{0}=l_{C}\circ (f\otimes C)\circ \rho _{L}^{r},\text{\qquad }%
b_{1}^{0}=r_{C}\circ \left( C\otimes f\right) \circ \rho _{L}^{l}
\end{equation*}%
and, for $n>0$:
\begin{equation*}
b_{i}^{n}(f)=\left\{
\begin{array}{ll}
(f\otimes C)\circ \rho _{L}^{r} & i=0; \\
(C^{\otimes n-i}\otimes \Delta \otimes C^{\otimes i-1})\circ f, &
i=1,\ldots
,n; \\
\left( C\otimes f\right) \circ \rho _{L}^{l} & i=n+1.%
\end{array}%
\right.
\end{equation*}%
One has that
\begin{equation*}
b^{n}(f)=\sum\nolimits_{i=0}^{n+1}(-1)^{i}b_{i}^{n}(f),\text{ for
every }n\geq0.
\end{equation*}%
In particular, for $n\in \{0,1,2\}$ the differentials $b^{n}$ are
given by:
\begin{eqnarray*}
b^{0}(f) &=&l_{C}\circ (f\otimes C)\circ \rho _{L}^{r}-r_{C}\circ
(C\otimes
f)\circ \rho _{L}^{l}; \\
b^{1}(f) &=&(f\otimes C)\circ \rho _{L}^{r}-\Delta \circ
f+(C\otimes f)\circ
\rho _{L}^{l}; \\
b^{2}(f) &=&(f\otimes C)\circ \rho _{L}^{r}-(C\otimes \Delta
)\circ f+(\Delta \otimes C)\circ f-(C\otimes f)\circ \rho
_{L}^{l}.
\end{eqnarray*}%
Further details can be found in \cite{AMS}.
\end{definition}

\begin{claim}
\label{cl:E-Proj} Let $\mathcal{M}$ be an abelian category and let $\mathcal{%
H}$ be a class of monomorphisms in $\mathcal{M}$. We recall that an object $%
I $ in $\mathcal{M}$ is called \emph{injective rel} $\lambda $, where $%
\lambda :X\rightarrow Y$ is a monomorphism in $\mathcal{H}$, if $\mathcal{M}%
(\lambda ,I):\mathcal{M}(Y,I)\rightarrow \mathcal{M}(X,I)$ is
surjective. $I$ is called $\mathcal{H}$-injective if it is
injective rel $\lambda $ for every $\lambda $ in $\mathcal{H}$.
The\emph{\ closure} of $\mathcal{H}$ is
the class $\mathcal{C(H)}$ containing all monomorphisms $\lambda $ in $%
\mathcal{M}$ such that every $\mathcal{H}$-injective object is
also injective rel $\lambda $. The class $\mathcal{H}$ is called
\emph{closed} if $\mathcal{H}=\mathcal{C(H)}$. A closed class
$\mathcal{H}$ is called \emph{injective} if for any object $M$ in
$\mathcal{M}$ there is a
monomorphism $\lambda :M\rightarrow I$ in $\mathcal{H}$ such that $I$ is $%
\mathcal{H}$-injective.
\end{claim}

\begin{claim}
\label{claim injective class}We fix a coalgebra $C$ in an abelian
monoidal category $\mathcal{M}$. Let
$\mathbb{U}:{^{C}\mathcal{M}^{C}}\rightarrow \mathcal{M}$ be the
forgetful functor. Then
\begin{equation}
\mathcal{I}:=\{f\in {^{C}\mathcal{M}^{C}}\mid \mathbb{U}(f)\text{
cosplits in }\mathcal{M}\}.  \label{inj class}
\end{equation}%
is an injective class of monomorphisms.\newline Now, for any $C$-
bicomodule $M\in {}^{C}\!\mathcal{M}^{C}$, we define the
Hochschild cohomology of $C$ with coefficients in $M$ by:
\begin{equation*}
\mathbf{H}^{\bullet }(M,C)=\mathbf{Ext}_{\mathcal{I}}^{\bullet
}(M,C),
\end{equation*}%
where $\mathbf{Ext}_{\mathcal{I}}^{\bullet }(M,-)$ are the
relative left derived functors of $^{C}\mathcal{M}^{C}(M,-)$. The
notion of Hochschild cohomology for algebras and coalgebras in
monoidal categories has been deeply investigated in \cite{AMS}.
Here we quote some results that will be needed afterwards.
\end{claim}

\begin{definition}
A coalgebra in $\mathcal{M}$ is called \emph{coseparable} if the
comultiplication $\Delta :C\rightarrow C\otimes C$ has a retraction in ${}^{C}%
\mathcal{M}^{C}$.
\end{definition}

\begin{theorem}
\label{teo separable} \cite[Theorem 4.4]{AMS} Let $C$ be a
coalgebra in an abelian monoidal category $\mathcal{M}$. Then the
following assertions are equivalent:

(a) $C$ is coseparable.

(b) $C$ is $\mathcal{I}-$injective.

(c) $\mathrm{H}^{1}(M,C)=0$, for all $M\in
{}^{C}\!\mathcal{M}^{C}$.

(d) $\mathrm{H}^{n}(M,C)=0$, for all $M\in {}^{C}\!\mathcal{M}^{C}$ and $n>0$%
.

(e) Any morphism in $^{C}\mathcal{M}^{C}$ cosplits in
$^{C}\mathcal{M}^{C}$ whenever it cosplits in $\mathcal{M}$.

(f) The category $^{C}\mathcal{M}^{C}$ is
$\mathcal{I}-$cosemisimple (i.e. every object in
$^{C}\mathcal{M}^{C}$ is $\mathcal{I}$-injective).
\end{theorem}

\begin{definition}
An extension $\sigma :C\rightarrow E$ is a \emph{trivial extension }%
whenever it admits a retraction that is a coalgebra homomorphism.
\end{definition}

\begin{definition}
We say that the sequence $\left( \left( E_{i}\right) _{i\in
\mathbb{N}
},\left( \eta _{i}^{j}\right) _{i,j\in
\mathbb{N}
}\right) $ of morphisms in an abelian monoidal category $\mathcal{M}$%
\begin{equation}
E_{1}\rTo^{\eta _{1}^{2}}E_{2}\rTo^{\eta _{2}^{3}}\cdots
\rTo^{\eta _{n-1}^{n}}E_{n}\rTo^{\eta _{n}^{n+1}}E_{n+1}\rTo^{\eta
_{n+1}^{n+2}}\cdots \label{ec:sistem direct}
\end{equation}%
is a \emph{direct system of extensions} if $\eta _{n}^{n+1}$ is a
coalgebra homomorphism and $E_{n}\wedge _{E_{n+1}}E_{n}=E_{n+1},$
for any $n\geq 1.$ We say that a direct system of extensions
$\left( \left( E_{i}\right) _{i\in
\mathbb{N}
},\left( \eta _{i}^{j}\right) _{i,j\in
\mathbb{N}
}\right) $ is a \textit{direct system of Hochschild extensions} if
each $\eta _{n}^{n+1}$ has a retraction in $\mathcal{M}$.
\end{definition}

\begin{example}
Let $C$ and $E$ be coalgebras in an abelian monoidal category
$\mathcal{M}$. Let $\delta :C\rightarrow E$ be a monomorphism
which is a homomorphism of coalgebras in $\mathcal{M}$. The
sequence $\left( \left( C^{\wedge _{E}^{i}}\right) _{i\in
\mathbb{N}
},\left( \xi _{i}^{j}\right) _{i,j\in
\mathbb{N}
}\right) $ is a direct system of coalgebras extensions that will
be called the $C$-\emph{adic direct system in }$E.$
\end{example}

\begin{definition}
Let $\mathcal{M}$ be an abelian monoidal category. We say that the
direct system of extensions $\left( \left( E_{i}\right) _{i\in \N
},\left( \eta _{i}^{j}\right) _{i,j\in \N }\right) $ \emph{has a
direct limit} if $\underrightarrow{\lim }E_{i}$ exists in the category $%
\mathfrak{Coalg}(\mathcal{M})$ of coalgebras in $\mathcal{M}.$
\end{definition}

\begin{remark}
If $\left( \left( E_{i}\right) _{i\in \N },\left( \eta
_{i}^{j}\right) _{i,j\in \N }\right) $ is a direct system of
Hochschild extensions, then, for any $n\geq
1,$%
\begin{equation*}
0\rTo E_{n}\rTo^{\eta
_{n}^{n+1}}E_{n+1}\rTo\mathrm{\mathrm{Coker}}\left( \eta
_{n}^{n+1}\right) \rTo0
\end{equation*}%
is a Hochschild extension of $E_{n}$ with cokernel
$\mathrm{Coker}\left( \eta _{n}^{n+1}\right).$
\end{remark}

\begin{theorem}
\cite[Theorem 4.16]{AMS}\label{teo: smooth}\label{X cof teo}
\label{X
coformally smooth teo}Let $(C,\Delta ,\varepsilon )$ be a coalgebra in $%
\mathcal{M}$. Then the following conditions are equivalent:

(a) The canonical map $\mathrm{Hom}_{\mathrm{coalg}}(\delta,C):\mathrm{Hom}_{\mathrm{coalg}}(E,C)\rightarrow \mathrm{%
Hom}_{\mathrm{coalg}}(D,C)$ is surjective for every coalgebra homomorphism $%
\delta :D\rightarrow E$ that cosplits in $\mathcal{M}$ and such that $%
D\wedge _{E}D=E.$

(b) The canonical map $\mathrm{Hom}_{\mathrm{coalg}}(\underrightarrow{\lim }%
E_{n},C)\rightarrow \mathrm{Hom}_{\mathrm{coalg}}(E_{1},C)$ is
surjective for every direct system of Hochschild extensions
$\left( \left( E_{i}\right) _{i\in\N },\left( \eta _{i}^{j}\right)
_{i,j\in\N }\right) $ which has direct limit
$\underrightarrow{\lim }E_{n}$.

(c) The canonical map $\mathrm{Hom}_{\mathrm{coalg}}(\underrightarrow{\lim }%
D^{\wedge _{E}^{n}},C)\rightarrow
\mathrm{Hom}_{\mathrm{coalg}}(D,C)$ is surjective for any
coalgebra $E$ in $\mathcal{M}$ and any subcoalgebra $D$ of $E$
such that $\left( \left( D^{\wedge _{E}^{i}}\right) _{i\in
\mathbb{N}
},\left( \xi _{i}^{j}\right) _{i,j\in
\mathbb{N}
}\right) $ is a direct system of Hochschild extensions which has
direct limit $\underrightarrow{\lim }D^{\wedge _{E}^{n}}.$

(d) Any Hochschild extension of $C$ is trivial.

(e) $\mathrm{H}^{2}\left( M,C\right) =0$, for any $M\in {}^{C}{{\mathcal{M}}}%
^{C}{.}$
\end{theorem}

\begin{definition}\label{def: formally smooth}
Any coalgebra $(C,\Delta,\varepsilon)$ in $(\mathcal{M},\otimes
,\mathbf{1})$ satisfying one of the conditions of Theorem
\ref{teo: smooth}, is called \emph{formally smooth}.
\end{definition}

\begin{theorem}
\cite[Corollary 4.21]{AMS}\label{X formally smooth teo2} Let
$(C,\Delta
,\varepsilon )$ be a coalgebra in an abelian monoidal category $\mathcal{M}$%
. Then the following assertions are equivalent:

(a) $C$ is formally smooth.

(b) $\mathrm{Coker}\left( \Delta \right) $ is
$\mathcal{I}$-injective, where $\Delta $ is the comultiplication
of $C$.
\end{theorem}

\begin{remark}
  We point out that, in view of Proposition \ref{pro: D^2}, our
  definition of $D^{\wedge _{E}^{n}}$ and the one given in
  \cite{AMS} agree.
\end{remark}

\begin{theorem}
\label{teo Smooth} Let $(C,\Delta ,\varepsilon )$ be a formally
smooth coalgebra in a cocomplete and complete abelian monoidal
category $\mathcal{M} $ satisfying $AB5,$ with left and right
exact tensor functors. Assume that denumerable coproducts commute
with $\otimes $. Let $(M,\rho _{M}^{r},\rho _{M}^{l})$ be an
$\mathcal{I}$-injective $C$-bicomodule. Then the cotensor
coalgebra $T_{C}^{c}(M)$ is formally smooth.
\end{theorem}

\begin{proof}
We will prove that any Hochschild extension of $T:=T_{C}^{c}(M)$
is trivial. Let $\sigma :T\rightarrow B$ be a Hochschild extension
of $T.$ Since the
canonical projection $p_{0}:T\rightarrow C$ is a coalgebra homomorphism and $%
C$ is formally smooth, by $a)$ of Theorem \ref{teo: smooth},  there exists a coalgebra homomorphism $%
f_{C}:B\rightarrow C$ such that $f_{C}\sigma =p_{0}.$ Then $B$ is a $C$%
-bicomodule via $f_{C}$. Moreover $\sigma $ becomes a morphism of $C$%
-bicomodules. Since $M$ is $\mathcal{I}$-injective and the
canonical projection $p_{1}:T\rightarrow M$ is a morphism of
$C$-bicomodules, then
there is a morphism of $C$-bicomodules $f_{M}:B\rightarrow M$ such that $%
f_{M}\sigma =p_{1}.$ Since, by Proposition \ref{pro: wedge limit},
we have
\begin{equation*}
\widetilde{C}_{T}=\underrightarrow{\lim }C^{\wedge
_{T}^{n}}=\underrightarrow{\lim }{C^n(M)} =T
\end{equation*}%
since $\sigma _{n}:C^{n}\rightarrow T$ cosplits and since, by
definition of Hochschild extension, $B=T^{\wedge _{B}^{2}}$ and
$\sigma $ cosplits, then by Theorem \ref{teo alfa beta sigma}
applied to the case $"\alpha "=i_{0}:C\rightarrow T$ the canonical
inclusion and $"\sigma "=\sigma ,$ we have $\widetilde{C}_{B}=B.$
Now we have
\begin{equation*}
f_{M}\sigma i_{0}=p_{1}i_{0}=0.
\end{equation*}%
Therefore we can apply Theorem \ref{coro: univ property of
cotensor coalgebra} in the case when $"C"="D"=C$, $"M"=M$, $"E"=B$
and $"\delta "=\sigma i_{0}$ in order to obtain a unique coalgebra
homomorphism $f:B\rightarrow T $ such that $p_{0}f=f_{C}$ and
$p_{1}f=f_{M}.$ Then we have
\begin{equation*}
p_{0}f\sigma =f_{C}\sigma =p_{0},\text{\qquad and\qquad
}p_{1}f\sigma =f_{M}\sigma =p_{1}.
\end{equation*}%
By Theorem \ref{coro: univ property of cotensor coalgebra} applied to $%
\delta =i_{0}:C\rightarrow T,$ the morphism $\text{Id}_{T}$ is the
unique
coalgebra homomorphism such that $p_{0}\text{Id}_{T}=p_{0}$ and $p_{1}\text{%
Id}_{T}=p_{1}$, Therefore, we conclude that $f\sigma
=\text{Id}_{T}$.
\end{proof}

\begin{corollary}
Let $(C,\Delta ,\varepsilon )$ be a coalgebra in a cocomplete and
complete abelian monoidal category $\mathcal{M}$ satisfying $AB5,$
with left and right exact tensor functors. Assume that denumerable
coproducts commute with $\otimes $. \newline a) If $C$ is formally
smooth, then

a1) $T_{C}^{c}(C\otimes X\otimes C)$ is formally smooth, for any
$X\in \mathcal{M}$. In particular $T_{C}^{c}(C^{\otimes n})$ is
formally smooth, for any $n>1$.

a2) The cotensor coalgebra $T_{C}^{c}(\text{Coker}(\Delta ))$ is
formally smooth. \newline b) If $C$ is coseparable, then the
cotensor coalgebra $T_{C}^{c}(M)$ is formally smooth, for any
$C$-bicomodule $M$.
\end{corollary}

\begin{proof}
We will apply Theorem \ref{teo Smooth}.\newline
a1) By the dual of \cite[Theo 1.16]{AMS}, all the objects of the form $%
C\otimes X\otimes C$, where $X\in \mathcal{M}$, are
$\mathcal{I}$-injective.
\newline
a2) By Theorem \ref{X formally smooth teo2}, since $C$ is formally
smooth, we have that $\mathrm{Coker\,}\Delta $ is
$\mathcal{I}$-injective. \newline
b) By Theorem \ref{teo separable}, any $C$-bicomodule $M$ is $\mathcal{I}$%
-injective.
\end{proof}

\section{Applications.}

\begin{lemma}
\label{lem union}Let $R$ be a ring and let $\mathcal{M}$ be the
category of right $R$-modules. Let $((X_{i})_{i\in
\mathbb{I}},(\xi _{i}^{j})_{i,j\in
\mathbb{I}})$ be a direct system in $\mathcal{M}$, where, for $i\leq j$, $%
\xi _{i}^{j}:X_{i}\rightarrow X_{j}$, and let $(\xi
_{i}:X_{i}\rightarrow
C)_{i\in \mathbb{I}}$ be a compatible family of monomorphism in $\mathcal{M}%
. $ Assume that $(C,(\xi _{i})_{i\in \mathbb{I}})=\underrightarrow{\lim }%
((X_{i})_{i\in \mathbb{I}},(\xi _{i})_{i\in \mathbb{I}}).$ Then $%
C=\bigcup_{i\in \mathbb{I}}\mathrm{Im}(\xi _{i}).$
\end{lemma}

\begin{proof}
Let $i\in \mathbb{I}$ and let $\tau _{i}:\mathrm{Im}(\xi
_{i})\hookrightarrow C,$ $\zeta _{i}:\mathrm{Im}(\xi
_{i})\rightarrow \bigcup_{i\in \mathbb{I}}\mathrm{Im}(\xi _{i}),$
$\tau :\bigcup_{i\in \mathbb{I}}\mathrm{Im}(\xi _{i})\rightarrow
C$ be the canonical inclusions. We have that $\tau \zeta _{i}=\tau
_{i},$ for any $i\in \mathbb{I.}$ Note that, for $i\leq j,$ we
have $\mathrm{Im}(\xi _{i})\subseteq \mathrm{Im}(\xi _{j})$. Let
$\alpha _{i}:X_{i}\rightarrow \mathrm{Im}(\xi _{i})$ be the
corestriction to $\mathrm{Im}(\xi _{i})$ of $\xi _{i}.$ Then
$\alpha _{i}$ is an isomorphism. Note that $\tau _{i}\alpha
_{i}=\xi _{i},$ for any $i\in \mathbb{I.}\newline
$Since $(\zeta _{i}\alpha _{i}:X_{i}\rightarrow \bigcup_{i\in \mathbb{I}}%
\mathrm{Im}(\xi _{i}))_{i\in \mathbb{I}}$ is a compatible family
of
morphisms, there exists a unique morphism $\lambda :C=\underrightarrow{\lim }%
((X_{i})_{i\in \mathbb{I}},(\xi _{i})_{i\in
\mathbb{I}})\rightarrow \bigcup_{i\in \mathbb{I}}\mathrm{Im}(\xi
_{i})$ such that $\lambda \xi _{i}=\zeta _{i}\alpha _{i},$ for any
$i\in \mathbb{I}$. Therefore we have
\begin{equation*}
\tau \lambda \xi _{i}=\tau \zeta _{i}\alpha _{i}=\tau _{i}\alpha
_{i}=\xi _{i},\text{ for any }i\in \mathbb{I,}
\end{equation*}
so that $\tau \lambda =\mathrm{Id}_{C}.$ Therefore $\tau $ is
surjective, i.e. $C=\bigcup_{i\in \mathbb{I}}\mathrm{Im}(\xi
_{i})$.
\end{proof}

\begin{lemma}
\label{lem cicciobello}Let $\delta :D\rightarrow C$ be a
monomorphism which
is a coalgebra homomorphism in the monoidal category $(\mathfrak{M}%
_{K},\otimes_{K},K)$ of vector spaces over a field $K$. Assume that $%
\widetilde{D}_{C}=C$. Then $\mathrm{Corad}(C)\subseteq
\text{Im}(\delta ).$
\end{lemma}

\begin{proof}
Denote by $(L,p)$ the cokernel of $\delta $ in $\mathcal{M}=(\mathfrak{M}%
_{K},\otimes_{K},K)$. Recall that
\begin{equation*}
(D^{\wedge _{C}^{n}},\delta _{n}):=\text{ker}(p^{\otimes
{n}}\Delta _{C}^{n-1})
\end{equation*}
for any $n\in \mathbb{N}\setminus \{0\}$, where $\Delta
^{n}:C\rightarrow C^{\otimes {n+1}}$ is the $n^{\text{th}}$
iterated comultiplication of $C.$ By Proposition \ref{pro: limit
of delta}, there are suitable morphisms $(\xi
_{i}^{j}:D^{\wedge _{C}^{i}}\rightarrow D^{\wedge _{C}^{j}})_{i,j\in \mathbb{%
N}}$ such that $((D^{\wedge _{C}^{i}})_{i\in \mathbb{N}},(\xi
_{i}^{j})_{i,j\in \mathbb{N}})$ is a direct system in $\mathcal{M}$ and $%
(\delta _{n})_{n\in \mathbb{N}}$ is a compatible family. Moreover
recall that $\widetilde{D}_{C}=\underrightarrow{\lim }((D^{\wedge
_{C}^{i}})_{i\in
\mathbb{N}},(\xi _{i})_{i\in \mathbb{N}})$ in $\mathcal{M}$. By Lemma \ref%
{lem union}, $C=\bigcup_{i\in \mathbb{N}}\mathrm{Im}(\delta _{i}).$ Let $%
D_{i}:=\mathrm{Im}(\delta _{i+1})$ and note that
\begin{equation*}
D_{i}=\left\{ c\in C\mid p^{\otimes {i+1}}\Delta
_{C}^{i}(c)=0\right\} .
\end{equation*}
for any $i\geq 0.$ Note that
\begin{equation}
\mathrm{Im}(\delta )=D_{0}\subseteq D_{1}\subseteq \cdots
\subseteq C. \label{coalgebra filtration}
\end{equation}
By \cite[Theorem 9.1.6, page 191]{Sw}, $C$ is a filtered
coalgebra, with
coalgebra filtration given by (\ref{coalgebra filtration}) so that, by \cite[%
Proposition 11.1.1, page 226]{Sw}, we have that
$\mathrm{Corad}(C)\subseteq D_{0}$.
\end{proof}

\begin{theorem}
\label{teo Nichols} Let $\mathcal{M}=(\mathfrak{M}_{K},\otimes
_{K},K)$ be the monoidal category of vector spaces over a field
$K$, let $E$ be a coalgebra in $\mathcal{M}$ and let $M$ be a
vector space. The following
assertions are equivalent for a morphism $g:E\rightarrow M$ in $\mathcal{M}$:%
$\newline \left( i\right) $ $g(\mathrm{Corad}(E))=0.\newline
\left( ii\right) $ There is a monomorphism $\delta :D\rightarrow
E$ which is a coalgebra homomorphism in $\mathcal{M}$ such that
$\widetilde{D}_{E}=E$ and $g\delta =0.$
\end{theorem}

\begin{proof}
$(i)\Longrightarrow (ii)$ Choose $D:=\mathrm{Corad}(E)$ and let
$\delta :D\rightarrow E$ be the canonical inclusion.

$(ii)\Longrightarrow (i)$ By Lemma \ref{lem cicciobello}, we have $\mathrm{%
Corad}(E)\subseteq \mathrm{Im}(\delta ).$
\end{proof}

\begin{remark}
By means of Theorem \ref{coro: univ property of cotensor
coalgebra} and Theorem \ref{teo Nichols}, in the particular case
when $\mathcal{M}$ is the
category $(\mathfrak{M}_{K},\otimes _{K},K)$, of vector spaces over a field $%
K$, the universal property of our cotensor coalgebra is equivalent
to Nichols's one.
\end{remark}

\begin{claim}
\label{cl:MonFun}A \emph{monoidal functor} between two monoidal categories $(%
\mathcal{M},\otimes ,\mathbf{1},a,l,r\mathbf{)}$ and $\mathbf{(}\mathcal{M}%
^{\prime },\otimes ,\mathbf{1},a,l,r\mathbf{)}$ is a triple
$(F,\phi _{0},\phi _{2}),$ where $F:\mathcal{M}\rightarrow
\mathcal{M}^{\prime }$ is a functor, $\phi
_{0}:\mathbf{1}\rightarrow F(\mathbf{1})$ is an isomorphism such
that the diagrams
\begin{equation*}
\begin{diagram}[h=2em,w=1.5em] \mathbf{1}\otimes F(U) & \rTo^{ l_{F(U)} } &
F(U) \\ \dTo<{ \phi_0\otimes F(U) } & & \uTo>{ F(l_U) } \\
F(\mathbf{1})\otimes F(U) & \rTo_{\phi_2( \mathbf{1} ,U) } & F(
\mathbf{1}\otimes U) \\ \end{diagram}\hspace*{1cm}%
\begin{diagram}[h=2em,w=1.5em] F(U) \otimes \mathbf{1} & \rTo^{ r_{F(U)} } &
F(U) \\ \dTo<{F(U) \otimes\phi_0 } & & \uTo>{ F(r_U) } \\
F(U)\otimes
F(\mathbf{1})& \rTo_{\phi_2(U, \mathbf{1}) } & F(U\otimes \mathbf{1}) \\
\end{diagram}
\end{equation*}%
are commutative, and $\phi _{2}(U,V):F(U)\otimes F(V)\rightarrow
F(U\otimes V)$ is a family of functorial isomorphisms such that
the following diagram
\begin{equation*}
\begin{diagram}[h=2em] (F(U)\otimes F(V))\otimes
F(W)&\rTo^{\phi_2(U,V)\otimes F(W)}&F(U\otimes V)\otimes
F(W)&\rTo^{\phi_2(U\otimes V,W)}&F((U\otimes V)\otimes W)\\
\dTo>{a_{F(U),F(V),F(W)}}&&&&\dTo<{F(a_{ U,V, W})}\\ F(U)\otimes
(F(V)\otimes F(W))&\rTo_{F(U)\otimes \phi_2(V,W)}&F(U)\otimes
F(V\otimes W)& \rTo_{\phi_2(U,V\otimes W)}&{F(U\otimes (V\otimes
W))} \end{diagram}
\end{equation*}%
is commutative. A monoidal functor $(F,\phi _{0},\phi _{2})$ is called \emph{%
strict} if both $\phi _{0}$ and $\phi _{2}$ are identities.
\end{claim}

The following Proposition states that the image of an (co)algebra
through a monoidal functor carries a natural (co)algebra
structure.

\begin{proposition}
\label{pro: algebras through monoidal functors} Let $\mathcal{M}$ and $%
\mathcal{M}^{\prime}$ be monoidal categories. Let $(F,\phi
_{0},\phi _{2}),$
be a monoidal functor between the categories $\mathcal{M}$ and $\mathcal{M}%
^{\prime }$. Then:

1) If $(A,m,u)$ is an algebra in $\mathcal{M}$, then $%
(F(A),m_{F(A)},u_{F(A)})$ is an algebra in $\mathcal{M}^{\prime
}$, where
\begin{equation*}
\begin{tabular}{l}
$m_{F(A)}:=F(A)\otimes F(A)\overset{\phi _{2}(A,A)}{\longrightarrow }%
F(A\otimes A)\overset{F(m)}{\longrightarrow }F(A)$ \\
$u_{F(A)}:=\mathbf{1}^{\prime }\overset{\phi _{0}}{\longrightarrow }F(%
\mathbf{1})\overset{F(u)}{\longrightarrow }F(A)$.%
\end{tabular}%
\end{equation*}%
Moreover $F\left( f\right) :F\left( A\right) \rightarrow F\left(
B\right) $
is an algebra homomorphism in $\mathcal{M}^{\prime }$ whenever $%
f:A\rightarrow B$ is an algebra homomorphism in $\mathcal{M}$. The
converse also holds true if $F$ is faithful.

2) If $(C,\Delta ,\varepsilon )$ is a coalgebra in $\mathcal{M}$, then $%
(F(C),\Delta _{F(C)},\varepsilon _{F(C)})$ is a coalgebra in $\mathcal{M}%
^{\prime }$, where
\begin{equation*}
\begin{tabular}{l}
$\Delta _{F(C)}:=F(C)\overset{F(\Delta )}{\longrightarrow }F(C\otimes C)%
\overset{\phi _{2}^{-1}(C,C)}{\longrightarrow }F(C)\otimes F(C)$ \\
$\varepsilon _{F(C)}:=F(C)\overset{F(\varepsilon )}{\longrightarrow }F(%
\mathbf{1})\overset{\phi _{0}^{-1}}{\longrightarrow }\mathbf{1}^{\prime }.$%
\end{tabular}%
\end{equation*}%
Moreover $F\left( f\right) :F\left( C\right) \rightarrow F\left(
D\right) $
is a coalgebra homomorphism in $\mathcal{M}^{\prime }$ whenever $%
f:C\rightarrow D$ is a coalgebra homomorphism in $\mathcal{M}$.
The converse also holds true if $F$ is faithful.
\end{proposition}

\begin{proof}
follows directly from the definitions.
\end{proof}

Let $K$ be a field. From now on $\mathfrak{M}_{K}$ will denote the
category of vector spaces over $K$.

\begin{theorem}
\label{teo: monoidal functor}Let $(\mathcal{M},\otimes
,\mathbf{1)}$ be a complete and cocomplete abelian monoidal
category. Let $K$ be a field and let
\begin{equation*}
(F,\phi _{0},\phi _{2}):(\mathcal{M},\otimes ,\mathbf{1)}\rightarrow (%
\mathfrak{M}_{K},\otimes _{K},K)
\end{equation*}%
be an additive monoidal functor. Assume that:

\begin{enumerate}
\item $F$ is faithful.

\item $F$ preserves kernels and cokernels.

\item $F$ preserves denumerable products and coproducts.
\end{enumerate}

Let $C$ and $E$ be coalgebras in $\mathcal{M}$. Let
$f_{C}:E\rightarrow C$
be a coalgebra homomorphism and let $f_{M}:E\rightarrow M$ be a morphism of $%
C$-bicomodules, where $E$ is a bicomodule via $f_{C}$. Assume that
$$F\left( f_{M}\right) (\mathrm{Corad}(F\left( E\right) ))=0.$$ Then
there is a unique coalgebra homomorphism $f:E\rightarrow
T_{C}^{c}(M)$ such that $p_{0}f=f_{C}$ and $p_{1}f=f_{M}$, where
$p_{n}:T_{C}^{c}(M)\rightarrow M^{\square n}$ denotes the
canonical projection.
\begin{equation*}
\begin{diagram}[h=2em,w=3em]
T_{C}^{c}(M)&\rTo^{p_1}&M\\\dTo^{p_0}&\luDotsto^{f}&\uTo>{f_M}\\C&%
\lTo_{f_C}&E \end{diagram}
\end{equation*}
\end{theorem}

\begin{proof}
In order to simplify the computations, we will omit the
isomorphisms $\phi _{0},\phi _{2}.$\newline For every object $X\in
\mathcal{M}$, we will denote by $X^{\prime }$ the vector
space $F\left( X\right) .$ Also, for every morphism $f\in \mathcal{M}$%
, we will denote by $f^{\prime }$ the $K$-linear map $F\left(
f\right) .$\newline Recall that $E^{\prime }$ and $D^{\prime }$
carry the coalgebra structures described in Proposition \ref{pro:
algebras through monoidal functors}.\newline Since $F$ preserves
denumerable coproducts and cokernels, it preserves
denumerable direct limits and finite coproducts so that%
\begin{eqnarray*}
F\left[ (T_{C}^{c}(M),(\sigma _{n})_{n\in \mathbb{N}})\right]
&=&F\left[ \underrightarrow{\lim }\left( \left( C^{n}(M)\right)
_{n\in \N},\left( \sigma _{m}^{n}\right) _{m,n\in\N}\right) \right] \\
&=&\underrightarrow{\lim }\left( \left( F\left[ C^{n}(M)\right]
\right) _{n\in
\mathbb{N}
},\left( F\left( \sigma _{m}^{n}\right) \right) _{m,n\in
\mathbb{N}
}\right) \\
&=&\underrightarrow{\lim }\left( \left[ \left( C^{\prime }\right)
^{n}(M^{\prime })\right] _{n\in
\mathbb{N}
},\left( \sigma _{m}^{n}\right) _{m,n\in
\mathbb{N}
}\right) .
\end{eqnarray*}%
\newline
Therefore we get $\left( T_{C}^{c}(M)\right) ^{\prime
}=T_{C^{\prime }}^{c}(M^{\prime }),$ that is the cotensor
coalgebra in the category of vector spaces of the $C^{\prime
}$-bicomodule $M^{\prime }$.\newline Now, by Proposition \ref{pro:
algebras through monoidal functors}, $f_{C}^{\prime }:E^{\prime
}\rightarrow C^{\prime }$ is a coalgebra homomorphism and
$f_{M}^{\prime }:E^{\prime }\rightarrow M^{\prime }$ is a morphism
of $C^{\prime }$-bicomodules, where $E^{\prime }$ is a bicomodule
via $f_{C}^{\prime }$. By hypothesis, $f_{M}^{\prime }(\mathrm{Corad}%
(E^{\prime }))=0.$

By Theorem \ref{teo Nichols}, there is a coalgebra $D$ in
$\mathfrak{M}_{K}$ and a monomorphism $\delta :D\rightarrow
E^{\prime }$ which is a coalgebra
homomorphism in $\mathfrak{M}_{K}$ such that $\underrightarrow{\lim }%
D^{\wedge _{E^{\prime }}^{n}}=E^{\prime }$ and $f_{M}^{\prime }\delta =0.$%
\newline
By Theorem \ref{coro: univ property of cotensor coalgebra}, there
exists a unique coalgebra homomorphism $g:E^{\prime }\rightarrow
T_{C^{\prime
}}^{c}(M^{\prime })$ such that $p_{0}^{\prime }g=f_{C}^{\prime }$ and $%
p_{1}^{\prime }g=f_{M}^{\prime }$, where
$p_{n}:T_{C}^{c}(M)\rightarrow
M^{\square n}$ denotes the canonical projection:%
\begin{equation*}
\begin{diagram}[h=2em,w=3em]
T_{C'}^{c}(M')&\rTo^{p'_1}&M'\\\dTo^{p'_0}&\luDotsto^{g}&\uTo>{f'_M}\\C'&%
\lTo_{f'_C}&E' \end{diagram}
\end{equation*}%
We will prove that $g=f^{\prime }$ for a suitable morphism
$f:E\rightarrow T_{C}^{c}(M)$ in $\mathcal{M}$.\newline By
construction (see (\ref{relation f})):
\begin{equation*}
p_{k}^{\prime }g=\left( f_{M}^{\prime }\right) ^{\square k}\overline{\Delta }%
_{E^{\prime }}^{k-1}\text{ for any }k\in \mathbb{N}.
\end{equation*}%
\newline
Let
\begin{equation*}
q_{k}:\dprod\limits_{n\in\N }M^{\square n}\rightarrow M^{\square
k}
\end{equation*}%
be the canonical projection and let%
\begin{equation*}
\omega :T_{C}^{c}(M)=\dbigoplus\limits_{n\in\N }M^{\square
n}\rightarrow \dprod\limits_{n\in\N }M^{\square n}
\end{equation*}%
be the diagonal morphism of the $\left( p_{k}\right) _{k\in \N }.$
This is uniquely defined by
\begin{equation*}
q_{k}\omega =p_{k},\text{ for every }k\in \N .
\end{equation*}%
Since $F$ preserves kernels and denumerable (co)products, it is
clear that
\begin{equation*}
\omega ^{\prime }:T_{C^{\prime }}^{c}(M^{\prime
})=\dbigoplus\limits_{n\in \N }\left( M^{\prime }\right) ^{\square
n}\rightarrow \dprod\limits_{n\in \N }\left( M^{\prime }\right)
^{\square n}
\end{equation*}%
is a monomorphism. Since $F$ is faithful, we get that $\omega $ is
a monomorphism too. \newline Let
\begin{equation*}
\alpha :E\rightarrow \dprod\limits_{n\in \N }M^{\square n}
\end{equation*}%
be the diagonal morphism of the $\left( f_{M}^{\square k}\circ \overline{%
\Delta }_{E}^{k-1}\right) _{k\in \N
}.$ This is uniquely defined by%
\begin{equation*}
q_{k}\alpha =f_{M}^{\square k}\circ \overline{\Delta
}_{E}^{k-1},\text{ for every }k\in \N .
\end{equation*}%
We have that
\begin{equation*}
\omega ^{\prime }g=\alpha ^{\prime }.
\end{equation*}%
In fact%
\begin{equation*}
q_{k}^{\prime }\omega ^{\prime }g=\left( q_{k}\omega \right)
^{\prime
}g=p_{k}^{\prime }g=\left( f_{M}^{\prime }\right) ^{\square k}\overline{%
\Delta }_{E^{\prime }}^{k-1}=\left( f_{M}^{\square k}\circ \overline{\Delta }%
_{E}^{k-1}\right) ^{\prime }=\left( q_{k}\alpha \right) ^{\prime
}=q_{k}^{\prime }\alpha ^{\prime },
\end{equation*}%
for every $k\in \N .$\newline Let $\left( \mathrm{Coker}\left(
\omega \right) ,\tau \right) $ be the
cokernel of $\omega $ in $\mathcal{M}$. We have%
\begin{equation*}
\left( \tau \alpha \right) ^{\prime }=\tau ^{\prime }\alpha
^{\prime }=\tau ^{\prime }\omega ^{\prime }g=\left( \tau \omega
\right) ^{\prime }g=0.
\end{equation*}%
Since $F$ is faithful, we get $\tau \alpha =0.$ Since $\omega$ is
a monomorphism, it is clear that $(T^c_C(M),\omega)=\K(\tau)$. By
the universal property of
kernels, there exists a unique morphism $f:E\rightarrow T_{C}^{c}(M)$ in $%
\mathcal{M}$ such that%
\begin{equation*}
\omega f=\alpha .
\end{equation*}%
We have%
\begin{equation*}
\omega ^{\prime }f^{\prime }=\left( \omega f\right) ^{\prime
}=\alpha ^{\prime }=\omega ^{\prime }g.
\end{equation*}%
Since $\omega ^{\prime }$ is a monomorphism we get
\begin{equation*}
f^{\prime }=g.
\end{equation*}%
By Proposition \ref{pro: algebras through monoidal functors}, we have that $%
f $ is a coalgebra homomorphism. \newline
It remains to prove that $f$ is uniquely defined by the relations $%
p_{0}f=f_{C}$ and $p_{1}f=f_{M}.$ Let $h$ be another coalgebra
homomorphism such that $p_{0}h=f_{C}$ and $p_{1}h=f_{M}.$ Then
\begin{equation*}
p_{0}^{\prime }h^{\prime }=f_{C}^{\prime }=p_{0}^{\prime }g\text{
\qquad and }\qquad p_{1}^{\prime }h^{\prime }=f_{M}^{\prime
}=p_{1}^{\prime }g.
\end{equation*}%
By uniqueness of $g,$ we get%
\begin{equation*}
h^{\prime }=g=f^{\prime }.
\end{equation*}%
Since $F$ is faithful we conclude that $h=f.$
\end{proof}

\begin{claim}
Let $H$ be a Hopf algebra over a field $K.$ Consider the following
examples of monoidal categories.\smallskip

$\bullet $ \emph{The category $_{H}{\mathfrak{M}}=(_{H}{\mathfrak{M}}%
,\otimes _{K},K)$, of all left modules over $H$}. The tensor
$V\otimes W$ of two left $H$-modules is an object in
$_{H}{\mathfrak{M}}$ via the diagonal
action; the unit is $K$ regarded as a left $H$-module via $\varepsilon _{H}$%
.\smallskip

$\bullet $ \emph{The category ${_{H}\mathfrak{M}_{H}}=({_{H}\mathfrak{M}_{H}}%
,\otimes _{K},K)$, of all two-sided modules over $H$}. The tensor
$V\otimes W $ of two $H$-bimodules carries, on both sides, the
diagonal action; the unit is $K$ regarded as a $H$-bimodule via
$\varepsilon _{H}$.\smallskip

$\bullet $ \emph{The category ${^{H}\mathfrak{M}}=({^{H}\mathfrak{M}}%
,\otimes _{K},K)$, of all left comodules over $H$}. The tensor product $%
V\otimes W$ of two left $H$-comodules is an object in
$^{H}{\mathfrak{M}}$ via the diagonal coaction; the unit is $K$
regarded as a left $H$-comodule via the map $k\mapsto 1_{H}\otimes
k$.\smallskip

$\bullet $ \emph{The category ${^{H}\mathfrak{M}^{H}}=({^{H}\mathfrak{M}^{H}}%
,\otimes _{K},K)$ of all two-sided comodules over $H$}. The tensor
$V\otimes W$ of two $H$-bicomodules carries, on both sides, the
diagonal coaction; the unit is $K$ regarded as a $H$-bicomodule
via the maps $k\mapsto 1_{H}\otimes k$ and $k\mapsto k\otimes
1_{H}$.\smallskip

$\bullet $ \emph{The category ${_{H}^{H}\mathcal{YD}}=({_{H}^{H}\mathcal{YD}}%
,\otimes _{K},K)$ of left Yetter-Drinfeld modules over $H$}.
Recall that an
object $V$ in ${_{H}^{H}\mathcal{YD}}$ is a left $H$-module and a left $H$%
-comodule satisfying, for any $h\in H,v\in V$, the compatibility
condition:
\begin{equation*}
\sum ({}h_{(1)}v)_{<-1>}{h_{(2)}}\otimes (h_{(1)}v)_{<0>}=\sum
h_{(1)}v_{<-1>}\otimes h_{(2)}v_{<0>}
\end{equation*}%
or, equivalently,
\begin{equation*}
\rho (hv)=\sum h_{(1)}v_{<-1>}S(h_{(3)})\otimes h_{(2)}v_{<0>},
\end{equation*}%
where $\Delta _{H}\left( h\right) =\sum h_{(1)}\otimes h_{(2)}$
and $\rho \left( v\right) =\sum v_{<-1>}\otimes v_{<0>}$ denote
the comultiplication of $H$ and the left $H$-comodule structure of
$V$ respectively (we used Sweedler notation).

The tensor product $V\otimes W$ of two Yetter-Drinfeld modules is
an object in ${_{H}^{H}\mathcal{YD}}$ via the diagonal action and
the codiagonal
coaction; the unit in $_{H}^{H}\mathcal{YD}$ is $K$ regarded as a left $H$%
-comodule via the map $x\mapsto 1_{H}\otimes x$ and as a left
$H$-module via the counit $\varepsilon _{H}$.\medskip\newline
Let $\mathcal{M}$ denote one of the categories above and let $F:\mathcal{M}%
\rightarrow \mathfrak{M}_{K}$ be the forgetful functor. Then Theorem \ref%
{teo: monoidal functor} applies.
\end{claim}

\begin{claim}
Let $\left( H,\phi \right) $ be a quasi-Hopf algebra over a field
$K$. Then the category of right $H$-modules
\emph{${\mathfrak{M}}$}$_{H}$ is a monoidal category. By
\cite[Example 9.1.4, page 422]{Maj2}, the forgetful functor
\emph{${\mathfrak{M}}$}$_{H}\rightarrow \mathfrak{M}_{K}$ is
monoidal if and only if $H$ is twisted-equivalent to an ordinary
Hopf algebra. Therefore theorem \ref{teo: monoidal functor} does
not apply in general to \emph{${\mathfrak{M}}$}$_{H}$. Anyway
Theorem \ref{coro: univ property of cotensor coalgebra} still
holds. In fact, since $H$ is in
particular an ordinary $K$-algebra, it is clear that \emph{${\mathfrak{M}}$}$%
_{H}$ is a cocomplete and complete abelian monoidal category
satisfying AB5.
\end{claim}

\noindent\begin{minipage}[t]{15cm}\vspace*{2mm}\sc \footnotesize
A. Ardizzoni.\\ University of Ferrara, Department of
Mathematics, Via Machiavelli 35, Ferrara,  I-44100, Italy .\\
{\small\it email:} {\small \rm
alessandro.ardizzoni@unife.it}\\{\small\it URL:} {\small \rm
http://www.unife.it/utenti/alessandro.ardizzoni}\vspace*{1mm}%

C. Menini.\\ University of Ferrara, Department of Mathematics,
Via Machiavelli 35, Ferrara, I-44100, Italy.\\
{\small\it email:} {\small\ttfamily \rm
men@dns.unife.it}\\{\small\it URL:} {\small \rm
http://web.unife.it/utenti/claudia.menini}\vspace*{1mm}%

D. \c{S}tefan.\\ University of Bucharest, Faculty of
Mathematics, Str. Academiei 14, Bucharest, Ro--70109, Romania.\\
{\small\it email:} {\small\ttfamily \rm dstefan@al.math.unibuc.ro}
\end{minipage}
\end{document}